\documentstyle{amlts}
\begin{document}
\annalsline{155}{2002}
\received{May 25, 2000}
\startingpage{491}
\def\bye{\end{document}}
 \font\tenrm=cmr10

\input boxedeps.tex 
\SetepsfEPSFSpecial 
\HideDisplacementBoxes
\def\figin#1#2{
$$
 {\BoxedEPSF{#1.eps scaled
#2}%
}%
$$
\noindent}

\def\ritem#1{\item[{\rm #1}]}
\def\tag#1{\eqno({#1})}
\def\joinrel{\mathrel{\mkern-4mu}}
\def\relbar{\mathrel{\smash-}}
\def\lrar{\relbar\joinrel\relbar\joinrel\relbar\joinrel\rightarrow}
\def\llar{\leftarrow\joinrel\relbar\joinrel\relbar}
\def\nn{\nonumber}

\catcode`\@=11
\font\twelvemsb=msbm10 scaled 1100
\font\tenmsb=msbm10
\font\ninemsb=msbm10 scaled 800
\newfam\msbfam
\textfont\msbfam=\twelvemsb  \scriptfont\msbfam=\ninemsb
  \scriptscriptfont\msbfam=\ninemsb
\def\msb@{\hexnumber@\msbfam}
\def\Bbb{\relax\ifmmode\let\next\Bbb@\else
 \def\next{\errmessage{Use \string\Bbb\space only in math
mode}}\fi\next}
\def\Bbb@#1{{\Bbb@@{#1}}}
\def\Bbb@@#1{\fam\msbfam#1}
\catcode`\@=12

 \catcode`\@=11
\font\twelveeuf=eufm10 scaled 1100
\font\teneuf=eufm10
\font\nineeuf=eufm7 scaled 1100
\newfam\euffam
\textfont\euffam=\twelveeuf  \scriptfont\euffam=\teneuf
  \scriptscriptfont\euffam=\nineeuf
\def\euf@{\hexnumber@\euffam}
\def\frak{\relax\ifmmode\let\next\frak@\else
 \def\next{\errmessage{Use \string\frak\space only in math
mode}}\fi\next}
\def\frak@#1{{\frak@@{#1}}}
\def\frak@@#1{\fam\euffam#1}
\catcode`\@=12

\input amssym.def
\input amssym.tex

\def\operatorname#1{{\rm #1}}
  \newcommand{\km}{{\frak g}}
  \newcommand{\nskip}{\vspace{6pt}\noindent}
  \newcommand{\secskip}{\vspace{12pt}\noindent}
  \newcommand{\Fr}{\operatorname{Fr}}
  \newcommand{\Hom}{\operatorname{Hom}}
  \newcommand{\Char}{\operatorname{Char}}
  \newcommand{\End}{\operatorname{End}}
  \newcommand{\Lie}{\operatorname{Lie}}
  \newcommand{\sgn}{\operatorname{sgn}}
  \newcommand{\Frp}{\Fr'}
  \newcommand{\Frps}{\Fr^{\prime\ast}}
  \newcommand{\Frb}{\overline{\Fr}}
  \newcommand{\Frs}{\Fr^*}
  \newcommand{\K}{{\Bbb K}}
  \newcommand{\Ocal}{{\cal O}}

  \newcommand{\Fcal}{{\cal F}}  
\newcommand{\Dcal}{{\cal D}} 
  
   \newcommand{\lam}{\lambda}
   \newcommand{\half}{\frac{1}{2}}
   \newcommand{\veps}{\varepsilon}
   \newcommand{\om}{(\omega )}
   \newcommand{\bspm}{\left({\begin{array}{c}}} 
   \newcommand{\espm}{\end{array}\right)}
   \newcommand{\bsbm}{\left[\begin{smallmatrix}} 
   \newcommand{\esbm}{\end{smallmatrix}\right]}
   \newcommand{\bpm}{\left( \begin{array}{c}}
 \newcommand{\epm}{\end{array}\right)}
   \newcommand{\bbm}{\left[ \begin{array}{c}} 
\newcommand{\ebm}{\end{array}\right]}
                                           
   \renewcommand{\bf}{\bfseries}
   \renewcommand{\it}{\itshape}
   \newcommand{\nzsum}[1]{\sum_{n=0}^{#1}}
   \newcommand{\nosum}[1]{\sum_{n=1}^{#1}}
 \newcommand{\Fgil}{{\Bbb F}}
   \newcommand{\Cgil}{{\Bbb C}}   \newcommand{\Zgil}{{\Bbb Z}}
   \newcommand{\Pgil}{{\Bbb P}} 
   \newcommand{\Ngil}{{\Bbb N}}
 \newcommand{\zi}{{\Bbb Z}_{\xi}}  \newcommand{\N}{{\Bbb N}}
   \newcommand{\Qgil}{{\Bbb Q}}
   \newcommand{\Zt}{\tilde{{\Bbb Z}}}

 \newcommand{\ufk}{{\frak u}}    \newcommand{\Lfk}{{\frak L}}
   \newcommand{\hfk}{{\frak h}}  \newcommand{\gfk}{{\frak g}}
   \newcommand{\bfk}{{\frak b}}  \newcommand{\Ufk}{{\frak U}}
   \newcommand{\nfk}{{\frak n}}
   \newcommand{\pfk}{{\frak p}}   \newcommand{\dfk}{{\frak d}}
   \newcommand{\wfk}{{\frak w}}
\newcommand{\Xfk}{{\frak X}}
   \newcommand{\Acal}{{\cal A}}   \newcommand{\Bcal}{{\cal B}}
   \newcommand{\Ccal}{{\cal C}}   \newcommand{\Lcal}{{\cal L}}
             \newcommand{\Cbar}{ \bar{{\cal C}} }
 \newcommand{\Ucal}{{\cal U}}  \newcommand{\ucal}{\mathsf{u}}
   \newcommand{\Ebar}{\bar{E}} 
 \newcommand{\Nbar}{\bar{N}} 
  \newcommand{\Fbar}{\bar{F}}
   \newcommand{\Ubar}{\bar{U}}   \newcommand{\Qbar}{\bar{Q}}
   \newcommand{\Dbar}{\bar{D}}   \newcommand{\Mbar}{\bar{M}}
\newcommand{\Xbar}{\bar{X}}
   \newcommand{\Btimes}{\operatorname*{\otimes}_{\Bcal}}
  \newcommand{\pch}{\psi_{\gamma}}  
\newcommand{\mapleft}[1]{\smash{\mathop{\longleftarrow}\limits^{#1}}}
 \newcommand{\tX}{\tilde{X}}
\newcommand{\Ical}{{\cal I}}
 
  \newcommand{\Frd}{\Fr_{\Delta}}
  \newcommand{\Frds}{\Fr_{\Delta}^{\ast}}
  \newcommand{\Frdps}{\Fr_{\Delta}^{\prime\ast}}
  
  \newcommand{\simarr}{\stackrel{\sim}{\rightarrow}}

  \newcommand{\Zfk}{{\frak Z}}
  \newcommand{\Zbar}{\bar{Z}}

 \newcommand{\Zcal}{{\cal Z}}
 
  \newcommand{\vfk}{{\frak v}}

\def\ssni#1{\smallbreak\noindent{#1}. }

\def\sni#1{\smallbreak\noindent{\phantom{0}#1}. }

\title{Algebraization of Frobenius splitting\\ via quantum
groups}
\shorttitle{Algebraization of Frobenius splitting}  
\acknowledgement{The first author was 
supported by NSF grants and the second author by the EEC program TMR ERB 
FMRX--CT97--0100, Algebraic Lie Representations.}

 \twoauthors{Shrawan Kumar}{Peter Littelmann}
 \institutions{University of North Carolina at Chapel Hill, Chapel Hill, NC 27599-3250, USA\\
{\eightpoint {\it E-mail address\/}: shrawan@email.unc.edu}\\
\vglue6pt
Universit\"at Wuppertal, 42097 Wuppertal, Germany\\
{\eightpoint {\it E-mail address\/}: littelmann@math.uni-wuppertal.de
}}

\vglue24pt \centerline{\bf Abstract}
\vglue16pt

An important breakthrough in understanding the geometry of Schubert
varieties
was the introduction of the notion of  Frobenius split varieties and the result that the flag varieties $G/P$ are  Frobenius split.
The aim of this article is to give in this case a  complete and self contained
representation theoretic approach to this method. The geometric Frobenius
method (in char  $k=p>0$) will here be  replaced by Lusztig's Frobenius
maps for
quantum groups at roots of unity  (which exist not only for primes but any odd
integer $\ell > 1$).

\vglue24pt \centerline{\bf Table of Contents}
\vglue8pt
\sni{0} Introduction
\sni{1} Notation, preliminaries and review of certain results of Lusztig and\hfill\break
\noindent \hglue24pt  Andersen-Polo-Wen
\sni{2} Definition of   quantized Frobenius homomorphism
\sni{3} Definition of quantized Frobenius splitting
\sni{4} Stronger quantized Frobenius splitting
\sni{5} The Kempf vanishing theorem
\sni{6} Sheafification: Frobenius splitting of $G/B$ and Schubert varieties
\sni{7} Splitting of the diagonal in $G/B\times G/B$
\sni{8} Frobenius splitting of quantized Bott-Samelson desingularization
\sni{9} Extension of results to the parabolic case
\ssni{10} Appendix: Applications

\advance\sectioncount by -1 
\section{Introduction}

The passage from representations of quantum groups at roots of unity in 
characteristic zero to 
representations of algebraic groups in characteristic $p$ is extremely 
important in view of Lusztig's conjectures. 
(Recall that Andersen-Jantzen-Soergel 
confirmed the conjectured link in Lusztig's program but only asymptotically.) 
The aim of the present article is to establish that Lusztig's two Frobenius 
maps in characteristic zero lead naturally to two familiar objects in 
characteristic $p$. One Frobenius map leads simply to the Frobenius map in 
characteristic $p$, while the other leads to the so-called canonical Frobenius
splitting on $G/B$ and related varieties.

Let $k$ be  an algebraically closed field of arbitrary characteristic and let 
 $G$ be a semisimple simply-connected algebraic group over $k$ with 
 a fixed  Borel subgroup  $B$ and the associated Weyl group $W$. 
Let $X=G/B$ be the flag variety and let (for any $w \in W$)  $X(w) \subset X$ 
be
 the Schubert subvariety,
which is  the closure of the
$B$-orbit $BwB$ in $G/B$. For a homogeneous  line
bundle $\Lcal$ on $X$, the cohomology groups
$H^i(X,\Lcal)$ are $G$-modules and the corresponding
groups $H^i(X(w),\Lcal\vert_{X(w)})$ inherit  naturally the structure
of  $B$-modules. These modules have been extensively studied
from  algebro-geometric as well as representation theoretic points of view.

An important breakthrough in understanding  the geometry of\break Schubert
varieties
was the introduction, by  Mehta-Ramanathan and Ramanan-\break Ramanathan,
 of the notion of a Frobenius $\Dcal$--split variety $X$
(defined over $k$ of char~$p > 0$) and
compatibly split subvarieties, where  $\Dcal$ is a line bundle together with a 
nonzero section   
(cf.\ Definition 6.1). 
 `Very few' projective 
varieties turn out to be split but those which do have rather remarkable 
geometric and cohomological properties. The most important class of 
examples of varieties which are Frobenius split arise in group theory. In 
particular, the flag varieties $X = G/B$ are  Frobenius split (in fact are 
$\Dcal$-split
for   the homogeneous line bundle $\Dcal$ corresponding to the character 
$-2 (p-1) \rho$, where $\rho$ is half the sum of positive roots) compatibly 
splitting all the Schubert subvarieties and so are the product varieties 
$X \times X$ compatibly splitting all the $G$-Schu\-bert subvarieties. This 
leads to various important 
geometric facts about them\break (normality, projective normality,  Cohen-Macaulay,
 projective  Cohen-Macaulay, rational singularity etc.) and various  
representation theoretic
results (vanishing theorems, Demazure character formula, good filtrations, 
etc.) (see the papers  
 [MR1],  [RR], [R1], [R2], [MR2], [M] etc.). However, this geometric 
method   does not provide 
 explicit representation  theoretic information directly.

The aim of this article is to give a
complete and self contained representation theoretic approach to these
methods.
The algebro--geometric ``Frobenius methods'' will here be replaced by
Lusztig's two Frobenius maps for quantized enveloping algebras  at roots of 
unity. 

Let $\km$ be a  complex semisimple Lie algebra
with triangular decomposition $\gfk=\nfk\oplus\hfk\oplus\nfk^-$ and
denote by $\bfk=\hfk\oplus\nfk$ the corresponding Borel subalgebra. Assume 
that  $(\gfk , \bfk )$ corresponds to the pair $(G,B)$. 
Fix an odd integer $\ell>1$ which is assumed to be coprime to $3$ if $\gfk$
has a  component of type $G_2$. Choose a primitive $\ell^{\rm th}$
root of unity $\xi$ and let $\zi$ be the corresponding ring of
cyclotomic integers. Let  $\Ufk_{\zi}(\gfk )$
denote the corresponding quantized enveloping algebra over $\zi$ 
obtained by the base 
change $ \Zgil 
[v, v^{-1}] \to \zi, v \mapsto \xi$,  from  Lusztig's  $ \Zgil 
[v, v^{-1}]$-form of the  quantized enveloping algebra $U_{\Qgil(v)} (\gfk)$
(divided by the ideal generated by the
central elements $K_i^\ell-1$)  and let 
$\bar{U}_{\zi}(\gfk )$ be the corresponding classical
 enveloping algebra  over $\zi$ (obtained via the base change 
$\Zgil \hookrightarrow \zi $ from Kostant's $\Zgil$-form of the classical 
enveloping 
algebra 
$\bar{U} (\gfk)$ over $\Cgil$). The subalgebras   $\Ufk_{\zi}(\bfk ), 
\Ufk_{\zi}(\nfk ), \Ufk_{\zi}(\nfk^- )$ (resp.\ $\bar{U}_{\zi}(\bfk ), \bar{U}_{\zi}(\nfk ), \bar{U}_{\zi}(\nfk^- )$) 
are defined similarly. 

Lusztig [Lu2] has defined Frobenius homomorphisms
$\Fr : \Ufk_{\zi}(\gfk ) \to \bar{U}_{\zi}(\gfk )$ which maps
  the generators  by dividing  the exponents by
$\ell$ when possible  (cf.\ Theorem 1.1), and
$\Fr' : \bar{U}_{\zi}(\nfk ) \to \Ufk_{\zi}(\nfk )$ which maps the generators
 by multiplying the exponents by $\ell$. In fact, $\Fr'$ extends
to a homomorphism (still denoted by) $\Fr' : \bar{U}_{\zi}(\bfk ) \to
\Ufk_{\zi}(\bfk )$ (cf.\ Theorem 1.2). We make crucial use of
these maps $\Fr$ and $\Fr'$ together with the homological machinery developed 
by Andersen-Polo-Wen [APW] to define 
natural functors (for any  $\bar{U}_{\zi}(\bfk)$-module $\Mbar$): 
$$
\Fr^* : H^i \bigl(\bar{U}(\gfk )/\bar{U}(\bfk ), \Mbar\bigr)^{ \Fr}
\longrightarrow H^i\bigl(\Ufk (\gfk )/\Ufk (\bfk ), \Mbar^{\Fr}\bigr)
$$
of $\Ufk (\gfk )$-modules (cf.\ Theorem 2.3), and (for any  $\Ufk (\bfk )$-module $M$),
$$
\Frps : H^i\bigl(\Ufk (\gfk )/\Ufk (\bfk ),M\bigr)^{\Fr'} \rightarrow
H^i\bigl(\Ubar(\gfk )/\Ubar (\bfk ), M^{\Fr'}\bigr) 
$$
of $\bar{U}_{\zi}(\bfk)$-modules (cf.\ Theorem 3.8), where 
we have abbreviated  $\Ufk_{\zi}(\gfk )$ by  $\Ufk(\gfk )$ etc. and $\Mbar^{\Fr}$ is a 
$\Ufk (\bfk )$-module under $\Fr$ and the superscript $\Fr'$ has a similar meaning.
Moreover, the composition of these two maps $ \Frps \circ \Fr^*$ is the 
identity map (cf.\ Corollary 3.9). The first map is our representation theoretic replacement of the Frobenius morphism $F$ (which corresponds to the $p^{\rm th}$
power map  on functions) and the second corresponds to a splitting map. (For
one dimensional representations, these maps are given more explicitly in 
[KL] for $i=0$.)

To define a  representation theoretic analogue of the $\Dcal$-splitting in 
our setting,
consider the element $F_o\in\Ufk(\nfk^-)$ defined as the product of
all divided $(\ell-1)$-powers $F_\beta^{(\ell-1)}$ of Lusztig's root vectors,
where the ordering has been chosen relative to a reduced decomposition of
the longest element in the Weyl group. The idea is then to ``twist''
the splitting $\Frps$ by $F_o$ to get, for any $\Ubar(\bfk)$-module $\Mbar$, 
a functorial $\Ubar (\bfk )$-module map
$$
\Frps_{\gamma} : H^i\bigl(\Ufk (\gfk )/\Ufk (\bfk ) ,
\chi^{\xi}_{\gamma}
\otimes\Mbar^{\Fr}\big)^{^{\Fr'}} \rightarrow
H^i(\Ubar(\gfk )/\Ubar (\bfk ) , \Mbar ), 
$$
where $\chi^{\xi}_{\gamma}$ stands for the one-dimensional
$\Ufk(\bfk)$-representation
of weight $\gamma=-2(\ell-1)\rho$ (cf.\ Theorem 4.7). Moreover,
 both the maps $\Frps$ and $\Frps_\gamma$  also commute  with the action (induced by $\Fr'$) of
$\Ubar(\nfk^-)$. Further, all the above three maps are compatible with any base change.  (Note that though $\Ubar(\bfk)$ and
$\Ubar(\nfk^-)$ both 
act on $H^i\bigl(\Ufk (\gfk )/\Ufk (\bfk ) ,
M\big)$ and similarly on $H^i\bigl(\Ufk (\gfk )/\Ufk (\bfk ) ,\break \chi^{\xi}_{\gamma}
\otimes\Mbar^{\Fr}\big)$ via $\Fr'$, these actions do not in general
 glue together to provide a $\Ubar(\gfk )$-action.)

Now assume  that $\ell = p$ is a prime and $k$ is an algebraically 
closed 
field of characteristic $p$.   Since the
constructions
of $\Fr^*,\Frps$ and $\Frps_{\gamma}$  are compatible with any base
change, we consider them under the base change $\zi \to k$ taking $\xi 
\mapsto 1$. Recall  that, for any $\Ufk (\bfk )$-module $M$, there is 
a canonical isomorphism (cf.\ Proposition 5.1),
$$H^i(\Ufk_k(\gfk )/\Ufk_k(\bfk ), M_k) \simeq H^i(G/B, \Lcal (M)),$$
 and similarly, for a  $\Ubar (\bfk )$-module $\Mbar$, there is a canonical
isomorphism
$$H^i(\Ubar_k(\gfk )/\Ubar_k(\bfk ), \Mbar_k) \simeq H^i(G/B, \Lcal
(\Mbar)), $$
 where $M_k := M \otimes_{\zi} k$ etc.  and $\Lcal (M)$ denotes the homogeneous
vector bundle on the flag variety $G/B$ associated to the $\Ufk_k(\bfk)$
(and hence $B$) module $M_k$. 
 Using these identifications and the usual Serre vanishing theorem, one 
readily  deduces from the above functors $\Fr^*$ and $\Frps_\gamma$
 the standard Kempf vanishing theorem  asserting that for any weight $\lam$ 
such that $\lam + \rho$ is dominant, $H^i(G/B, \Lcal (-\lam)) = 0$ for all $i 
> 0$ 
(cf.\ Theorem 5.2).

In  Section~6  we establish a precise connection between our representation
theoretic approach with the 
algebro-geometric Frobenius splitting  mentioned earlier. 
Actually, by an appropriate  `sheafification'   we obtain 
from the functors $\Frps$ and $\Frps_\gamma$ an entirely new proof 
(purely from the  representation theory of quantum groups) of the Frobenius 
splitting and the stronger Frobenius $\Dcal = 
\Lcal (\gamma)$-splitting respectively   of the flag variety $G/B$, and these
compatibly split all the Schubert subvarieties $X_w$ (cf.\ Theorems 6.4, 6.5
and 6.7). In fact, from our constructions, it is immediately clear that the 
splitting
thus obtained is canonical (in the sense of Mathieu), a property which is
not so transparent (though true) from the original (geometric) proof of the 
splitting of $G/B$ given by  Mehta-Ramanathan. In particular, from the uniqueness (noted by Mathieu) of the canonical splitting on 
$G/B$, we deduce that our splitting coincides with 
the original splitting given in  [MR1]. 

Since our 
constructions  live at the level of algebras over $\zi$, we can view them
as  `characteristic zero lift' of the (characteristic $p$) Frobenius splitting and Frobenius $\Dcal$-splitting of $G/B$. Also, they are
defined for any odd integer 
$\ell > 1$ which is   coprime to $3$ if $\gfk$ has a 
 component of type $G_2$ (not only for primes). It is possible that the 
above restriction on $\ell$ can be removed by using certain results of Kaneda
and Andersen-Paradowski (cf.\ Remark A.8). 

We extend the above constructions and results to cover the case of the product
flag variety $G/B\times G/B$ and deduce, by methods as above, that\break
 $G/B\times G/B$ is Frobenius split such that all the $G$-Schubert 
subvarieties (in particular the diagonal) are compatibly split (cf.\ 
Theorem 7.5). 

We also deduce the splitting of the Bott-Samelson desingularization 
from an analogous quantum setup (cf.\ Section 8).

Even though in a large part of the paper, for notational convenience,  we considered the case of the (full) flag 
variety $G/B$, most of the results can easily be generalised to cover the 
case of $G/P$ for any parabolic subgroup $P$. We formulate  the extensions
in Section~9 but omit the proofs as they are similar to that of $G/B$.

For completeness and convenience of the reader, we collect various important
(and by now standard) 
consequences of the above Frobenius splitting results in the appendix. This 
includes normality,  the Demazure character formula, projective 
normality, Cohen-Macaulay, projective 
Cohen-Macaulay, and rational singularity of Schubert varieties in $G/P$. 
In particular, with our setup, these results follow from the representation 
theory of quantum groups  and
 the Serre vanishing theorem. There are other algebraic proofs of the Demazure character formula
using quantum groups:  by Kashiwara using  his crystal base  [Kas] and by
Littelmann using his LS path model [Li]. 

We believe that many other results (concerning the Frobenius splitting property 
of varieties arising in group theory) are amenable to the methods of this 
paper: e.g., we believe that one can deduce  the `good filtration property'
originally due to Donkin in most cases (and proved by Mathieu in general).

This paper is a sequel to our paper [KL], where a weaker form of some of the 
results of this paper are proved. However, we have kept the exposition 
of this paper almost self-contained (with the exception of [KL, Lemma 3 and 
Th.~1],
which we use here without including a proof). 

We thank H.\ H.\ Andersen for some helpful correspondences, and D. Prasad and  P. Polo
for a helpful conversation.  We are grateful to the referee for the comments,
in particular, for pointing out an inaccuracy.

    \section{Notation, preliminaries and review of certain  results\\ of 
Lusztig and
Andersen-Polo-Wen}

Let us fix a Cartan matrix of finite type $A = (a_{i,j})_{1\leq
i,j\leq n}$.  Then there is a diagonal matrix $D$ with
positive integral diagonal entries $(d_1,\cdots ,d_n)$ with  
$d_i\in\{ 1,2,3\}$ such that
g.c.d. $(d_1,\cdots ,d_n)=1$ and $DA$ is symmetric and positive
definite. 

Let $\gfk = \gfk (A)$ be the semisimple Lie algebra over $\Cgil$
associated to the Cartan matrix $A$.  Recall that $\gfk$ is generated
by its Cartan subalgebra $\hfk$ and  positive root vectors $\{
\bar{E}_1,\cdots ,\bar{E}_n\}$ and negative root vectors
$\{\bar{F}_1,\cdots ,\bar{F}_n\}$ subject to certain relations.  Let
$\bfk ; \nfk ; \bfk^- ; \nfk^-$ be the Lie subalgebras of $\gfk$
generated respectively by $\{\hfk , \bar{E}_1,\cdots ,\bar{E}_n\}$;
$\{\bar{E}_1,\cdots ,\bar{E}_n\}$; $\{\hfk ,\bar{F}_1,\cdots
,\bar{F}_n\}$; $\{\bar{F}_1,\cdots ,\bar{F}_n\}$.  Let
$\bar{U}_{\Zgil}(\gfk )$ be the Kostant $\Zgil$-form of the
enveloping algebra $\bar{U}(\gfk )$, i.e., the\break $\Zgil$-subalgebra of
$\bar{U}(\gfk )$ generated by $\{\bar{E}_i^{(m)}, \bar{F}_i^{(m)}, 
{h_i\choose m} ; 1\leq i\leq n, m\in\Zgil_+\}$, where
$\bar{E}_i^{(m)} := \frac{\bar{E}^m_i}{m!}$,  $h_i := [\bar{E}_i,
\bar{F}_i]$ and ${h_i\choose m} := \frac{h_i(h_i-1) \cdots
(h_i-m+1)}{m!}$. 
Let $\bar{U}_{\Zgil}(\bfk )$, $\bar{U}_{\Zgil}(\nfk )$ be the
subalgebras of $\bar{U}_{\Zgil}(\gfk )$ generated by $\big\{^{\phantom{|}}\hskip-4pt
\bar{E}_i^{(m)}, {h_i\choose m}   ; 1\leq i\leq n,
m\in\Zgil_+\big\}$ and $\big\{\bar{E}_i^{(m)}; 1\leq i\leq n,
m\in\Zgil_+\big\}$ respectively.  Similarly, let $\bar{U}_{\Zgil}^0$
be the subalgebra of $\bar{U}_{\Zgil}(\gfk )$ generated by $\big\{^{\phantom{|}}\hskip-4pt
{h_i\choose m}   ; 1\leq i\leq n, m\in \Zgil_+\big\}$.  For any
$\Zgil$-algebra $\Bcal$, by $\bar{U}_{\Bcal}(\gfk )$ we mean
  $$
\bar{U}_{\Bcal}(\gfk ) = \bar{U}_{\Zgil}(\gfk )\otimes_{\Zgil}\Bcal ,
  $$
with a similar meaning for $\bar{U}_{\Bcal}(\bfk )$ etc. 

Now we come to the corresponding quantized algebras.  Let $\Acal :=
\Zgil [v, v^{-1}]$, where $v$ is an indeterminate and let $\Qgil (v)$
be
the quotient field of $\Acal$.  Let $U_{\Qgil (v)}(\gfk )$ be the
quantized enveloping algebra over the field $\Qgil (v)$ generated by
$\big\{ E_i, F_i, K_i^{\pm 1}; 1\leq i\leq n\big\}$ and subject to
certain relations (cf.\ [Lu2, \S1.1]).  
Let $U_{\Acal}(\gfk )$ be the $\Acal$-form of
$U_{\Qgil (v)}(\gfk )$ defined by Lusztig, i.e., $U_{\Acal}(\gfk )$
is the
$\Acal$-subalgebra of $U_{\Qgil (v)}(\gfk )$ generated by $\big\{
E_i^{(m)}, F_i^{(m)}, K_i^{\pm 1};1\le i\le n, m\in\Zgil_+\big\}$, where 
  $$
E_i^{(m)} := \frac{E^m_i}{[m]!_{d_i}} ,\; [m]!_{d_i} :=
\prod_{h=1}^m \frac{v^{d_ih}- v^{-d_ih}}{v^{d_i} - v^{-d_i}} \in
\Acal .
  $$
Let $U_{\Acal}(\nfk )$, $U^0_{\Acal}$ be the $\Acal$-subalgebras of
$U_{\Acal}(\gfk )$ generated by $\big\{ E_i^{(m)}; 1\leq i\leq~n,\break
m\in\Zgil_+\big\}$ and $\big\{ K_i^{\pm 1},   {\ninepoint \left[\!\!\begin{array}{c} {K_i; c}\\  
{m}\end{array}\!\! \right]}; 1\leq i\leq n, m\in\Zgil_+, c\in\Zgil\big\}$ respectively, where
  $$
\bbm K_i ; c\\ m\ebm := \prod^m_{h=1} \frac{ K_iv^{d_i(c-h+1)} -
K_i^{-1} v^{-d_i(c-h+1)} }{ v^{d_ih} - v^{-d_ih} } .
  $$

\noindent Also let $U_{\Acal}(\bfk )$ be the subalgebra of $U_{\Acal}(\gfk )$
generated by $U^0_{\Acal}$ and $U_{\Acal}(\nfk )$.  We similarly
define $U_{\Acal}(\nfk^-)$ and $U_{\Acal}(\bfk^-)$. Then $U_{\Qgil (v)}(\gfk )$
is a Hopf algebra with comultiplication $\Delta$, antipode 
(an anti-automorphism) $S$ and co-unit $\epsilon$ given by:
$$ \Delta E_i =  E_i \otimes 1 +  K_i  \otimes E_i,\,\, \Delta F_i =  
F_i \otimes K_i^{-1} + 1 \otimes  F_i, \,\,\Delta K_i =  
K_i \otimes  K_i \,, $$
$$S E_i =  - K_i^{-1} E_i, \,\, S F_i = - F_i K_i, \,\, S K_i =  K_i^{-1} \,,$$
$$\epsilon E_i = \epsilon F_i = 0, \,\, \epsilon K_i = 1 \,.$$
Then $\Delta$ and $S$ keep  $U_{\Acal}(\gfk ), U_{\Acal}(\bfk ), 
U_{\Acal}(\bfk^-)$ and $
U^0_{\Acal}$ stable. \pagebreak

{\it We now fix an odd integer $\ell >1$ which is{\rm ,} in addition{\rm ,}  assumed to be
coprime to $3$ if $G_2$ is a component of $\gfk = \gfk (A)$.}  This
will be our tacit restriction on $\ell$.  Now choose a primitive
$\ell^{\rm th}$ root of unity $\xi$ and let $\zi$ be the corresponding ring
of cyclotomic integers with quotient field $\Qgil_{\xi}$; i.e.,
$\Qgil_{\xi}$
(resp.\ $\zi$) is obtained from $\Qgil$ (resp.\ $\Zgil$) by attaching
$\xi$.

For any $\Acal$-algebra $\Bcal$, by $U_{\Bcal}(\gfk )$ we mean
  $$
U_{\Bcal}(\gfk ) = U_{\Acal}(\gfk )\otimes_{\Acal} \Bcal ,
  $$
and a similar meaning for $U_{\Bcal}(\bfk )$ etc.  In particular,
taking $\Bcal = \zi$ with the homomorphism $\Acal \to\Bcal$,
$v\mapsto \xi$, we get $U_{\zi}(\gfk )$ etc.

We recall the following result due to Lusztig [Lu2, Cor.\ 8.14].

  \proclaim{Theorem}  There is a unique $\zi$\/{\rm -}\/algebra homomorphism
  $$
\Fr : U_{\zi}(\gfk ) \to \bar{U}_{\zi}(\gfk ),
  $$
taking $E_i^{(m)} \mapsto \bar{E}_i^{(m/\ell )}${\rm ,} $F_i^{(m)} \mapsto
\bar{F}_i^{(m/\ell )}${\rm ,} $K_i^{\pm 1}\mapsto 1${\rm ;} where
$\bar{E}_i^{(h)}$ and $ \bar{F}_i^{(h)}$ mean $0$ if $h\notin \Zgil_+$ and{\rm ,}
moreover{\rm ,}  $\Zgil
\hookrightarrow \zi$ is the canonical inclusion{\rm .}  

Then $\Fr$ takes
${K_i; 0 \choose m} \mapsto {h_i\choose 
m/\ell}$ if $\ell$ divides $m${\rm ,} and $0$ otherwise{\rm .}\endproclaim

We also have the following theorem [KL, Lemma 3].  The corresponding result 
with
$\bfk$ replaced by $\nfk$ was proved by Lusztig [Lu2, Lemma 8.6].

  \proclaim{Theorem}  There is a unique $\zi$\/{\rm -}\/algebra homomorphism
  $$
\Fr'=\Fr'_{\bfk} : \bar{U}_{\zi}(\bfk ) \to U_{\zi}(\bfk )/ \langle K_i^{\ell}-1;
1\leq i\leq n\rangle \,,
  $$
taking $\bar{E}_i^{(m)} \mapsto E_i^{(m\ell )}${\rm ,} ${h_i\choose m}
\mapsto {K_i; 0\choose \ell m}$; where $\langle\; \rangle$
denotes the (two sided) ideal of $U_{\zi}(\bfk )$ generated by the
central elements $\{ K_i^{\ell}-1\}$.

Similarly{\rm ,} we have the ${\Bbb Z}_{\xi}$-algebra homomorphism
 $$
\Fr'_{\bfk^-} : \bar{U}_{\zi}(\bfk^- ) 
\to U_{\zi}(\bfk^- )/ \langle K_i^{\ell}-1;
1\leq i\leq n\rangle \,.
  $$
\endproclaim

Let $X:=\{ \lam\in\hfk^* : \lam (h_i)\in\Zgil \mbox{ for all }1\leq
i\leq n\}$ be the set of weights and $X^+ := \{ \lam\in\hfk^* : \lam
(h_i)\in\Zgil_+ \mbox{ for all }1\leq i\leq n\}$ be the set of dominant
weights.  For any $\lam\in X$, define the character $\chi_{\lam} :
U^0_{\Acal}\to\Acal$ by
  \begin{eqnarray*}
  \chi_{\lam}\bigl( K_i^{\pm 1}\bigr) &= &v^{\pm d_i\lam
(h_i)}, \qquad\mbox{and} \\[5pt]
\chi_{\lam}\bbm K_i; c\\m\ebm &=& \bbm \lam (h_i)+c\\ m\ebm_{d_i} ,
  \end{eqnarray*}
where, for $m\in \Zgil_+$ and $c\in\Zgil$ ,
  $$
\bbm c\\ m\ebm_{d_i} := \prod^m_{h=1} \frac{
v^{d_i(c-h+1)}-v^{-d_i(c-h+1)} }{ v^{d_ih}-v^{-d_ih} } \in \Acal .
  $$  

In particular, $\chi_{\lam}$ gives rise to a homomorphism by extending the scalars (denoted
by)
  $$
\chi^{\xi}_{\lam} : U^0_{\zi} \to \zi .
  $$

Moreover, $\chi^{\xi}_{\lam}$ descends to give a homomorphism (again
denoted by)
  $$
\chi^{\xi}_{\lam} : U^0_{\zi}/\langle K_i^{\ell} -1; 1\leq i\leq
n\rangle \to \zi .
  $$
Similarly, there is a homomorphism
  \begin{eqnarray*}
\bar{\chi}_{\lam} : \bar{U}^0_{\Zgil} &\rightarrow \Zgil ,
\qquad\mbox{taking} \\
\bpm h_i\\ m\epm &\mapsto \bpm \lam (h_i)\\ m\epm ,
  \end{eqnarray*}
where for $h\in\Zgil$, $m\in\Zgil_+$, ${h\choose m}$ is the
ordinary binomial coefficient
  $$
\bpm h\\ m\epm := \frac{ h(h-1) \cdots (h-m+1) }{m!} .
  $$
By extending the scalars $\Zgil\hookrightarrow\zi$, we get a
$\zi$-algebra homomorphism (still denoted by)
  $$
\bar{\chi}_{\lam} : \bar{U}^0_{\zi} \to \zi .
  $$

Let us denote, by the corresponding Gothic letter,
  $$ 
\Ufk_{\zi}(\gfk )  := U_{\zi}(\gfk )/\langle K_i^{\ell}-1; 1\leq
i\leq n\rangle
  $$
with a similar meaning for $\Ufk_{\zi}(\bfk )$ and $\Ufk^0_{\zi}$.  By
its definition,   Fr descends to a homomorphism (again denoted
by)
  $$
\Fr : \Ufk_{\zi}(\gfk )\to \bar{U}_{\zi}(\gfk ).
  $$

  \numbereddemo{Definition}  A $\Ufk^0_{\zi}$-module $M$ is called a {\it weight module} if 
  \begin{eqnarray*}
M &= &\oplus_{\lam\in X} \,M_{\lam}\,,\qquad\mbox{where} \\
M_{\lam} &:=& \{ v\in M: av = \chi^{\xi}_{\lam}(a)v, \mbox{ for all
}a\in\Ufk^0_{\zi}\} .
  \end{eqnarray*}

Let $\Ccal_{\zi}(\bfk )$ (resp.\ $\Ccal_{\zi}(\gfk )$) be the category
of those $\Ufk_{\zi}(\bfk )$ (resp.\ $\Ufk_{\zi}(\gfk )$) -modules
$M$ such that
$M = F_{\bfk}(M)$ (resp.\ $M = F_{\gfk}(M)$), where
  \begin{eqnarray*}
F_{\bfk}(M) &:=& \biggl\{ v\in\oplus_{\lam\in X} \,M_{\lam}:
E_i^{(m)}v=0,\, \mbox{for}\,
m\geq m(v), 
\mbox{for some }m(v)\in\Zgil_+ \biggr\} ,\\
F_{\gfk}(M) &:=& \biggl\{ v\in\oplus_{\lam\in X} \,M_{\lam} :
E_i^{(m)}v=F_i^{(m)}v=0, \, \mbox{for}\,
 m\geq m(v) \biggr\} .
  \end{eqnarray*}
In particular, any $M\in {\cal C}_{\Bbb Z_\xi}(\frak b)$ is a weight module under $\Ufk^0_{\zi}$.

We similarly define the notion of weight modules for $\bar{U}^0_{\Zgil}$
(resp.\ $\bar{U}^0_{\zi}$ ) and define the categories
$\Cbar_{\Zgil}(\bfk )$ and $\Cbar_{\Zgil}(\gfk )$ 
(resp.\ $\Cbar_{\zi}(\bfk )$ and $\Cbar_{\zi}(\gfk )$).  Then these are
abelian categories (see [APW, $\S$2.2]).
  \enddemo

  \advance\theoremcount by 1
 (1.4) {\it The induction functor}.
Following [APW, $\S$2], define the induction functor
$H^0\bigl(\Ufk_{\zi}(\gfk )/\Ufk_{\zi}(\bfk ), -\bigr) :
\Ccal_{\zi}(\bfk
)\to \Ccal_{\zi}(\gfk )$ as follows: $\phantom{\sum^\int}$

Take $M\in \Ccal_{\zi}(\bfk )$.  Consider the space $\tilde{M} :=
\Hom_{\Ufk_{\zi}(\bfk )}\bigl(\Ufk_{\zi}(\gfk ), M\bigr)$ of
$\Ufk_{\zi}(\bfk )$-module maps from $\Ufk_{\zi}(\gfk )$ to $M$,
where $\Ufk_{\zi}(\bfk )$ acts on $\Ufk_{\zi}(\gfk )$ via   left
multiplication.  Then $\tilde{M}$ is a $\Ufk_{\zi}(\gfk )$-module
under $(a\cdot f)(b) = f(ba)$, for $a,b\in\Ufk_{\zi}(\gfk )$.  Now
set
  $$
H^0\bigl(\Ufk_{\zi}(\gfk )/\Ufk_{\zi}(\bfk ), M\bigr) :=
F_{\gfk}(\tilde{M}). 
  $$

Then this is a left exact covariant functor.  
\smallbreak

There is a natural
$\Ufk_{\zi}(\bfk )$-module homomorphism
  $$
ev : H^0\bigl(\Ufk_{\zi}(\gfk )/\Ufk_{\zi}(\bfk ), M\bigr) \to M,
  $$
defined by
  $$
ev(f) = f(1).
  $$

For $M \in  \Ccal_{\zi}(\gfk )$  we have a  
$\Ufk_{\zi}(\gfk )$-module isomorphism:
$$ \theta : M \to  H^0\bigl(\Ufk_{\zi}(\gfk )/\Ufk_{\zi}(\bfk ), M\bigr) $$
given by $\theta (m) a = a.m ,$ for $m \in M $ and $a \in \Ufk_{\zi}(\gfk )$.

For a $\Ufk_{\zi}(\gfk )$-module $M$, the dual space $M^* := 
\Hom \,_{\zi}(M,\zi)$ is a\break $\Ufk_{\zi}(\gfk )$-module under 
$$ (af)m= f(S(a)m) \,, \,\hbox{for}\, a\in  \Ufk_{\zi}(\gfk ) \,,
f \in M^* \, \,\hbox{and}\, m\in M \,.$$

For $\lambda \in X^+$, there is a  $\Ufk_{\zi}(\gfk )$-module isomorphism
$$
\beta : V_\xi (\lambda)^\ast \to  
H^0\bigl(\Ufk_{\zi}(\gfk )/\Ufk_{\zi}(\bfk ), \chi_{-\lambda}^\xi\bigr)
 \tag{1}
$$
given by $\beta (f) (a) = f(S(a)v_\lambda) \,v_{\lambda}^*$, for $f\in V_\xi(\lam)^\ast$ and $ a \in 
\Ufk_{\zi}(\gfk )$, where  $V_\xi (\lambda)$ is the Weyl module over
$\zi$ with highest weight $\lambda$ (cf.\ [APW, Prop.~1.20~(ii)] or 
[KL, \S1]), $v_\lambda $ is a 
highest weight primitive vector of 
$ V_\xi (\lambda)$ and $v_{\lambda}^* \in 
\Hom ( V_\xi (\lambda)_\lambda , \zi)$ is given by  $v_{\lambda}^* ( 
v_{\lambda}) =1.$ (Observe that
 $\beta$ does not depend upon the choice of $v_\lam$.)

Exactly the same way we define the functor 
  $$
H^0\bigl(\bar{U}_{\Zgil}(\gfk )/\bar{U}_{\Zgil}(\bfk ),-\bigr) :
\Cbar_{\Zgil}(\bfk )\to \Cbar_{\Zgil}(\gfk ). 
  $$

We also need to consider the induction functor
  $$
H^0\bigl(\Ufk_{\zi}(\bfk )/\Ufk^0_{\zi}, -\bigr) : \Ccal^0_{\zi} \to
\Ccal_{\zi}(\bfk )
  $$
defined by
  $$
H^0\bigl(\Ufk_{\zi}(\bfk )/\Ufk^0_{\zi}, M\bigr) = F_{\bfk}
\Bigl(\Hom_{\Ufk^0_{\zi}}\bigl(\Ufk_{\zi}(\bfk ), M\bigr) \Bigr)\,,
  $$
where $\Ccal^0_{\zi}$ is the category of weight modules of
$\Ufk^0_{\zi}$.

Similarly, one defines the functor
  $$ H^0\bigl(\bar{U}_{\Zgil}(\bfk )/\bar{U}^0_{\Zgil}, -\bigr) : 
\bar{\Ccal}_{\Zgil}^0 \to \bar{\Ccal}_{\Zgil}(\bfk ). 
  $$

\proclaimtitle{[APW, Prop.\ 2.11
and Cor.\ 2.13]} 
   \proclaim{Proposition} 
\begin{itemize}
\ritem{a)} All the abelian categories $\Ccal^0_{\zi}$, $\Ccal_{\zi}(\bfk
)${\rm ,} $\Ccal_{\zi}(\gfk )$ have enough injective objects{\rm .}

\ritem{b)} The induction functors $H^0\bigl( \Ufk_{\zi}(\bfk
)/\Ufk^0_{\zi}, -\bigr)$ and $H^0\bigl(\Ufk_{\zi}(\gfk
)/\Ufk_{\zi}(\bfk ),-\bigr)$ take injective objects to injective
objects{\rm .}

\ritem{c)} The induction functor $H^0\bigl(\Ufk_{\zi}(\bfk
)/\Ufk^0_{\zi},-\bigr)$ is an exact functor{\rm ,} which takes $\zi$\/{\rm -}\/free
modules to $\zi$\/{\rm -}\/free modules{\rm .}
     \end{itemize}

An analogous result is true for the categories $\Cbar^0_{\zi},
\Cbar_{\zi}(\bfk ),\Cbar_{\zi}(\gfk )$ (and $ \Cbar^0_{\Zgil},\break
\Cbar_{\Zgil}(\bfk ), \Cbar_{\Zgil}(\gfk )$) as well{\rm . }
  \endproclaim

(Actually  the setting of the above proposition in [APW] is slightly different,
but the same proof works.)

Fix $M\in\Ccal_{\zi}(\bfk )$.  We need a certain specific resolution
of $M$ in the category $\Ccal_{\zi}(\bfk )$.  By the  Frobenius
reciprocity [APW, Prop.~2.12], $M$ is a $\Ufk_{\zi}(\bfk
)$-submodule of $Q_0 :=H^0\bigl(\Ufk_{\zi}(\bfk )/\Ufk^0_{\zi},
M\bigr)$ under $m\mapsto i_m$ (where  $i_m (X)=X.m)$ and, moreover, $M$ is 
a $\Ufk^0_{\zi}$ direct summand of
$Q_0$.  Apply the same to $Q_0/M$ and set $Q_1 :=
H^0\bigl(\Ufk_{\zi}(\bfk )/\Ufk^0_{\zi}, Q_0/M\bigr)$, etc.  This
gives a resolution of $M$ in $\Ccal_{\zi}(\bfk )$:
$$
0\to M\to Q_0 \to Q_1 \to \cdots .  \eqno(\ast)
  $$

If $M$ is $\zi$-free then so are each of the  $Q_i$.  We refer to ($*$)
as
the {\it standard resolution} of $M$ in the category
$\Ccal_{\zi}(\bfk )$.

  \numbereddemo{Definition}  Since the category $\Ccal_{\zi}(\bfk )$ has
enough injectives (by Proposition 1.5), the right derived functors of
$H^0\bigl(\Ufk_{\zi}(\gfk )/\Ufk_{\zi}(\bfk ),-\bigr)$ are defined.
Denote them by $H^i\bigl(\Ufk_{\zi}(\gfk )/\Ufk_{\zi}(\bfk
),-\bigr)$.

Similarly the $H^i\bigl(\bar{U}_{\Zgil}(\gfk )/\bar{U}_{\Zgil}(\bfk),-\bigr)$ are defined.

{We will abbreviate $H^i\bigl(\Ufk_{\zi} (\gfk )/\Ufk_{\zi} (\bfk ), -\bigr)$ 
(resp.\ $H^i\bigl(\Ubar_{\zi}
(\gfk )/\Ubar_{\zi} (\bfk ), -\bigr) $) by $H^i\bigl(\Xfk , -\bigr)$ (resp.\ $ H^i\bigl(\Xbar
, -\bigr)$).}
  \enddemo

\proclaimtitle{[APW, Prop.~2.19]}
  \proclaim{Proposition}   For any
$M\in\Ccal_{\zi}(\bfk )${\rm ,} the modules $Q_j$ in the standard
resolution ($*$) of $M$ have
  $$
H^i\bigl(\Xfk ,Q_j\bigr) = 0
\quad\mbox{for all }i>0, j\geq 0.
  $$
A similar result is true for 
$H^i\bigl(\Xbar ,\bar{Q}_j\bigr)${\rm .}
  \endproclaim

  \section{Definition of   quantized Frobenius homomorphism}

From now on we drop the subscript $\zi$ from $\Ufk^0_{\zi}$,
$\Ufk_{\zi}(\bfk )$ etc.; i.e., $\Ufk^0$ means $\Ufk^0_{\zi}$ etc.
Similarly, the category $\Cbar_{\zi}(\bfk )$ is abbreviated as $\Cbar
(\bfk )$ etc.

  \numbereddemo{Definition}  For any $\bar{M}\in\Cbar (\bfk )$, let
$\bar{M}^{\Fr}\in\Ccal (\bfk )$ be defined by taking $\bar{M}^{\Fr}=\bar{M}$ as a
$\zi$-module and the action of $a\in\Ufk (\bfk )$ on $\bar{M}^{\Fr}$ is
defined as
  \begin{equation}
a\odot m = \Fr (a)\cdot m.
  \end{equation}

Observe that, for $\lam\in X$,
  \begin{equation}
\bar{\chi}_{\lam}\circ\Fr_{|\Ufk^0} = \chi_{\ell\lam}^{\xi}.
  \end{equation}

To prove (2), use [Lu3, Lemma 34.1.2(c)].
\smallbreak
From (2) it is easy to see that $\bar{M}^{\Fr}$ is a weight module and
hence\break $\bar{M}^{\Fr}\in\Ccal (\bfk )$.

Clearly, for any $\bar{U}(\bfk )$-module morphism $f: \bar{M}\to \bar{N}$
($\bar{M},\bar{N}\in\Cbar (\bfk )$), the same map $f: \bar{M}^{\Fr}\to \bar{N}^{\Fr}$ is a
$\Ufk (\bfk )$-module morphism.

Exactly the same way  we can define $\bar{M}^{\Fr}\in\Ccal (\gfk )$ (resp.\ $\in\Ccal^0$) for $\bar{M}\in\Cbar (\gfk )$ (resp.\
$\in\Cbar^0$).
  \enddemo

 \proclaim{Lemma}     {\rm a)} For any $\bar{M}\in\Cbar^0${\rm ,} there
is a functorial $\Ufk (\bfk )$\/{\rm -}\/module  map
  $$
\Fr^*_{\bfk} :\bigl( H^0(\bar{U}(\bfk )/\bar{U}^0, \bar{M})\bigr)^{ \Fr}
\longrightarrow H^0\bigl(\Ufk (\bfk )/\Ufk^0, \bar{M}^{\Fr}\bigr),
  $$
defined by
 $$
(\Fr^*_{\bfk}f)(a) = f(\Fr (a)),  \tag{1}
  $$
for $f\in \bigl(H^0 (\bar{U}(\bfk )/\bar{U}^0, \bar{M} )\bigr)^{\Fr} :=
\left(F_{\bfk}\bigl(\Hom_{\bar{U}^0} (\bar{U}(\bfk ),\bar{M} )\bigr)\right)^{\Fr}$ and 
$a \in \Ufk (\bfk )${\rm . }

 \smallbreak
Similarly{\rm ,}

  {\rm b)} For any $\bar{M}\in\Cbar (\bfk )${\rm ,} there is a functorial $\Ufk
(\gfk )$\/{\rm -}\/module map
  $$
\Fr^* :  H^0(\bar{X} , \bar{M})^{\Fr}
\longrightarrow  H^0\bigl(\Xfk ,
\bar{M}^{\Fr}\bigr), 
 $$  
defined by 
$$
(\Fr^* f)(a) = f(\Fr (a)), \, \hbox{for}\, f\in   H^0(\bar{X} , \bar{M})^{\Fr}
\, \hbox{ and }\,  a \in   \Ufk (\gfk )\,. \tag{2}
  $$
  \endproclaim

{\it Proof}.  We prove (a); the proof of (b) is identical.
 
Clearly, $\Fr^*_{\bfk}f$ is a $\Ufk^0$-module map.  Moreover, for
$a,b\in\Ufk (\bfk )$,
  \begin{eqnarray*}
\bigl( a\cdot (\Fr^*_{\bfk}f)\bigr)\, b = (\Fr^*_{\bfk}f)\, (ba) &=&
f\,
(\Fr\, b\, \Fr\, a) \\[5pt]
 &= &\bigl( (\Fr\, a)\cdot f\bigr) (\Fr\, b).
  \end{eqnarray*}
    
This implies that
  $$
a\cdot (\Fr^*_{\bfk}f) = \Fr^*_{\bfk}\bigl( (\Fr\, a)\cdot f\bigr).
 \tag{3}
  $$

By (3) it is easy to see that $\Fr^*_{\bfk}f\in H^0\bigl(\Ufk (\bfk
)/\Ufk^0, \bar{M}^{\Fr}\bigr)$ and, moreover, $\Fr^*_{\bfk}$ is a $\Ufk
(\bfk
)$-module map.
  \hfill\qed\medbreak

Now we extend the above  $\Ufk (\gfk )$-module map $\Fr^*$ to an arbitrary 
$H^i$ , still keeping the same notation.

  \proclaim{Theorem}  For any $\bar{M}\in\Cbar (\bfk )${\rm ,} there exists a
 functorial $\Ufk (\gfk )$\/{\rm -}\/module map
 $$
\Fr^* : H^i (\bar{X} , \bar{M})^{ \Fr}
\longrightarrow H^i\bigl(\Xfk , \bar{M}^{\Fr}\bigr)
\eqno(1)
$$
in the sense that for any $\bar{U}(\bfk )$\/{\rm -}\/module morphism $\bar{M}\to \bar{N}$
{\rm (}$\bar{M},\bar{N}\in\Cbar (\bfk )${\rm ),} the following diagram is commutative
  $$  
    \begin{array}{ccc}
 H^i\bigl(\bar{X} , \bar{M}\bigr)^{ \Fr}&
\stackrel{\Fr^*}{\lrar}& H^i\bigl(\Xfk , \bar{M}^{\Fr}\bigr) \\ [5pt]
 \Big\downarrow&&\Big\downarrow\\ [9pt]
 H^i\bigl(\bar{X} , \bar{N}\bigr)^{ \Fr}&
\stackrel{\Fr^*}{\lrar}& H^i\bigl(\Xfk , \bar{N}^{\Fr}\bigr) \,,
     \end{array}  \tag{\rm D}
   $$
where the vertical maps are induced maps in cohomology{\rm .}\endproclaim   

  \demo{Proof} Let
  \begin{eqnarray}
&&0 \longrightarrow  \bar{M} \longrightarrow \bar{Q}_0
\stackrel{\bar{\veps}_0}{\longrightarrow}
\bar{Q}_1 \stackrel{\bar{\veps}_1}{\longrightarrow} \cdots ,
\speqnu{2}  \\[5pt]
&&0 \longrightarrow  \bar{M}^{ \Fr} \longrightarrow Q_0
\stackrel{\veps_0}{\longrightarrow}
Q_1 \stackrel{\veps_1}{\longrightarrow}\cdots ,
\speqnu{3}
  \end{eqnarray}  
be the standard resolutions in the  categories $\Cbar (\bfk )$ and $\Ccal
(\bfk )$ respectively (cf.\ ($*$) of $\S$1.5).  We construct by induction on 
$i$,
using  Lemma  2.2,  a $\Ufk (\bfk )$-module morphism
  $$
\theta_i :  \bar{Q}_i^{ \Fr} \to Q_i \qquad\mbox{for each }i\geq 0,
  $$
making the following squares commutative:
 \begin{equation} \begin{array}{cccccc}
   \bar{M}^{ \Fr} &\lrar& \bar{Q}_0^{ \Fr}&\lrar& \bar{Q}_1^{ \Fr} &\lrar \cdots\phantom{.}\\[5pt]
\Big\downarrow{\scriptstyle I}&& \Big\downarrow{\scriptstyle \theta_0}&& \Big\downarrow {\scriptstyle \theta_1} \\ [9pt]
  \bar{M}^{\Fr} &\lrar& Q_0 &\lrar& Q_1&\lrar\cdots  \,.
     \end{array}   \speqnu{D$'$}
  \end{equation}

  We first take
$i=0$.  Then, by definition,
  \begin{eqnarray*}
\bar{Q}_0 &=& H^0 (\bar{U}(\bfk )/\bar{U}^0, \bar{M}),\quad\mbox{
and}\\
Q_0 &= &H^0 (\Ufk (\bfk )/\Ufk^0, \bar{M}^{\Fr}).
  \end{eqnarray*}

Then $\theta_0$ is the map $\Fr^*_{\bfk}$ of Lemma 2.2(a).  Assume
now that we have constructed $\theta_j$ ($j\leq i$) making all the
squares in (D$'$) commutative up to $\theta_i$.  Now we construct
$\theta_{i+1}$ as follows: By definition
  \begin{eqnarray*}
\Qbar_{i+1} &=& H^0 (\Ubar (\bfk )/\Ubar^0,
\Qbar_i/\mbox{Image }\bar{\veps}_{i-1}),\quad\mbox{
and}\\
Q_{i+1} &= &H^0 (\Ufk (\bfk )/\Ufk^0, Q_i/\mbox{Image
}\veps_{i-1}).
  \end{eqnarray*}
By Lemma 2.2(a), we have the $\Ufk (\bfk )$-module map
  $$
\Fr^*_{\bfk}:(\bar{Q}_{i+1})^{\Fr} \longrightarrow H^0\bigl(\Ufk (\bfk
)/\Ufk^0, (\bar{Q}_i/\mbox{Image }\bar{\veps}_{i-1})^{\Fr}\bigr) .
  $$
From the commutativity of (D$'$) for the square containing
$\theta_{i-1}$ and $\theta_i$, we get a $\Ufk (\bfk )$ (and hence
$\Ufk^0$)-module map $\bar{\theta}_i: (\bar{Q}_i/\mbox{Image
}\bar{\veps}_{i-1})^{\Fr} \to Q_i/\mbox{Image }\veps_{i-1}$ induced by
$\theta_i$.  Inducing the map $\bar{\theta}_i$ via the functor
$H^0(\Ufk (\bfk )/\Ufk^0, -)$ and composing this with $\Fr^*_{\bfk}$
we get the desired $\Ufk (\bfk )$-module map $\theta_{i+1}$.  This
completes the induction.

The resolution (2) gives rise to a complex by taking
$H^0(\bar{X} ,-)^{\Fr}$:
  $$
0\to  H^0(\Xbar , \Qbar_0)^{\Fr} \to
 H^0(\Xbar , \Qbar_1)^{\Fr} \to
\cdots .  \tag{4}
  $$
Similarly, the resolution (3) gives rise to the complex:
  $$
  0\to H^0(\Xfk , Q_0) \to H^0 (\Xfk , Q_1) \to \cdots .
  \tag{5}  
  $$
Define the $\Ufk (\gfk )$-module map (for any $i\geq 0$)
  $$
\beta_i : H^0 (\Xbar ,
\Qbar_i)^{\Fr} \to H^0 \bigl(\Xfk , Q_i\bigr)
  $$
as the composite of 
  $\Fr^* : H^0(\Xbar , \Qbar_i)^{\Fr}
\to H^0\bigl(\Xfk , \Qbar_i^{\Fr}\bigr)$
(guaranteed by Lemma 2.2(b)) and the map $\theta^*_i : 
H^0 \bigl(\Xfk , \Qbar_i^{\Fr}\bigr)
\to H^0
\bigl(\Xfk , Q_i\bigr)$
  induced by the $\Ufk (\bfk
)$-module map $\theta_i : \Qbar_i^{\Fr} \to Q_i$\,.
  \smallbreak

The $\Ufk (\gfk )$-module maps $\beta_i$ give rise to a cochain map
from the cochain complex (4) to the cochain complex (5).  Taking
cohomology, we get the existence of the map (1).  (Observe that
$\Fr^*$ being an exact functor, the  $i^{\rm th}$ cohomology of the complex (4)
is the same as $H^i(\Xbar ,
\bar{M})^{\Fr}$.)  The functoriality of $\Fr^*$ follows from the
functoriality of all the constructions involved.
   \enddemo

\numbereddemo{{R}emark} (a) As we will see in a subsequent section, $\Fr^*$ is a 
quantization of the map induced on the cohomology of homogeneous vector 
bundles from the Frobenius morphism of the flag varieties $G/B$.
\smallbreak
(b) As informed by H.\ H.\ Andersen, for an extension $\zi \to k$  where
$k$ is a field, Theorem~2.3 
can also be deduced from [AW, Prop.\ 2.4]. 
\enddemo \pagebreak

  \section{Definition of quantized Frobenius splitting}

We continue to use the same abbreviation $\Ufk^0$ for $\Ufk_{\zi}^0$
etc. as given in the beginning of Section 2.

Recall the definition of the algebra homomorphism $\Fr'$ from Theorem~1.2.

Analogous to the Definition~2.1, we make the following.

  \numbereddemo{Definition}  For any $M\in\Ccal (\bfk )$, let
$M^{\Fr '}\in\Cbar (\bfk )$ be defined by taking $M^{\Fr'}=M$ as a
$\zi$-module and the action of $a\in\Ubar (\bfk )$ on $M^{\Fr'}$ is
defined by
  $$
  a\odot m = \Fr' a\cdot m.  \tag{1}
  $$

Observe that for $\lam\in X$, with $\lam =\lam_0 +\ell\lam_1$, where
$0\leq \lam_0 (h_i)\leq \ell -1$ for all $1\leq i\leq n$, and $\lam_1\in
X$, 
  $$
\chi_{\lam}^\xi\circ\Fr'_{|\Ubar^0} = \bar{\chi}_{\lam_1} . \tag{2}
  $$

To prove (2), again use [Lu3, Lemma 34.1.2(c)] or [KL, Lemma 3].

Similarly, for any $M\in\Ccal^0$, we define $M^{\Fr'}\in\Cbar^0$.

Clearly for any $\Ufk (\bfk )$-module morphism $f: M\to N$
($M,N\in\Ccal (\bfk )$), the same map $f: M^{\Fr'}\to N^{\Fr'}$ is a
$\Ubar (\bfk )$-module morphism.
  \enddemo

Exactly by the same proof as that of Lemma 2.2(a), we get the
following:

  \proclaim{Lemma}  For any $M\in\Ccal^0${\rm ,} there is a functorial $\Ubar
(\bfk )$\/{\rm -}\/module map
  $$
\Fr_{\bfk}^{\prime *}: H^0(\Ufk (\bfk )/\Ufk^0,
M)^{\Fr'} \to H^0\bigl(\Ubar (\bfk )/\Ubar^0, M^{\Fr'}\bigr),
  $$
defined by 
$$
(\Fr^{\prime*}_{\bfk}f)(a) = f(\Fr' (a)),  \tag{1}
  $$
for $f\in  H^0(\Ufk (\bfk )/\Ufk^0,
M)^{\Fr'},$ and $  a \in \Ubar (\bfk ).$ \hfill\qed
 \endproclaim

For a $\Ufk^0$-module $V$, $V^{\frac{1}{\ell}}$ denotes the sum of weight 
spaces corresponding to the weights $\lambda \in \ell X$.

 \proclaim{Proposition}  For any $M\in \Ccal (\bfk )${\rm ,} there is a
functorial $\Ubar (\bfk^- )$\/{\rm -}\/module map
    $$
\Frps : \Bigl( H^0\bigl(\Xfk ,
M\bigr)^{\frac{1}{\ell}}\Bigr)^{\Fr'_{\bfk^-}} \rightarrow
H^0\bigl(\Xbar , M^{\Fr'}\bigr),
 $$ 
defined by 
$$
(\Fr^{\prime*} f)(a) = f(\Fr' (a)),\,\hbox{for} \, f\in 
 \Bigl( H^0\bigl(\Xfk ,
M\bigr)^{\frac{1}{\ell}}\Bigr)^{\Fr'_{\bfk^-}} \, \hbox{and }\, a \in 
 \Ubar (\bfk^- ) . \tag{1}
 $$
 
Moreover{\rm ,} for any $m\geq 0$ and $f\in H^0(\Xfk , M)^{\frac{1}{\ell}}${\rm ,}
  $$
\Ebar_i^{(m)}\cdot (\Frps f) = \Frps (E_i^{(m\ell )}\cdot f).\tag{2}
  $$

Consider the projection 
$$ \pi : H^0(\Xfk , M)  \rightarrow  
H^0(\Xfk , M)^{\frac{1}{\ell}}$$
obtained by decomposing 
$$ H^0(\Xfk , M) =  
H^0(\Xfk , M)^{\frac{1}{\ell}} \oplus 
\bigg(H^0(\Xfk , M)^{\frac{1}{\ell}}\bigg)^\perp , $$
where $\perp$ is the sum of weight spaces corresponding to the weights 
$\lambda \notin \ell X${\rm .}
{\rm (}\/Observe that $\pi$ is a  $\Ubar (\bfk^- )$\/{\rm -}\/module map if   the module 
structures of both the domain and range are twisted by ${\Fr'}_{\bfk^-}$ and similarly $\pi$ is a  $\Ubar (\bfk )$-module map\/{\rm
.)} Composing $\Fr^{\prime*}$ with $\pi${\rm ,} we get a $\Ubar (\bfk^- )$\/{\rm -}\/module map 
{\rm (}\/again denoted by\/{\rm )} 
 $$
\Frps : H^0\bigl(\Xfk ,
M\bigr)^{\Fr'_{\bfk^-}} \rightarrow
H^0\bigl(\Xbar , M^{\Fr'}\bigr) .
  $$
\endproclaim

 {\it Proof}.  We have by the triangular decomposition 
[Lu2, Th.\ 6.7(d)],
 $$
\begin{array}{c}
   H^0 \bigl(\Xfk , M\bigr)^{\Fr_{\bfk^-}'}
\hookrightarrow \Hom_{\Ufk (\bfk )}\bigl(\Ufk (\gfk ),M\bigr)^{\Fr'_{\bfk^-}}
\,\begin{array}{c}{\scriptstyle \delta}\\ \noalign{\vskip-8pt} {\scriptstyle \simeq}\\[6pt]\end{array}
\Hom_{\Ufk^0}
\bigl(\Ufk (\bfk^-), M\bigr)^{\Fr_{\bfk^-}'} \\
 \hskip2.in \psi\downarrow    \\
\phantom{-}H^0 \bigl(\Xbar , M^{\Fr'}\bigr)
\hookrightarrow \Hom_{\Ubar (\bfk )} \bigl(\Ubar (\gfk
),M^{\Fr'}\bigr) 
\begin{array}{c}{\scriptstyle \bar\delta}\\ \noalign{\vskip-8pt} {\scriptstyle \simeq}\\[6pt] \end{array}
\Hom_{\Ubar^0}
\bigl(\Ubar (\bfk^-), M^{\Fr'}\bigr), \end{array} \tag{\rm D}
 $$  
where $\psi =\Fr^{\prime\ast}_{\bfk^-}$ is the $\Ubar
(\bfk^-)$-module map of Lemma~3.2  with $\bfk$ replaced by $\bfk^-$,
and $\delta, \bar{\delta}$ are the restriction maps.

We now show that 
  $$
\psi \Bigl(\Bigl( H^0\bigl(\Xfk ,
M\bigr)^{\frac{1}{\ell}}\Bigr)^{{\Fr'}_{\bfk^-}} \Bigr) \subset H^0\bigl(
\Xbar , M^{\Fr'}\bigr) . \tag{\ast}
  $$

Take a
weight vector $f\in H^0 (\Xfk ,M)$
 with respect to the $\Ufk^0$-action.  Then since $\psi$ is a 
$\Ubar (\bfk^-)$-module map, by (2) of $\S$(3.1), $\psi (f)$ is a
weight vector with respect to the $\bar{U}^0$-action.  Next, 
  $$
\psi \bigl( F_i^{(m\ell )}\cdot f\bigr) = \psi\bigl(
\bar{F}_i^{(m)}\odot f\bigr) = \bar{F}_i^{(m)}\cdot \psi (f)\, ,\,
\qquad \hbox{for any }\, m\geq 0\,.
  $$
In particular,  $\bar{F}_i^{(m)}\cdot\psi (f)=0$ for all large enough $m$.  So,
to prove the assertion ($*$), it suffices to show that for $f \in 
\Hom\,_{\Ufk (\bfk )}\bigl(\Ufk (\gfk ),M\bigr)^{\frac{1}{\ell} }$,
  $$
\psi \bigl( E_i^{(m\ell )}\cdot f\bigr) = \bar{E}_i^{(m)}\cdot \psi (f)\, ,\,
\hbox{for \, any}\, m\geq 0\, .\tag{3}
 $$
This will also prove (2).

As a  preparation, we prove the following lemmas.

  \proclaim{Lemma}   For any nonnegative integers $m, m_1, \cdots ,
m_r$ and $1\leq i,\break i_1,\cdots, i_r\leq n${\rm ,} let $1\leq
j_1 < \cdots < j_s \leq r$ be precisely the indices such that
$i_{j_p}=i${\rm .}

Then in the quantized enveloping algebra $U_{\Acal}(\gfk )$\/{\rm :}
  \begin{eqnarray}
F_{i_r}^{(m_r)} \cdots F_{i_1}^{(m_1)} E_i^{(m)}
& =& \sum_{t=(t_1,\cdots ,t_s)\in\Zgil_+^s}
E_i^{(m-\sum^s_{k=1}t_k)} F_{i_r}^{(m_r)} \speqnu{1}\\
&& \cdots\ 
F_{i_{j_s}}^{(m_{j_s}-t_s)} \cdots F_{i_{j_1}}^{(m_{j_1}-t_1)} \cdots
F_{i_1}^{(m_1)}\,A_t\,\, , \nonumber
  \end{eqnarray}
where $E_i^{(m' )}$ and  $F_i^{(m' )}$ are interpreted as 0 if $m' < 0${\rm ,}  
  \begin{eqnarray*}
 &&\hskip-20pt A_t :=  \bbm K_i^{-1}; m_{j_1}-m-\sum_{j<j_1}a_jm_j \\ t_1\ebm 
\bbm K_i^{-1}; m_{j_2}-m- \sum_{j<j_2} a_jm_j-t_1 \\ t_2\ebm \\ 
&&\hskip12pt \cdots\ \bbm K_i^{-1}; m_{j_s}-m-\sum_{j<j_s}
a_jm_j-(t_1+ \cdots +
t_{s-1}) \\ t_s\ebm , \\  
&&\hskip-20pt a_j   :=  -\alpha_{i_j}(h_i),\;\; \hbox{and}\,\, \bbm K_i^{-1}; c\\ t\ebm :=
\prod^t_{s=1} \frac{ K_i^{-1} v^{d_i(c-s+1)}-K_i v^{-d_i(c-s+1)} }
{ v^{d_is}-v^{-d_is} } .
  \end{eqnarray*}
  \endproclaim

  \demo{Proof}
Prove the lemma by induction on $r$, using the commutation
relations [Lu2, $\S$6.5], the following lemma and the $\Acal$-algebra 
automorphism
$\omega$ of $\Ufk_{\cal A} (\gfk )$ as in [Lu3, $\S$3.1.3].
  \hfill\qed

 \proclaim{Lemma}  For $t\in\Ngil := \{1,2, \cdots \}$ and $c\in\Zgil${\rm ,} as elements of $\Ufk_{\cal A}(
\frak J)${\rm ,}
  $$
  \bbm K_i^{-1}; c \\ t\ebm \, E_j^{(m)} = E_j^{(m)}
\bbm K_i^{-1}; c - a_{ij}m \\ t \ebm
  $$
and
  $$
  \bbm K_i^{-1}; c \\ t\ebm \, F_j^{(m)} = F_j^{(m)}
\bbm K_i^{-1}; c + a_{ij}m \\ t \ebm \,,
  $$
where $a_{ij} := \alpha_j(h_i)${\rm .}
  \endproclaim

  \demo{Proof}  Apply the automorphism $\omega$ of [Lu3, $\S$3.1.3]
to the identities [Lu2, $\S$6.5].
  \enddemo

  \proclaim{Lemma}  For any $m\in\Zgil$ and $t\geq 0$
  $$
 \Fr\bbm K_i; \ell m\\t\ebm =0
  $$
if $t$ is not divisible by $\ell${\rm ,}   where $[ \begin{array}{c}
{\scriptstyle K_i; \ell m}\\ \noalign{\vskip-5pt} {\scriptstyle t}\end{array}]$ is interpreted as an element of $\frak U_{\Bbb
Z_\xi}(\frak J)${\rm .}
  \endproclaim 

{\it Proof}.  First assume that $m<0$.  Then the lemma follows
from [Lu2,\break $\S$6.4-b3].  By [Lu2, $\S$6.5--a6], we have $[ \begin{array}{c}
{\scriptstyle K_i; \ell m}\\ \noalign{\vskip-5pt} {\scriptstyle t}\end{array}] F_i^{(\ell m)} = F_i^{(\ell m)}  [ \begin{array}{c}
{\scriptstyle  K_i; -\ell
m}\\ \noalign{\vskip-5pt} {\scriptstyle t}\end{array}]$.  From this the case $m>0$ also follows.
  \hfill\qed 

  \proclaim{Lemma}   For $t\in \Bbb Z_+${\rm ,} $c,a\in \Bbb Z${\rm ,}
  \begin{eqnarray}
\bbm K_i^{-1}; c \\ t\ebm& =& (-1)^t \bbm K_i; t-1-c\\ t \ebm \hbox{ as elements of } \frak U_{\cal A}(\frak J),\,\, \hbox{and} 
\speqnu{1}
\\
\bbm a\\ t\ebm^\xi_{d_i}  &=& (-1)^t \bbm -a+t-1\\ t \ebm_{d_i}^{\xi} , \speqnu{2}
\end{eqnarray} 
where  $[ \begin{array}{c} {\scriptstyle a}\\ \noalign{\vskip-5pt} {\scriptstyle t}\end{array}]^{\xi}_{d_i}$ denotes $[
\begin{array}{c} {\scriptstyle  a}\\ \noalign{\vskip-5pt} {\scriptstyle t}\end{array}]_{d_i} \in\Acal$ evaluated at $v=\xi${\rm .}

We also recall the $q$\/{\rm -}\/binomial identity (for $0\leq a_0 , t_0\leq \ell -1,
a_1\in\Zgil ,\break t_1\in\Zgil_+$) from {\rm [Lu3,} Lemma {\rm 34.1.2]:}
  $$
\bbm a_0+\ell a_1\\ t_0+\ell t_1\ebm^{\xi}_{d_i} = \bpm a_1\\
t_1\epm\, \bbm a_0\\ t_0\ebm^{\xi}_{d_i}, \tag{3}
  $$
where ${a_1\choose t_1} \in\Zgil_+$ is the ordinary binomial
coefficient{\rm .}
  \endproclaim 

  {\it Proof}. The first identity follows from the definition. 
  For the second
identity see [Lu3, p.~266].
\hfill\qed \vglue6pt

{\it Proof of} (3)  {\it of} \S 3.3  {\it continued}.
 First take $f\in\Hom_{\Ufk (\bfk )}(\Ufk (\gfk ),M)$.  Then, by Lemma~3.4 (following the same notation),
  \begin{eqnarray}
\psi\bigl( E_i^{(m\ell )}\cdot f\bigr)\bigl(\Fbar_{i_r}^{(m_r)}
\cdots \Fbar_{i_1}^{(m_1)}\bigr) 
&= &f\bigl( F_{i_r}^{(\ell m_r)}
\cdots F_{i_1}^{(\ell m_1)} E_i^{(m\ell )}\bigr)  \speqnu{1}\\
&=& \sum_t E_i^{(\ell m-\sum^s_{k=1}t_k)} \cdot
\bigl(( A_t\cdot f)(F_t)\bigr),  \nonumber
  \end{eqnarray} 
where $F_t := F_{i_r}^{(\ell m_r)} \cdots
F_{i_{j_s}}^{(\ell m_{j_s}-t_s)} \cdots F_{i_{j_1}}^{(\ell
m_{j_1}-t_1)} \cdots F_{i_1}^{(\ell m_1)}.$

Now assume that $f\in\Hom_{\Ufk (\bfk )}(\Ufk (\gfk
),M)^{\frac{1}{\ell}}$, and $f$ is of weight $\ell\lam$, for $\lam\in
X$.
Then (by Lemma 3.7) the above sum reduces to 
 \begin{eqnarray}
&& \sum_t\, E_i^{(\ell m-\sum^s_{k=1}t_k)} \speqnu{$\ast$}\\
&& \cdot \Bigl(\left[ \begin{array}{c}
-\ell\lam (h_i)+\ell p'_1\\ t_1\ebm^{\xi}_{d_i} \cdots 
\bbm -\ell\lam (h_i)+\ell p'_s-(t_1 + \cdots + t_{s-1})\\
t_s\end{array}\right]^{\xi}_{d_i} 
f(F_t)\Bigr),  \nonumber
\end{eqnarray}
for some $ p'_1, \cdots ,p'_s \in {\Zgil}$.
If at least one of $t_1,\cdots ,t_s$ is not divisible by $\ell$, say
$t_j$, and $t_j$ is the first one with this
property, then
  $$
\bbm -\ell\lam (h_i)+\ell p'_j-(t_1 + \cdots + t_{j-1})\\
t_j\ebm^{\xi}_{d_i} = 0,
  $$
by (3) of Lemma~3.7.

So the sum $(\ast )$ reduces to $t=(t_1,\cdots ,t_s)$ such that $\ell
\vert t$, i.e., each $t_k$ is divisible by $\ell$ giving
  $$
\psi\bigl( E_i^{(m\ell )}\cdot f\bigr)\bigl(\Fbar_{i_r}^{(m_r)}
\cdots \Fbar_{i_1}^{(m_1)}\bigr) 
= \sum_{\ell \vert t} E_i^{(\ell m-\sum^s_{k=1}t_k)} \cdot
\bigl(( A_t\cdot f)(F_t)\bigr)\,.
\tag{2}
  $$

On the other hand, applying $\Fr$ to the commutation relation as in 
Lemma~3.4, we get 

  \begin{eqnarray}
&&\bigl(\Ebar_i^{(m)}\cdot (\psi f)\bigr)  \bigl(\Fbar_{i_r}^{(m_r)}
\cdots \Fbar_{i_1}^{(m_1)}\bigr)  \speqnu{3}
 \\
&&\qquad\qquad= \sum_{\ell \vert t}\, \Bigl(\Fr' \bigl(\Fr\,  E_i^{(\ell
m-\sum^s_{k=1}t_k)}\bigr)\Bigr)\cdot \bigl( ((\Fr' \Fr\, A_t)\cdot f)(F_t)
\bigr) 
\nonumber \\
&&\qquad\qquad = \sum_{\ell \vert t}\,E_i^{(\ell m-\sum^s_{k=1}t_k)} \cdot
\bigl(( A_t\cdot f)(F_t)\bigr)\,. 
 \nonumber  \end{eqnarray}
Comparing (2) and (3), we get
  $$
\psi\bigl( E_i^{(m\ell )}\cdot f) = \Ebar_i^{(m)}\cdot (\psi f),
\mbox{  for all } f\in\Hom_{\Ufk (\bfk )} (\Ufk (\gfk ),
M)^{\frac{1}{\ell}} .
  $$
This proves (3) of Proposition~3.3  and hence  Proposition~3.3  itself.\hfill\qed
  \medbreak

Now we extend the $\bar{U}(\bfk^-)$-module map $\Fr^{\prime*}$ of Proposition~3.3  to an arbitrary cohomology $H^i$. 
\proclaim{Theorem} For any $M\in \Ccal (\bfk )${\rm ,} there exists a
functorial $\Ubar (\bfk^- )$\/{\rm -}\/module map for all $i\geq 0${\rm :}
 $$
\Frps : H^i\bigl(\Xfk ,
M\bigr)^{\Fr'_{\bfk^-}} \rightarrow
H^i\bigl(\Xbar , M^{\Fr'}\bigr) .
  $$ 
Moreover{\rm ,} for any $m\geq 0$ and $f\in H^i(\Xfk , M)${\rm ,}
  $$
\Ebar_i^{(m)}\cdot (\Frps f) = \Frps (E_i^{(m\ell )}\cdot f),\tag{1}
  $$
i.e.{\rm ,} $\Frps : H^i\bigl(\Xfk ,
M\bigr)^{\Fr'_{\bfk}} \rightarrow
H^i\bigl(\Xbar , M^{\Fr'}\bigr)$ is a  $\Ubar (\bfk )$\/{\rm -}\/module map as well{\rm ,} 
for all $i \geq 0${\rm .}\endproclaim 
  
  \demo{Proof}  The proof  is  parallel to the proof of 
Theorem~2.3.
Let
  \begin{eqnarray*}
&&0 \longrightarrow  M \longrightarrow Q_0
\stackrel{\veps_0}{\longrightarrow}
Q_1 \stackrel{\veps_1}{\longrightarrow} \cdots ,\\[5pt]
&&0 \longrightarrow  {M}^{\Fr'} \longrightarrow \bar{Q}_0
\stackrel{\bar{\veps}_0}{\longrightarrow}
\bar{Q}_1 \stackrel{\bar{\veps}_1}{\longrightarrow} \cdots ,
  \end{eqnarray*}
be the standard resolutions in categories $\Ccal (\bfk )$ and $\Cbar
(\bfk )$ respectively.  By induction,  we construct
   $\bar{U} (\bfk )$-module morphisms
  $\theta_j : Q_j^{ \Fr'} \to  \bar{Q}_j $
making the  squares commutative up to $\theta_j$.

  First define 
$$ \theta_0:  Q_0^{\Fr'} := H^0 (\Ufk (\bfk )/\Ufk^0, M)^{\Fr'} \to 
\bar{Q}_0:=  H^0 (\bar{U}(\bfk )/\bar{U}^0, M^{\Fr'})$$
as the map $ \Fr^{\prime*}_{\bfk}$ of Lemma 3.2. Having defined $\theta_j$, 
define (abbreviating Image by Im)
$$ \theta_{j+1}:  Q_{j+1}^{\Fr'} := H^0 (\Ufk (\bfk )/\Ufk^0,  
Q_j/\mbox{Im}\,\veps_{j-1})^{\Fr'} \to 
\bar{Q}_{j+1}:=  H^0 (\bar{U}(\bfk )/\bar{U}^0, 
\Qbar_j/\mbox{Im }\bar{\veps}_{j-1})$$
as the composite 
\begin{eqnarray*}
 H^0 (\Ufk (\bfk )/\Ufk^0,  
Q_j/\mbox{Im}\,\veps_{j-1})^{\Fr'}& \stackrel{ \Fr^{\prime*}_{\bfk}}
{\longrightarrow}& H^0 \bigl(\bar{U} (\bfk )/\bar{U}^0,  
(Q_j/\mbox{Im}\,\veps_{j-1})^{\Fr'}\bigr) \\
&\longrightarrow& 
H^0 (\bar{U}(\bfk )/\bar{U}^0, 
\Qbar_j/\mbox{Im }\bar{\veps}_{j-1})\,,
\end{eqnarray*}
where the second map is induced from the $\bar{U}(\bfk)$-module map $\theta_j$.
Finally, define a cochain map $H^0(\Xfk , Q_{\bullet})^{\Fr'_{\bfk^-}} 
\to
H^0(\bar{X}, \bar{Q}_{\bullet})$ as the composite map
$$ H^0(\Xfk , Q_{\bullet})^{\Fr'_{\bfk^-}} \stackrel{\Fr^{\prime*}}
{\longrightarrow}
H^0(\bar{X}, Q_{\bullet}^{\Fr'}) \to
H^0(\bar{X}, \bar{Q}_{\bullet})\,,
$$
where the second map is induced from the $\bar{U}(\bfk )$-module maps 
$\theta_{\bullet}$. This proves the theorem.
\enddemo

\proclaim{{C}orollary} For any $\bar{M}\in \bar{\Ccal} (\bfk )${\rm ,} the composite map
 $$
\Frps \circ \Fr^* : H^i\bigl(\Xbar ,
\bar{M}\bigr) \rightarrow
H^i\bigl(\Xbar , \bar{M}\bigr) 
  $$ 
is the identity map  for all $i\geq 0${\rm .}
\endproclaim 

{\it Proof}. It is easy to see that the corollary holds for $H^0$.
To prove the result for general $i$, take an exact sequence in  $
\bar{\Ccal} (\bfk ):\,
0 \to \bar{M} \to  \bar{N} \to \bar{Q} \to 0 $  such that
$H^i(\bar{X}, \bar{N}) = 0 $ , for all $ i \geq 1.$
 Then, from the surjective map  $H^{i-1}(\bar{X}, \bar{Q}) \twoheadrightarrow 
H^i(\bar{X}, \bar{M})$ and the functoriality of $ \Fr^*$ and $\Frps$, the 
corollary for $i$ follows by induction.
 \hfill\qed
\medbreak
{\it Remark} 3.10.
 (a) As we will see in a subsequent section, $\Frps$ is a 
quantization of the map induced on the cohomology of homogeneous vector 
bundles from the `canonical' Frobenius splitting of the flag variety $G/B$
obtained by Mehta-Ramanathan. Thus the key lemma of Mathieu (asserting that 
a\break $B$-canonical splitting of a $B$-variety $Y$ sends any $B$-submodule of 
$H^0(Y, \Lfk^{\otimes p})$ to a $B$-submodule of $H^0(Y, \Lfk)$, for any 
$B$-equivariant line bundle $\Lfk$ on $Y$;\break cf.\ [M, Lemma 2.4 and the remark following it]) in this case 
follows from the fact in Theorem~3.8  that the splitting is a 
$\bar{U}(\bfk^- )$-module (as well as a $\bar{U}(\bfk )$-module) map.
 
\smallbreak
(b) For $\lam \notin \ell X$, the map $\Frps : H^0\bigl(\Xfk ,
\chi^\xi_{\lam}\bigr)^{\Fr'_{\bfk^-}} \rightarrow
H^0\bigl(\Xbar ,(\chi^\xi_{\lam})^{\Fr'}\bigr)$ is identically zero. To
see this, write $\lam = \lam_0 + \ell \lam_1$ with $0 \leq \lam_0(
h_i) \leq \ell - 1$ for all simple coroots $h_i$ and
$ 0< \lam_0(h_{i_o})$ for some  $h_{i_o}$. Now
take any $f \in  H^0\bigl(\Xfk ,
\chi^\xi_{\lam}\bigr)^{1/\ell} $. Then, by Lemma~3.7,  
 $ \bbm K_{i_o}; 0\\ 
\lam_0(h_{i_o})\ebm \cdot f = 0\,.$ This gives that  $\phantom{\displaystyle\sum^{\raise2pt\hbox{$\int$}}}$
$$\left(  \bbm K_{i_o}; 0\\ 
\lam_0(h_{i_o})\ebm \cdot f \right)\bigl( F_{i_r}^{(\ell m_r)}
\cdots F_{i_1}^{(\ell m_1)}\bigr) =0 $$ for any nonnegative $m_1, \cdots, 
m_r$. Hence
\begin{eqnarray*}
\left(  \bbm K_{i_o}; 0\\ 
\lam_0(h_{i_o})\ebm \cdot f \right)\bigl( F_{i_r}^{(\ell m_r)}
\cdots F_{i_1}^{(\ell m_1)}\bigr) &= &f\left( F_{i_r}^{(\ell m_r)}
\cdots F_{i_1}^{(\ell m_1)}  \bbm K_{i_o}; 0\\ 
\lam_0(h_{i_o})\ebm \right)\\
&=&cf\bigl( F_{i_r}^{(\ell m_r)}
\cdots F_{i_1}^{(\ell m_1)}\bigr)=0\,,
\end{eqnarray*}
for some nonzero $c$. 
Thus we conclude that 
$f( F_{i_r}^{(\ell m_r)}
\cdots F_{i_1}^{(\ell m_1)}) = 0\,,$ and hence $\Frps f = 0$ \,(by the 
definition of $\Frps ).$

  \section{Stronger quantized Frobenius splitting}

{\it In this section we abbreviate the homomorphism $\Fr'_{\bfk^-}$ 
of Theorem~{\rm 1.2} to} $\Fr'$. We also continue to 
abbreviate 
$H^i (\Ufk (\gfk )/\Ufk (\bfk ), -)$ (resp.\ $H^i (\Ubar
(\gfk )/\Ubar (\bfk ), -) $) to $H^i (\Xfk , -)$ (resp.\ $ H^i (\Xbar
, -) $). Any reduced decomposition $w_o = s_{i_1} \cdots s_{i_N}$ of the 
longest element $w_o$ of the Weyl 
group in terms of the simple reflections gives an indexing of the set $\Delta_+$ of positive roots $\{\beta_1, \cdots , \beta_N\},$
where $\beta_j := s_{i_1} \cdots s_{i_{j-1}}(\alpha_{i_j})$ ($\alpha_i$ being the simple root corresponding to the simple
reflection $s_i$).

  \numbereddemo{Definition} For any $\Ufk^0$-module $M \in \Ccal^0$, define the map
(abbreviating $-2(\ell - 1)\rho$ to $\gamma$, where $\rho$ is the half sum
of positive roots)
  $$
\pch : \Hom\,_{\Ufk^0} \bigl(\Ufk (\bfk^-), \chi^{\xi}_{\gamma}
\otimes M\bigr)^{\Fr'} \rightarrow \Hom\,_{\Ubar^0}
\bigl(\Ubar (\bfk^-), M^{\Fr'}\bigr)
  $$
by
  $$
(\pch f)(a) = f(F_o \Fr'(a))\otimes v_+
\quad\mbox{for }a\in\Ubar (\bfk^-)\,,
  $$
where $F_o := F^{(\ell -1)}_{\beta_N} \cdots F^{(\ell -1)}_{\beta_1}$ , 
$F_{\beta_i}$ are Lusztig's root vectors [Lu2, \S4], 
 $v_+$ is a ${\Bbb Z}_{\xi}$-basis vector of the one-dimensional 
representation
$\chi^{\xi}_{2(\ell -1)\rho}$ and we identify $\chi^{\xi}_{\gamma} \otimes M
\otimes \chi^{\xi}_{-\gamma}$ with $M$ ($F_o$ does not depend upon the choice 
of the reduced decomposition of $w_o$ up to a nonzero scalar multiple, since the 
corresponding weight space  in the quantized restricted enveloping algebra 
is of rank one).  By the following lemma,
$\pch (f)$ is indeed $\Ubar^0$-linear.
  \enddemo

  \proclaim{Lemma} 
With   notation and assumptions as above{\rm ,} $\pch (f)$ is $\Ubar^0$\/{\rm -}\/linear for
any $f\in\Hom\,_{\Ufk^0}\bigl(\Ufk (\bfk^-), 
\chi^{\xi}_{\gamma}\otimes M\bigr)^{\Fr'}${\rm .}
  \endproclaim

  \demo{Proof}  For any $a\in\Ubar (\bfk^-)$ and $h={h_i\choose
m} \in\Ubar^0$,
  \begin{eqnarray*}
  (\pch f)(ha) &\hskip-4pt = \hskip-4pt & f (F_o\, \Fr'(h)\, \Fr'(a))\otimes v_+ \\[5pt]
&\hskip-4pt = \hskip-4pt & f\Biggl( F_o \bbm K_i; 0\\[5pt] \ell m\ebm\,\Fr'(a)\Biggr)\otimes
  v_+\\ 
&\hskip-4pt = \hskip-4pt &f\Biggl(\bbm K_i; 2(\ell -1)\\[5pt] \ell m\ebm\, F_o \Fr'(a)\Biggr)
  \otimes v_+, \quad\mbox{by [Lu2, $\S$6.5]}\\ 
&\hskip-4pt = \hskip-4pt & \Biggl( \bbm K_i; 2(\ell -1)\\[5pt] \ell m\ebm \cdot f(F_o
 \Fr'(a))\Biggr)\otimes v_+\\ 
&\hskip-4pt = \hskip-4pt & \bbm K_i; 0 \\[5pt] \ell m\ebm \cdot \Bigl( f (F_o \Fr'(a))\otimes
 v_+\Bigr) ,  \mbox{ since $v_+$ is of weight $2(\ell -1)\rho$}\\ 
&\hskip-4pt = \hskip-4pt & \bpm h_i\\[5pt] m\epm \odot ((\pch f)(a)).\\ 
\noalign{\vskip-36pt}
  \end{eqnarray*}
  \enddemo

\phantom{What?} 
  
\proclaim{Lemma} \hglue-8pt 
With notation and assumptions as in  ${{\rm \S 4.1,}}$ $\pch$ is $\Ubar
(\bfk^-)$\/{\rm -}\/linear{\rm .}
  \endproclaim 

\vglue-8pt

 {\it Proof}.  For $a,b\in \Ubar (\bfk^-)$, 
  \begin{eqnarray*}
   \bigl(\pch\bigl(b\odot f\bigr)\bigr) (a) &=&
 \bigl(\pch\bigl(\Fr'(b)\cdot f\bigr)\bigr) (a)\\[5pt]
 &= &\bigl( (\Fr'(b)\cdot f)(F_o
\Fr'(a))\bigr) \otimes v_+ \\[5pt]
&=& f(F_o\Fr'(a) \Fr'(b))\otimes v_+\\[5pt]
&=& f(F_o\Fr'(ab))\otimes v_+\\[5pt]
&=& (\pch f)(ab)\\[5pt]
&=& (b\cdot\pch f)(a).
  \end{eqnarray*}

This proves that $\pch$ is $\Ubar (\bfk^-)$-linear.
  \hfill\qed\medbreak

We now prove the following crucial proposition.

  \proclaim{Proposition}  For any $\Ubar (\bfk )$\/{\rm -}\/module $\Mbar$ which is a
$\Ubar^0$\/{\rm -}\/weight module{\rm ,}  $m\geq 0$ and
$f\in\Hom\,_{\Ufk (\bfk )} \bigl(\Ufk (\gfk ), 
\chi^{\xi}_{\gamma}\otimes\Mbar^{\Fr}\bigr) \simeq \Hom\,_{\Ufk^0}\bigl(\Ufk
(\bfk^-), \chi^{\xi}_{\gamma}\otimes\Mbar^{\Fr}\bigr)${\rm ,}
  $$
\pch \bigl( E_i^{(\ell m)}\cdot f\bigr) = \Ebar_i^{(m)}\cdot (\pch
f),
\tag{1}
  $$
where the action of $ \Ebar_i^{(m)}$ on $\pch
f$ comes from the similar identification $$\Hom\,_{\Ubar (\bfk)}
\bigl(\Ubar (\gfk), \Mbar\bigr) \simeq \Hom\,_{\Ubar^0}
\bigl(\Ubar (\bfk^-), \Mbar\bigr).$$
  \endproclaim

\phantom{strange}
\vglue-36pt
 {\it Proof}.  Take $a=\Fbar^{(m_{r})}_{i_{r}} \cdots
\Fbar^{(m_1)}_{i_1}$.  Then 
  \begin{eqnarray}
\bigl(\pch ( E_i^{(\ell m)}\cdot
f)\bigr) (a) &= &(E_i^{(\ell m)}\cdot f) (F_o \Fr'(a))\otimes v_+\speqnu{2}\\
&=& f\bigl( F_o \, F^{(\ell m_{r})}_{i_{r}} \cdots F^{(\ell
m_1)}_{i_1}\, E_i^{(\ell m)}\bigr) \otimes v_+ . \nonumber
  \end{eqnarray}
Now, by (1) of Lemma  3.4  and  Lemma 3.5, we get
  \begin{eqnarray}
&&\speqnu{3}\\ \noalign{\vskip-5pt}
F_o\, F^{(\ell m_{r})}_{i_{r}} \cdots F_{i_1}^{(\ell m_1)}
E_i^{(\ell m)} &=& \sum_t\, \hat{A}_t\, F_o\, E_i^{ (\ell
m-\sum^s_{k=1}t_k )} F_{i_{r}}^{(\ell m_{r})} \cdots
F^{(\ell m_{j_s}-t_s)}_{i_{j_s}}\nn\\
&& \cdots  F_{i_{j_1}}^{(\ell
m_{j_1}-t_1)} \cdots F^{(\ell m_1)}_{i_1},  \nonumber
  \end{eqnarray}
where the summation is over $t=(t_1,\cdots ,t_s)\in\Zgil^s_+$ ,
$1\leq j_1 < \cdots < j_s\leq r$ are precisely the indices
such that $i_{j_{p}}=i$ and
  \begin{eqnarray*}
  &&\hskip-.35in\hat{A}_t  :=  \bbm K_i^{-1}; -\ell m_{j_1}+\ell m + (\sum_{j>j_1}
a_j\, \ell m_j) -2(\ell -1)\\ t_1\ebm \\ 
 && \,\times\ \bbm K_i^{-1}; -\ell
m_{j_2}+\ell m+(\sum_{j>j_2} a_j\ell m_j )-t_1-2(\ell -1)\\
t_2\ebm \\
 &&\, \cdots \ \bbm K_i^{-1}; -\ell m_{j_s}+\ell m+(\sum_{j>j_s}
a_j\ell m_j)-(t_1+\cdots +t_{s-1}) -2(\ell -1)\\ t_s\ebm 
  \end{eqnarray*}
(where $a_j := -\alpha_{i_j}(h_i)$).
Substituting (3) in (2), we get
  \begin{eqnarray}
&& \bigl( \pch (E_i^{(\ell m)}\cdot f)\bigr) (a)    = \sum_{t\in\Zgil_+^s}
\biggl(\hat{A}_t\cdot \biggl( f\biggl( F_o E_i^{(\ell
m-\sum^s_{k=1}t_k)}\, F_{i_{r}}^{(\ell m_{r})} \cdots
F_{i_{j_s}}^{(\ell m_{j_s}-t_s)} \speqnu{4}\\ 
&&\hskip2.5in\cdots F_{i_{j_1}}^{(\ell
m_{j_1}-t_1)} \cdots F_{i_1}^{(\ell m_1)}\biggr)\biggr)\biggr)\otimes
v_+ .  \nonumber
  \end{eqnarray}
Since Im $f\subset \chi^{\xi}_{\gamma}\otimes \Mbar^{\Fr}$ ,
using (3) of Lemma~3.7,   the sum in (4) reduces to $(t_1,\cdots
,t_s)\in\Zgil_+^s$ such that each $t_k$ is divisible by $\ell$, i.e.,
denoting $\ell t  = (\ell t_1,\cdots ,\ell t_s)$, we have
  \begin{eqnarray}   
 && \hskip-6pt \bigl( \pch (E_i^{(\ell m)}\cdot f)\bigr)
(a)  \speqnu{5}\\ 
&&  =  \sum_{t\in\Zgil_+^s}
\biggl(\hat{A}_{\ell t}\cdot \biggl( f\biggl( F_o E_i^{(\ell
m-\sum^s_{k=1}\ell t_k)}\, F_{i_{r}}^{(\ell m_{r})}\nonumber\\ 
&&\hskip1.25in \cdots
F_{i_{j_s}}^{(\ell m_{j_s}-\ell t_s)} \cdots F_{i_{j_1}}^{(\ell
m_{j_1}-\ell t_1)}    \cdots F_{i_1}^{(\ell
m_1)}\biggr)\biggr)\biggr)\otimes v_+\nonumber \\ 
&&  = \sum_{t\in\Zgil_+^s}
\biggl(\hat{A}_{\ell t} E_i^{(\ell m-\sum^s_{k=1}\ell t_k)} 
\cdot \biggl( f\biggl( F_o 
F_{i_{r}}^{(\ell m_{r})}\cdots
F_{i_{j_s}}^{(\ell m_{j_s}-\ell t_s)}\nonumber \\&&\hskip1in  \cdots F_{i_{j_1}}^{(\ell
m_{j_1}-\ell t_1)} \cdots F_{i_1}^{(\ell
m_1)}\biggr)\biggr)\biggr)\otimes v_+\,,\,\,
\hbox{by \,the \,next \,lemma}\nonumber \\
&&=  \sum_{t\in\Zgil_+^s}
\bigl(\tilde{A}_{\ell t} E_i^{(\ell m-\sum^s_{k=1}\ell t_k)} \bigr)\nonumber\\[5pt]
&&\quad  \odot \biggl( f\biggl( F_o 
F_{i_{r}}^{(\ell m_{r})} \cdots
F_{i_{j_s}}^{(\ell m_{j_s}-\ell t_s)}\cdots F_{i_{j_1}}^{(\ell
m_{j_1}-\ell t_1)} \cdots F_{i_1}^{(\ell
m_1)}\biggr)\otimes v_+\biggr)\nonumber\\[5pt] 
&&=  \sum_{t\in\Zgil_+^s}\Fr\biggl(
\tilde{A}_{\ell t} E_i^{(\ell m-\sum^s_{k=1}\ell t_k)}\biggr) \nonumber\\[5pt]
&&\quad
\cdot \biggl( f\biggl( F_o 
F_{i_{r}}^{(\ell m_{r})} \cdots
F_{i_{j_s}}^{(\ell m_{j_s}-\ell t_s)} \cdots F_{i_{j_1}}^{(\ell
m_{j_1}-\ell t_1)}  \cdots F_{i_1}^{(\ell
m_1)}\biggr)\otimes v_+\biggr), \nonumber 
   \end{eqnarray}
where 
\begin{eqnarray*}
 \tilde{A}_t& := &\bbm K_i^{-1}; -\ell m_{j_1}+\ell
m+{\displaystyle \sum_{j>j_1}} a_j\ell m_j\\ t_1 \ebm\\[5pt]
&& \cdots  \
\bbm K_i^{-1}; -\ell m_{j_s}+\ell m+{\displaystyle \sum_{j>j_s}}
a_j\ell m_j-(t_1+\cdots +t_{s-1})\\ t_s\ebm \,. 
\end{eqnarray*}
We now calculate the right side of (1):
  \begin{eqnarray}
&&\hskip-6pt
\bigl( \Ebar_i^{(m)}\cdot (\pch f)\bigr)(a) \speqnu{6}\\[5pt] 
&&=
(\pch f)\bigl(\Fbar_{i_{r}}^{(m_{r})} \cdots
\Fbar_{i_1}^{(m_1)} \Ebar_i^{(m)}\bigr)\nn \\[5pt]
&&=\sum_{t\in\Zgil_+^s} (\pch f)\biggl(
H_t\Ebar_i^{(m-\sum^s_{k=1}t_k)}\, \Fbar_{i_{r}}^{(m_{r})}
\cdots \Fbar_{i_{j_s}}^{(m_{j_s}-t_s)} \cdots
\Fbar_{i_{j_1}}^{(m_{j_1}-t_1)}\cdots \Fbar_{i_1}^{(m_1)}\biggr) \nn \\[5pt]
&&=\sum_{t\in\Zgil_+^s} H_t\Ebar_i^{(m-\sum^s_{k=1}t_k)}\nn\\[5pt]
&&\quad \cdot\biggl( (\pch f)  \biggl(\Fbar_{i_{r}}^{(m_{r})}
\cdots \Fbar_{i_{j_s}}^{(m_{j_s}-t_s)} \cdots
\Fbar_{i_{j_1}}^{(m_{j_1}-t_1)}\cdots
\Fbar_{i_1}^{(m_1)}\biggr)\biggr) \nn \\[5pt]
 &&=\sum_{t\in\Zgil_+^s} H_t\Ebar_i^{(m-\sum^s_{k=1}t_k)}\nn\\[5pt] &&\quad 
\cdot\biggl( f  \biggl(F_o\, F_{i_{r}}^{(\ell m_{r})}
\cdots F_{i_{j_s}}^{(\ell m_{j_s}-\ell t_s)} \cdots
F_{i_{j_1}}^{(\ell m_{j_1}-\ell t_1)}\cdots
F_{i_1}^{(\ell m_1)}\biggr)\otimes v_+\biggr)\,,  \nn
  \end{eqnarray}
where $H_t := \Fr (\tilde{A}_{\ell t})$. 

Comparing (5) and (6) we get (1).  This proves the proposition modulo the
next lemma.
  \hfill\qed
\pagebreak

  \proclaim{Lemma}  For any $\Ubar (\bfk )$\/{\rm -}\/module $\Mbar$ such that
$\Mbar$ is a $\Ubar^0$\/{\rm -}\/weight module\/{\rm ,}\/  $m\geq 0$ and
$f\in\Hom\,_{\Ufk (\bfk )} \bigl(\Ufk (\gfk ), \chi^{\xi}_{\gamma}
 \otimes\Mbar^{\Fr}\bigr)$\/{\rm ,}\/ 
  $$
f(F_o\, E_i^{(\ell m)}) = f(E_i^{(\ell m)}\, F_o). \tag{1}
  $$
Thus{\rm ,} replacing $f$  by $x\cdot f$\/{\rm ,} 
$$f(F_o\, E_i^{(\ell m)} x) = 
f(E_i^{(\ell m)}\, F_o x)\,\,\,\hbox{ for any }\, x \in \Ufk (\gfk).
$$
  \endproclaim

\vglue-12pt
 {\it Proof}.  Any such $\Mbar$ is a quotient of a $\Ubar (\bfk
)$-module $\Qbar$ such that $\Qbar$ is a $\Ubar^0$-weight module and
$\Qbar$ is $\zi$-free   (since, for any weight vector $v\in\Mbar$
of weight $\lam$, there exists a $\Ubar (\bfk )$-module map $\pi_v :
\Ubar (\bfk )\otimes_{\Ubar^0} \bar{\chi}_{\lam}\to\Mbar$ taking
$1\otimes 1\mapsto v$).  Now the surjective $\bar{U} (\bfk )$-module
map $\theta : \Qbar\to\Mbar$ induces a surjective map
  $$
\hat{\theta} : \Hom\,_{\Ufk (\bfk )}\bigl( \Ufk (\gfk ),
\chi^{\xi}_{\gamma} \otimes\Qbar^{\Fr}\bigr) \rightarrow
\Hom\,_{\Ufk (\bfk )} \bigl(\Ufk (\gfk ), \chi^{\xi}_{\gamma}
 \otimes\Mbar^{\Fr}\bigr) .
  $$
Hence, to prove (1), we can (and do) assume that $\Mbar$ is a
$\zi$-free module.

We first prove (1) for $m=1$.  Since $\Mbar$ is $\zi$-free (by
assumption), we can replace the ground ring $\zi$ by $\Qgil_{\xi}$.
For any $\dfk =(p_1,\cdots ,p_N)\in \{ 0,\cdots ,\ell -1\}^N$,
$N:=|\Delta_+|$, define
  $$
F^{\dfk} = F_{\beta_N}^{(p_N)} \cdots F^{(p_1)}_{\beta_1}.
   $$

By [Lu2, Lemma 8.5 and Th.~8.3] write
  $$
F_o\, E_i^{(\ell )} - E_i^{(\ell )} F_o = \sum_{0<m<\ell}
 E_i^{(m )} x_m+\sum_{\dfk\in\{ 0,\cdots ,\ell -1\}^N} c_{\dfk}\, F^{\dfk},
\tag{2}
  $$
for some $x_m\in\Ufk (\bfk^-)$ (in fact in the restricted quantized
enveloping algebra) and $c_{\dfk}\in\ucal^0_{\Qgil_{\xi}}$ (where
$\ucal^0_{\Qgil_{\xi}}\subset \Ufk^0_{\Qgil_{\xi}}$ is the
$\Qgil_{\xi}$-subalgebra generated by $\{ K_i^{\pm}; 1\le i\le n\}$).

Applying the anti-automorphism $S$ of $\Ufk (\gfk )$ to (2), we get
\begin{eqnarray*}
 (-K_i^{-1}E_i)^{(\ell )}  S(F_o)-S(F_o)(-K_i^{-1}E_i)^{(\ell )} 
&  =&
 \sum_{0<m<\ell}  S(x_m)(-K_i^{-1}E_i)^{(m)} \\
&&+\ \sum_{\dfk}
S(F^{\dfk})\, S(c_{\dfk}).
\end{eqnarray*}
Applying the above to a highest weight vector $v_+$ of
$V_{\xi}(2(\ell -1)\rho )$, we get
  $$
(-K_i^{-1} E_i)^{(\ell )} S(F_o) v_+ = \sum_{\dfk} S(F^{\dfk})\,
S(c_{\dfk})\, v_+ . \tag{3}
  $$

We next show that
  $$
E_i^{(m)}\, F_ov_+ = 0 \qquad\quad\mbox{for any } m>0. \tag{4}
  $$
Since $F_ov_+$ is a weight vector of weight 0, it suffices to show
that
  $$
F_i^{(m)}\, F_ov_+ = 0 \qquad\quad\mbox{for any } m>0: \tag{5}
  $$
For $0<m<\ell$, since $F_i^{(m)}\, F_o=0$, (5) follows in this case.
Further, $F_i^{(m)}$ commutes with $F_o$ (for any $m\geq 0$) as can
be 
seen from [Lu2, 5.8(c), Th.~8.3 and Lemma 8.5] by the weight
consideration.  Hence $F_i^{(m)}F_ov_+ = F_o\, F_i^{(m)}v_+$.  For
any
$1\leq i\leq n$, we can choose a reduced decomposition of $w_o$
starting in $s_i$ (and hence $\beta_1=\alpha_i$) resulting in the
expression $F_o=F_{\beta_N}^{(\ell -1)}  \cdots F_{\beta_2}^{(\ell -1)}\,
F_i^{(\ell -1)}$.  This gives  
$$
F_o F_i^{(m)} v_+=\left[\begin{array}{c} \ell +m-1\\ \ell -1
\end{array}\right]^{\xi}_{d_i}\, F_{\beta_N}^{(\ell -1)} \cdots F_{\beta_2}^{(\ell -1)}\,
F_i^{(\ell +m-1)} v_+ =0,  \quad\mbox{for}\,\, m\geq\ell ,
$$
which proves (5) and
hence (4). Substituting (4) in (3), we get
  \begin{equation}
\sum_{\dfk} S(F^{\dfk})\, S(c_{\dfk})v_+ = 0. \speqnu{6}  
  \end{equation}
Since $\{ F^{\dfk}v_+\}_{\dfk\in\{ 0,\cdots ,\ell -1\}^N}$ are
linearly independent, as the same is true already for the
Steinberg module $V_{\xi}((\ell -1)\rho )$ (cf.\ [Ku, Prop.~4.1]), from (6) we get
$\chi^\xi_{-\gamma} (S(c_{\dfk}))=0$, for all $\dfk$, i.e., 
  \begin{equation}
\chi^\xi_{\gamma} (c_{\dfk})=0, \quad\mbox{for all }\dfk . \speqnu{7}
  \end{equation}

By (2), 
  \begin{eqnarray*}
f(F_o E_i^{(\ell )}-E_i^{(\ell )} F_o) &= &\sum_{0<m<\ell}
E_i^{(m)}\cdot (f(x_m)) + \sum_{\dfk\in\{ 0,\cdots ,\ell -1\}^N}
c_{\dfk}\cdot (f(F^{\dfk})) \\[5pt]
 &=& \sum_{\dfk} c_{\dfk}\cdot (f(F^{\dfk}))\,,\,\hbox{since \, Im}\,f \,
\subset 
\chi^\xi_\gamma \otimes \bar{M}^{\Fr}, \\[5pt]
 &= &\sum_{\dfk} \chi^\xi_{\gamma} (c_{\dfk})\, f(F^{\dfk}),
\quad\mbox{since }c_{\dfk}\in\ucal^0_{\Qgil_{\xi}},\\[5pt]
 &= &0, \qquad\mbox{by (7)}.
  \end{eqnarray*}
This proves the identity (1) for $m=1$.
\smallbreak
We assume the validity of (1) for $m$ (by induction) and prove it for
$m$ replaced by $m+1$:  First of all
  \begin{eqnarray}
E_i^{(\ell m)} E_i^{(\ell )} &= \bbm \ell m+\ell\\[5pt] \ell
\ebm^{\xi}_{d_i}\, E_i^{(\ell m+\ell )}, \qquad\mbox{by [Lu2,
5.8(c)]} \speqnu{8}\\[5pt]
&= (m+1)\, E_i^{(\ell m+\ell )}, \qquad\mbox{by (3) of Lemma~3.7}.\nn
  \end{eqnarray}
Thus,
  \begin{eqnarray}
\enspace(m+1)\, f\bigl( F_o E_i^{(\ell m+\ell )}\bigr)&\hskip-7pt = \hskip-7pt & f\bigl( F_o
E_i^{(\ell m)} E_i^{(\ell )}\bigr) \speqnu{9}\\[5pt]
 &\hskip-7pt = \hskip-7pt &\bigl( E_i^{(\ell )}f\bigr) \bigl( F_o\, E_i^{(\ell m)}\bigr)\nn \\[5pt]
&\hskip-7pt = \hskip-7pt &  ( E_i^{(\ell )}f)\bigl( E_i^{(\ell m)}\, F_o\bigr) ,
\qquad\mbox{by induction},\nn \\[5pt]
 &\hskip-7pt = \hskip-7pt & E_i^{(\ell m)}\cdot \bigl( (E_i^{(\ell )}f)(F_o)\bigr)\nn \\[5pt]
 &\hskip-7pt = \hskip-7pt & E_i^{(\ell m)}\cdot \bigl( f(E_i^{(\ell )}\, F_o)\bigr),
\qquad\mbox{by the $m=1$ case},\nn \\[5pt]
 &\hskip-7pt = \hskip-7pt &f\bigl( E_i^{(\ell m)}\, E_i^{(\ell )}\, F_o\bigr)\nn \\[5pt]
 &\hskip-7pt = \hskip-7pt & (m+1)\, f\bigl( E_i^{(\ell m+\ell )}F_o\bigr), \qquad\mbox{by
(8)}. \nn
  \end{eqnarray}
Since $m+1$ is not a zero divisor in $\zi$ and (by assumption)
$\Mbar$ is $\zi$-free, we get the validity of (1) for $m+1$ (by
virtue of (9)).
  \hfill\qed

  \proclaim{Proposition} For any $\Mbar\in\Cbar (\bfk )${\rm ,} there exists
a functorial $\Ubar (\bfk ^- )$\/{\rm -}\/module map 
$$
\Frps_{\gamma} : H^0\bigl(\Xfk ,
\chi^{\xi}_{\gamma}\otimes\Mbar^{\Fr}\bigr)^{\Fr'}
\rightarrow H^0\bigl(\Xbar , \Mbar\bigr) ,
$$
defined by 
$$
(\Frps_{\gamma} f) (a) = f(F_o {\Fr'}(a))  \otimes  v_+ ,\tag{1}
$$
for $a \in \Ubar (\bfk ^- )$ and $f \in  H^0\bigl(\Xfk ,
\chi^{\xi}_{\gamma}\otimes\Mbar^{\Fr}\bigr)^{\Fr'}${\rm .}

Moreover{\rm ,} for any $ m\geq 0$ and  $f \in  H^0\bigl(\Xfk ,
\chi^{\xi}_{\gamma}\otimes\Mbar^{\Fr}\bigr)${\rm ,}
  $$
\Ebar_i^{(m)}\cdot \Bigl(\Frps_{\gamma} f\Bigr) = \Frps_{\gamma} 
({E}_i^{(m\ell )}\cdot f).\tag{2}
  $$
  \endproclaim

{\it Proof}.
As in the proof of Proposition 3.3, consider the diagram
 \smallbreak
{\ninepoint  \begin{eqnarray*}
\noalign{\vskip-12pt}
H^0\bigl( \Xfk ,\chi^{\xi}_{\gamma} 
\otimes\Mbar^{\Fr}\bigr)^{^{\Fr'}}
 &\hookrightarrow 
\Hom\,_{\Ufk(\bfk)}\bigl(\Ufk(\gfk ), \chi^{\xi}_{\gamma}
\otimes\Mbar^{\Fr}\bigr)^{^{\Fr'}} 
\simeq 
&\Hom\,_{\Ufk^0}\bigl( \Ufk(\bfk^-), \chi^{\xi}_{\gamma}
\otimes\Mbar^{\Fr}\bigr)^{^{\Fr'}}\\
&&\qquad\pch \downarrow \\
H^0(\Xbar , \Mbar ) &\hookrightarrow
\Hom\,_{\Ubar (\bfk )} (\Ubar (\gfk ),\Mbar ) \simeq &\Hom\,_{\Ubar^0} (
\Ubar (\bfk^-),\Mbar ),\\ 
\noalign{\vskip-24pt}
  \end{eqnarray*} 
}
\smallbreak\noindent 
where $\pch$ is as defined in \S4.1. 
By combining Lemma~4.3 and Proposition~4.4, we get  
  $$
\pch\Bigl( H^0\bigl(\Xfk , \chi^{\xi}_{\gamma}
 \otimes\Mbar^{\Fr}\bigr)^{^{\Fr'}}\Bigr) \subset
H^0(\Xbar ,\Mbar ).
  $$
So define $\Frps_{\gamma}$ as the restriction of $\pch$ to $H^0(\Xfk , 
\chi^{\xi}_{\gamma}\otimes\Mbar^{\Fr})^{^{\Fr'}}$. Since $\pch$ is\break
$\Ubar (\bfk^-)$-linear (by Lemma 4.3), so is $\Frps_{\gamma}$ and moreover (2)
follows from (1) of Proposition~4.4. \phantom{IMISSYOU}
  \hfill\qed

  \proclaim{Theorem}  
For any $\Mbar\in\Cbar (\bfk )${\rm ,} there exists a functorial $\Ubar
(\bfk^-)$\/{\rm -}\/module map {\rm (}\/for all $i\geq 0$\/{\rm )}
  $$
\Frps_{\gamma} : H^i\bigl(\Xfk ,
\chi^{\xi}_{\gamma}
\otimes\Mbar^{\Fr}\big)^{^{\Fr'}} \rightarrow H^i(\Xbar , \Mbar );
  $$
i.e.{\rm ,} the following diagram is commutative for any $\Ubar (\bfk
)$\/{\rm -}\/module map $\theta \!:\! \Mbar\! \to~\bar{N}${\rm :}
$$ 
    \begin{array}{ccc}
H^i\bigl(\Xfk , \chi^{\xi}_{\gamma}
\otimes\Mbar^{\Fr}\bigr)^{^{\Fr'}} &\stackrel{\Frps_{\gamma}}{\lrar}&
H^i(\Xbar,\Mbar )\\ [5pt]
 \Big\downarrow &&\Big\downarrow  \\
H^i\bigl(\Xfk , \chi^{\xi}_{\gamma}
\otimes\bar{N}^{\Fr}\bigr)^{^{\Fr'}} &\lower7pt\hbox{$\stackrel{\displaystyle\lrar}{\scriptstyle\Frps_{\gamma}}$}&
H^i(\Xbar ,\bar{N}),\\
  \end{array}  \eqno{\rm (D)}
   $$
where the vertical maps are the canonical maps induced from $\theta${\rm .}

Moreover{\rm ,} for any $m\geq 0$ and $f\in H^i(\Xfk 
, \chi^{\xi}_{\gamma}
\otimes\Mbar^{\Fr})${\rm , }
  $$
\Ebar_i^{(m)} \cdot (\Frps_{\gamma}\, f) = \Frps_{\gamma} (E_i^{(m\ell)}\cdot
f). \tag{1}
  $$\endproclaim

{\it Proof}.  Consider the
standard resolution in category $\bar{\Ccal} (\bfk )$:
  $$
0 \to \Mbar \to \Qbar_0 
\stackrel{\bar{\veps}_0\;}{\rightarrow} 
\Qbar_1 \stackrel{\bar{\veps}_1}{\rightarrow} \cdots . \tag{2}
  $$
Lifting $(2)$ by $\Fr$ and then tensoring 
 with $\chi^{\xi}_{\gamma}$, we get the
resolution in  $\Ccal (\bfk )$ :
  $$ 
0 \to \chi^{\xi}_{\gamma} \otimes\Mbar^{\Fr} \to
\hat{Q}_0 
\stackrel{\hat{\veps}_0\;}{\rightarrow} 
\hat{Q}_1 \stackrel{\hat{\veps}_1\;}{\rightarrow} \cdots , \tag{3}
  $$
where $\hat{Q}_k := \chi^{\xi}_{\gamma} \otimes \Qbar_k^{\Fr}$ and
$\hat{\veps}_k := \operatorname{Id} \otimes\bar{\veps}_k$.
By Proposition~4.6, we get the cochain map induced by the $\Ubar
(\bfk^-)$-module maps $\Frps_{\gamma}$:
  $$ \begin{array}{ccc}
H^0(\Xfk ,
\hat{Q}_0)^{\Fr'} &\lrar&  H^0(\Xfk , \hat{Q}_1)^{\Fr'}
{\longrightarrow}
\cdots \\[5pt]
\Big\downarrow {\scriptstyle{\Frps_{\gamma}}} &&\Big\downarrow{\scriptstyle\Frps_{\gamma}}  \\[10pt]
H^0(\Xbar ,
\bar{Q}_0)&\lrar& H^0(\Xbar , \bar{Q}_1) \longrightarrow \cdots \,\cdot
  \end{array} 
\eqno{(\ast)} 
$$ 

By the next lemma, for any $p\geq 0$ and $i>0$, $H^i(\Xfk , \hat{Q}_p)=0$.  
Hence the $i^{\rm th}$ cohomology of the
top cochain complex is equal to $H^i(\Xfk 
, \chi^{\xi}_{\gamma}
\otimes\Mbar^{\Fr})^{\Fr'}$, whereas the $i^{\rm th}$ cohomology of the bottom
cochain complex is $H^i(\Xbar ,\Mbar )$ (cf.\
[H, Prop.~1.2A, Chap. III]).  So, we define the $\Ubar
(\bfk^-)$-module map
  $$
\Frps_{\gamma}: H^i\bigl( \Xfk ,
\chi^{\xi}_{\gamma} \otimes
\Mbar^{\Fr}\bigr)^{\Fr'} \rightarrow H^i(\Xbar , \Mbar )
  $$
as the induced map in cohomology from ($\ast$).

Commutativity of   diagram (D) follows from the functoriality of
all the constructions involved and moreover (1) follows from (2) of
Proposition~4.6.  So the theorem is proved modulo the following
lemma.
  \hfill\qed

  \proclaim{Lemma}  For any $\Mbar\in\Cbar^0$ and $\lam\in X$
$$ H^i\Bigl( \Xfk , \chi^{\xi}_{\lam} \otimes
\bigl( H^0 (\Ubar (\bfk )/\Ubar^0 , \Mbar)\bigr)^{\Fr}\Bigr) =0,
$$
for all $i>0${\rm .}
\endproclaim
 
{\it Proof}.
The proof is similar to the proof of [APW, Th.~5.4].
By definition of the category $\Cbar^0$, $\Mbar = \bigoplus_{\mu\in
X} \Mbar_{\mu}$.  Since $H^*$ commutes with (possibly infinite)
direct sums (cf.\ [APW, Th.~1.31]), we can assume that $\Mbar
=\Mbar_{\mu}$.  Since $H^0(\Ubar (\bfk )/\Ubar^0, -)$ is an exact
functor (cf.\ Proposition 1.5) and any $\zi$-module $N$ admits
a free resolution for some $d\geq 0$ (since $\zi$ has finite global
homological dimension):
$$0 \to F_d \to \cdots \to F_1 \to F_0 \to N, 
$$
we can assume that $\Mbar_{\mu}$ is $\zi$-free of rank 1, i.e.,
$\Mbar_{\mu} = \bar{\chi}_{\mu}$.  By a result of
Cline-Parshall-Scott (cf.\ [CPS, Proposition 5.5], [APW, Lemma 5.3]), 
there is a $\Ubar (\bfk
)$-module isomorphism:
$$
H^0(\Ubar (\bfk )/\Ubar^0, \bar{\chi}_{\mu}) \simeq 
 {\displaystyle\lim_{{\textstyle \longrightarrow}\atop{m\geq 0}}}
 \, \, H^0(\Xbar ,
\bar{\chi}_{-m\rho}) \otimes \bar{\chi}_{-m\rho +\mu}, 
$$
where the right side is a directed union.  Since the cohomology
commutes with directed unions (cf.\ [APW, Proof of Th.~5.4]), to
prove the lemma, it suffices to show that
 $$
H^i\Bigl( \Xfk , H^0\bigl( \Xbar , \bar{\chi}_{-m\rho}\bigr)^{\Fr}\otimes
\chi^{\xi}_{-\ell m\rho +\ell\mu +\lam}\Bigr) =0,\,
 \forall i>0, \,m\gg 0. 
$$
 (We have used here the fact that $\Fr$
commutes with the coproduct.)

Now, by [APW, Prop.~2.19 (ii)] (since $H^0(\Xbar ,
\bar{\chi}_{-m\rho})$ is ${\Zgil}_\xi$-free; cf.\ [APW, Corollary 3.3(i)] for 
the corresponding result in the quantum case), 
  \begin{eqnarray*}
&&\hskip-1.5in H^i\bigl( \Xfk ,
H^0(\Xbar , \bar{\chi}_{-m\rho})^{\Fr} \otimes
\chi^{\xi}_{-\ell m\rho +\ell\mu +\lam}\bigr)  \\ 
& \simeq & H^0(\Xbar ,
\bar{\chi}_{-m\rho})^{\Fr} \otimes H^i (\Xfk ,
\chi^{\xi}_{-\ell m\rho +\ell\mu +\lam})\\ 
 & = & 0, \,\, \hbox{for} \,m\gg 0\,,
  \end{eqnarray*}
by the following quantized analogue of the Serre vanishing theorem.  

{\it For any $\lambda \in X^+ , \mu \in X$ and} $ i > 0$,
$$
 H^i(\Xfk , \chi^\xi_{-m(\lambda + \rho)+\mu} ) = 0 \,, \,\hbox{\it for\,
all}\,\, m \gg 0. \tag{1}   
$$
The ring $\zi$ has projective dimension one 
(cf.\ [Mi, Lemma 1.5]). Hence 
(1) follows by the base change [APW, (8) of \S3.6] and [AW, Th.~2.6].

This completes the proof of the lemma, thereby completing the proof
of Theorem~4.7.
  \hfill\qed \pagebreak

   \proclaim{Lemma}  For any $M\in \Ccal (\bfk )$ such that $M$ is
$\zi$\/{\rm -}\/flat and $\sigma\in H^0(\Xfk , M)${\rm ,}
there exists a functorial map{\rm ,} for any $N \in \Ccal (\bfk )${\rm ,}
$$
m_{\sigma}: H^0(\Xfk ,N) \rightarrow H^0(\Xfk , M\otimes N)
$$
given by
$$
(m_{\sigma}f)\, a = \sum_i\, \sigma (a'_i)\otimes f(a_i''), \quad
\mbox{for } a\in\Ufk (\gfk ), \tag{1}
$$
  where $\Delta a = \sum_i\, a_i'\otimes a_i''${\rm .} 

Moreover{\rm ,} $m_{\sigma}$ gives rise to a functorial map {\rm (}\/again denoted
by\/{\rm )}
\smallbreak\noindent
{\rm (2)} \hfill$
m_{\sigma}: H^i(\Xfk , N) \rightarrow H^i(\Xfk
, M\otimes N). 
$\hfill
  \endproclaim  
\vglue6pt
 {\it Proof}.  It is easy to see that (1) defines a map
  $$
m_{\sigma} : \Hom\,_{\Ufk (\bfk )} (\Ufk (\gfk ),N) \to \Hom\,_{\Ufk
(\bfk )}(\Ufk (\gfk ), M\otimes N).
  $$
Moreover, $m_{\sigma}(f)\in H^0(\Xfk , M\otimes
N)$ for any $f\in H^0(\Xfk , N)$ (as is easy to
see).  The existence of $m_{\sigma}$ at the higher cohomology follows
from the consideration of the standard resolution of $N$ in $\Ccal (\bfk )$:
$0 \to N \to Q_0\to Q_1 \to \cdots$.
(Observe that, since $M$ is
$\zi$-flat, by [APW, Prop.~2.16(i)],
 $0\to M\otimes N\to M\otimes Q_0 \to M\otimes Q_1 \to
\cdots$ is the standard resolution of $M\otimes N$.)
  \hfill\qed

\proclaim{Lemma}  Let  $ \zi v_+$ be the  highest weight subspace of 
$V_{\xi}(2(\ell -1)\rho ).$ Then  $S(F_o)\, \zi v_+$ is a $\zi$\/{\rm -}\/module 
direct summand of $V_{\xi}(2(\ell -1)\rho ).$
\endproclaim

\demo{Proof} It suffices to construct $f\in  \hbox{Hom}_{\zi}
\,(V_{\xi}(2(\ell -1)\rho ), \zi)$ such that $f(S(F_o) v_+)$ is an invertible
element of $\zi$. Consider the  $\Ufk (\gfk )$-module homomorphism 
$$\delta: V_{\xi}(2(\ell -1)\rho ) \to V_{\xi}((\ell -1)\rho ) \otimes
V_{\xi}((\ell -1)\rho ),$$
 taking $v_+ \mapsto w_+\otimes w_+$ , where $w_+$ is a primitive highest 
weight vector of $V_{\xi}((\ell -1)\rho )$.
Since $S$ is an (anti)automorphism which keeps the restricted enveloping
algebra stable, by the weight consideration, $S(F_o)= F_o a$,
for an invertible element $a \in u^0$, where $u^0$ is the $\Bbb Z_\xi$-subalgebra of $\frak U^0$ generated by $\{k_i;1\le i\le n\}$.
Thus 
$$S(F_o) v_+ = 
F_o a v_+ = x F_ov_+\,,$$
for an invertible element $x\in \zi$. Write
$$
\delta (S(F_o) v_+) = x w_+\otimes F_o w_+ + v\,, \tag{1}
$$
for $v\in  V_{\xi}((\ell -1)\rho )^+ \otimes
V_{\xi}((\ell -1)\rho )$, where $ V_{\xi}((\ell -1)\rho )^+$ (resp.\ $V_{\xi}((\ell -1)\rho )^-$) is the sum of all the weight spaces of $ 
V_{\xi}((\ell -1)\rho )$ except the highest (resp.\ lowest) weight space. Now,
 by
[Ku, Prop.\ 4.1] and [Lu2, Th.\ 8.3], $F_ow_+$ is a 
primitive vector in  
$V_{\xi}((\ell -1)\rho )$. The decomposition
\begin{eqnarray*}
&& V_{\xi}((\ell -1)\rho ) \otimes  V_{\xi}((\ell -1)\rho )\\ [5pt]
&&\qquad  =  
\biggl(V_{\xi}((\ell -1)\rho )^+ \otimes  V_{\xi}((\ell -1)\rho ) +  
V_{\xi}((\ell -1)\rho ) \otimes
 V_{\xi}((\ell -1)\rho )^-\biggr)\\ [5pt]
&&\qquad\qquad \oplus\  \zi (w_+ \otimes F_o w_+)
\end{eqnarray*}  
gives rise to the map $\tilde{f}: V_{\xi}((\ell -1)\rho ) \otimes 
V_{\xi}((\ell -1)\rho ) \to \zi$ by projecting on the last factor. Finally, let 
$f$ be the linear form $\tilde{f}\circ \delta : V_{\xi}(2(\ell -1)\rho ) \to \zi$. Then, by (1), $ f(S(F_o) v_+)= x$. Hence $\zi S(F_o)v_+$ is a $\zi$-module
direct summand in $V_{\xi}(2(\ell -1)\rho )$. 
\enddemo

Decompose 
  $$
V_{\xi}(2(\ell -1)\rho ) = S(F_o)\, \zi v_+ \oplus M \,,
  $$
where
$M$ is a weight subspace.  Define $\hat\sigma_o\in V_{\xi}
(2(\ell -1)\rho )^*$ by $\hat{\sigma}_o(S(F_o)v_+)=1$ and
$\hat{\sigma}_{o_{|M}}\equiv 0$.

  \proclaim{Proposition}  For any $\Mbar\in\Cbar (\bfk )${\rm ,} the composite
$$
H^i\bigl(\Xbar ,\Mbar \bigr)
 \stackrel{\Fr^*}{\longrightarrow} 
H^i\bigl(\Xfk ,\Mbar^{\Fr}\bigr) 
\stackrel{m_{\sigma_o}}{\longrightarrow} H^i\bigl(\Xfk , 
\chi^{\xi}_{\gamma} \otimes \Mbar^{\Fr}\bigr) 
 \stackrel{\Frps_{\gamma}}{\longrightarrow}  H^i\bigl(\Xbar , \Mbar \bigr)
\tag{1}
$$
is the identity map{\rm ,} where $\sigma_o\in H^0\bigl( \Xfk , 
\chi^{\xi}_{\gamma}\bigr)$ is given by
$\sigma_o:=\beta (\hat{\sigma}_o)${\rm ,} and $\beta$ is as defined in
Section~{\rm 1.}
  \endproclaim

  \demo{Proof} From the functoriality of all the maps involved, it
suffices to prove the lemma for $H^0$ (cf.\ the argument in the proof of 
Corollary 3.9).  Take
 $ f\in H^0(\Xbar ,\Mbar ), \bar{y}\in\Ubar
(\nfk^-)$
  and write $\Delta (\Fr' \bar{y}) = \sum_i y'_i\otimes y_i''$.
Also write $\Delta (F_o)=\sum_j F'_j\otimes F_j''$.  Then
  \begin{eqnarray*}
&&\hskip-18pt (\Frps_{\gamma}\, m_{\sigma_o} \Fr^* f)\, \bar{y} \\ [5pt]
&&\qquad= (m_{\sigma_o}\Fr^* f) (F_o \Fr'(\bar{y})) \otimes v_+ \\ [5pt]
&&\qquad= \sum_{i,j} \sigma_o(F_j' y_i')\otimes f(\Fr (F_j'')\, \Fr 
(y_i'')) \otimes v_+\\ [5pt]
&&\qquad= \sum_i \sigma_o(F_o y_i')\otimes f(\Fr (y_i''))\otimes v_+, 
  \mbox{  since }\Fr (F_j'')=0 \mbox{ unless } F_j''\in\Ufk^0 \\ [5pt]
&&\qquad= \sigma_o(F_o)\otimes f(\Fr\, \Fr' \bar{y})\otimes v_+,
  \mbox{  since }\sigma_o(F_oy_i')=0 \mbox{ unless }y'_i\in\Ufk^0 \\ [5pt]
&&\qquad= v_+^* \otimes f(\bar{y})\otimes v_+ \\ [5pt]
&&\qquad= f(\bar{y}).
  \end{eqnarray*}
This proves the proposition.
  \enddemo\pagebreak

\numbereddemo{{R}emark} (a) By the same argument as that in the proof of Proposition 
 4.6, we obtain the fact that  for any $\Mbar\in\Cbar (\bfk )$, there exists
a functorial $\Ubar (\bfk ^- )$-module map 
 $$
\Theta : H^0\bigl(\Xfk ,
\Mbar^{\Fr}\bigr)^{\Fr'}
\rightarrow H^0\bigl(\Xbar , \Mbar\bigr) ,
  $$ 
defined by 
$$
(\Theta f) (a) = f( {\Fr'}(a))  ,\,\hbox{for}\, a \in \Ubar (\bfk ^- ).
$$
From the following lemma, we see that the map $$\Fr^{\prime*} : 
 H^0\bigl(\Xfk ,
\bar{M}^{\Fr}\bigr)^{\Fr'}
\rightarrow H^0\bigl(\Xbar , \Mbar\bigr)$$ defined in Proposition~3.3 
coincides with $\Theta$. (Observe  however that in  Proposition~3.3, 
the map $\Fr^{\prime*}$ was defined for an arbitrary $M \in \Ccal (\bfk)$.)
\smallbreak
(b) It is not clear if the composite $\Frps_{\gamma} \circ\ m_{\sigma_o}$ of 
the last two maps in the above proposition is the map $\Frps$ of Proposition 
 3.3.   
\enddemo

  \proclaim{Lemma}  For any $\Mbar\in\Cbar (\bfk )${\rm ,} the map
  $$
\Frps : H^0\bigl(\Xfk ,\Mbar^{\Fr}
\bigr) \to H^0(\Xbar , \Mbar )
  $$
satisfies $\Frps\, f=0$ if $f$ is a weight vector of weight
$\lam\notin \ell X$.
  \endproclaim

  \demo{Proof}
Choose an $i$ such that $\lam (h_i) = \lam_0+\ell\lam_1$,
$0<\lam_0<\ell$.  For  $a=\left[\begin{array}{c} K_i; -\ell \lambda_1\\ \lam_0\end{array}\right]$ and
any
$\bar{y} = \bar{F}_{i_{r}}^{(m_r)} \cdots
\bar{F}_{i_1}^{(m_1)}$, 
  \begin{eqnarray}
 f\bigl( (\Fr'\bar{y}) a\bigr) &=& f\Bigl(
F_{i_{r}}^{(\ell m_r)} \cdots F_{i_1}^{(\ell m_1)}
\bbm K_i;  -\ell \lambda_1\\ \lam_0\ebm\Bigr) \speqnu{1}\\
&=& f\Bigl(
F_{i_{r}}^{(\ell m_r)}\cdots F_{i_1}^{(\ell m_1)}\Bigr) .\nn
  \end{eqnarray}
Also
  $$
F_{i_{r}}^{(\ell m_r)} \cdots F_{i_1}^{(\ell m_1)} 
\bbm K_i;  -\ell \lambda_1\\ \lam_0\ebm 
= \bbm K_i; -\ell \lambda_1+ \sum_k \ell m_k\, a_{i_k}\\ \lam_0\ebm 
F_{i_{r}}^{(\ell m_r)} \cdots F_{i_1}^{(\ell m_1)}, 
  $$
by [Lu2, \S6.5], where $a_{i_k} := \alpha_{i_k}(h_i).$
So
 $$
f\bigl( (\Fr' \bar{y}) a\bigr) =0 , \qquad\mbox{since \, Image\,}
f\subset \Mbar^{\Fr} .
 \tag{2}
$$

Comparing (1) and (2), we get 
$$(\Frps
f)(\bar{y}) = f\bigl( F_{i_{r}}^{(\ell m_r)} 
\cdots F_{i_1}^{(\ell m_1)}\bigr) =0\,.
$$  This proves the lemma.
  \enddemo

\numbereddemo{{R}emark} \,\,(a) Observe that to prove Theorem~4.7, 
we needed the quantized version of Serre vanishing (cf.\ (1) of Lemma 4.8),
whereas the corresponding results for $\Fr^*$ and $\Fr^{\prime*}$ (as in 
Theorems 2.3 and 3.8 respectively) did not require this. This is due to the 
fact that  Proposition~4.6  (and hence Theorem 4.7) is available only 
for $\Mbar^{\Fr}$, where 
 $\Mbar\in\Cbar (\bfk )$. In fact, the analogue of Proposition~4.6  for an 
arbitrary $M \in \Ccal (\bfk)$ (as in Proposition 3.3) is {\it false}, as 
can be seen already   in the case of $\gfk = sl(2)$. 
\smallbreak
(b) $\Fr^{\prime*}_\gamma$ is a quantization of the stronger 
{\it Frobenius $\Lcal(\gamma)$-splitting } of the flag variety $G/B$ (proved by 
Ramanan-Ramanathan) 
for the homogeneous line bundle   
$\Lcal (\gamma)$ on $G/B$ associated to the character $\gamma$ of $B$
(see \S 6 for further details). 
\enddemo

   \section{The Kempf vanishing theorem}

In this section we assume that $\ell = p$ is an odd prime.  We further
assume that $p\neq 3$ if $G_2$ is a factor of $\gfk$.  Let
$\Fgil_{p}$ be the prime field with $p$ elements and let
$k$ be any field containing  $\Fgil_p$.
 Let $G$ be the connected
simply-connected 
semisimple algebraic group defined and split over $k$ corresponding to $\gfk$ and let $B$ be its
Borel subgroup defined over $k$ (corresponding to $\bfk$). 
Consider the base change $c: \zi \to \Fgil_p\subset k$ which takes
$\xi \mapsto 1$.

 For $\lam\in X$, we denote
by $\Lcal (\lam )$ the line bundle on the flag variety $G/B$
corresponding to the character $e^{\lam}$ of $B$.  
More generally, for any $M\in\Ccal (\bfk )$, we denote by $\Lcal (M)$
the homogeneous vector bundle on $G/B$ associated to the\break $B$-module
induced by the $\Ufk_k (\bfk)$-module $M_k :=  M\otimes_{\zi}k $, where 
 $\Ufk_k (\bfk) :=  \Ufk (\bfk)\otimes_{\zi}k$ (cf.\ [Lu2, \S8.15] and
[CPS, Th.\ 9.4]).

  We recall the following result due to [APW,
Prop.~3.7].

  \proclaim{Proposition}  For any $M\in\Ccal (\bfk )${\rm ,} there exists a
canonical isomorphism\/{\rm :}
  \begin{equation}
H^i(\Ufk_k(\gfk )/\Ufk_k(\bfk ), M_k) \simeq H^i(G/B, \Lcal (M)) 
\qquad\mbox{for all } i\geq 0. 
  \end{equation}
Similarly{\rm ,} for $\Mbar\in\Cbar (\bfk )${\rm ,} there exists a canonical
isomorphism\/{\rm :}
  \begin{equation}
H^i(\Ubar_k(\gfk )/\Ubar_k(\bfk ), \Mbar_k) \simeq H^i(G/B, \Lcal
(\Mbar))  \qquad\mbox{for all } i\geq 0.
  \end{equation}
  \endproclaim

As a corollary of the strong Frobenius splitting as in Section~4  (cf.\
Theorem 4.7 and Proposition 4.11), we obtain the following Kempf vanishing
theorem [K] (if we use the usual Serre vanishing theorem).

  \proclaim{Theorem}  For any $\lam\in X$ such that $\lam +\rho \in X^+${\rm ,} 
  $$
H^i(G/B, \Lcal (-\lam ))=0 \qquad\mbox{ for all }i>0.
  $$\endproclaim

\phantom{goodsong}
\vglue-36pt
{\it Proof}.  The constructions and results of Section~4 are
compatible under base change.  Hence, by Proposition~4.11, the map
  $$
m_{\sigma_o}\Fr^* : H^i\bigl( \Ubar_k(\gfk )/\Ubar_k(\bfk ),
(\bar{\chi}_{-\lam})_k \bigr) \to H^i\bigl(\Ufk_k(\gfk )/\Ufk_k(\bfk ),
(\chi^{\xi}_{-2(p -1)\rho -p\lam})_k\bigr)
  $$
is injective.

Applying Proposition~5.1, we get an injective map
  $$
H^i(G/B, \Lcal (-\lam )) \hookrightarrow H^i(G/B, \Lcal (-2(p
-1)\rho -p\lam )).
 $$
Iterating $m$-times, we get an injection 
 $$
H^i(G/B, \Lcal (-\lam )) \hookrightarrow H^i(G/B, \Lcal (-p^m(\lambda
+2\rho ) + 2 \rho  )).
 $$
Now, using the Serre vanishing theorem for the cohomology of ample line 
bundles on $G/B$ [H, Chap. III, Prop.~5.3], we obtain that
$H^i(G/B, \Lcal (-\lam ))=0$, for all $i>0$.  
  \hfill\qed

 \vglue-6pt
  \section{Sheafification: Frobenius splitting of $G/B$  and Schubert
varieties}\label{frobsplitandschubi}
 \vglue-6pt

We follow the same notation and conventions as in Section~5.  In particular,
$\ell=p$ is an odd prime, and $p>3$ if $\km$ has a simple component
of type ${\tt G}_2$. Let $k$ be an algebraically closed field of
characteristic $p$. 

 \numbereddemo{Definition}  The {\it absolute Frobenius morphism } of a scheme $X$ over $k$
is the identity map on the underlying point space and arises to the $p^{\rm th}$ power locally 
on the  functions.  Observe that the absolute Frobenius morphism  is {\it not} 
a morphism of $k$-schemes. To remedy this, let $X'$ be the scheme with the 
same underlying topological space as that of $X$ and the same structure sheaf 
$\Ocal_X$ of rings, only the scalar multiplication of $k$ on   $\Ocal_{X'}$
is twisted as:
$$ z \odot f = z^p f\,, \,\,\hbox{for}\, z \in k \,\,\hbox{and} \, f \in 
\Ocal_{X'}\,. $$
Thus we get a morphism of $k$-schemes $F: X' \to X $, which at the point set 
level is the identity map and at the sheaf level corresponds to the morphism
$\Ocal_X \to \Ocal_{X'}, f \mapsto f^p$.

Following Mehta-Ramanathan [MR1],  a scheme $X$  is called {\it Frobenius split}
if the homomorphism $\Ocal_X \to F_*(\Ocal_{X'})$ of $\Ocal_X$-modules is 
split. 
Then clearly an $\Ocal_X$-module homomorphism $\sigma:  F_*(\Ocal_{X'}) \to 
\Ocal_X $ is a splitting of  $\Ocal_X \to F_*(\Ocal_{X'})$ if and only if $\sigma (1)
=1$. A splitting $\sigma$ is said to {\it split a closed 
subscheme $Y \subset X$
 compatibly}  if $\sigma (F_* (\Ical_{Y'} )) \subset \Ical_{Y} $, where
$\Ical_{Y}$ is the ideal sheaf of $Y$. 

Let $(\Dcal , \phi)$ be a line bundle on $X$  with a  section $\phi$. Then, 
pulling back via $F$, we get the line bundle (denoted) $\Dcal'$, and 
 the section 
$\phi$ clearly gives rise to an  $\Ocal_X$-linear morphism  $ \bar{\phi}: 
F_*(\Ocal_{X'}) \to 
F_*(\Dcal')$.
Following Ramanan-Ramanathan [RR, \S 2],  a scheme  $X$ is  called {\it Frobenius
$(\Dcal , \phi)$-split} (or less precisely {\it Frobenius
$\Dcal$\/{\rm -}\/split}) if there exists an  $\Ocal_X$-linear morphism 
$$\sigma^\Dcal: F_*(\Dcal') \to \Ocal_X\,, $$
such that $
\sigma^\Dcal \circ  \bar{\phi}$ is a splitting. A closed subscheme 
$Y \subset X$ is called {\it compatibly $(\Dcal , \phi)$-split} (or {\it 
compatibly $\Dcal$\/{\rm -}\/split}) if   $
\sigma^\Dcal \circ  \bar{\phi}$ compatibly  splits $Y$ and moreover on no
irreducible component of $Y$, $\phi$ is identically zero. 

Now, we come to the corresponding `local version'.
Let $K$  be a  $k$-algebra.
Recall (cf. [M], [Ka, \S4.3]) that a {\it Frobenius-linear  endomorphism} of $K$ 
 is a map
$\sigma: K\rightarrow K$ such that for all $f, g \in K :$

 a)  additivity: $\sigma(f+g)=\sigma(f)+\sigma(g)$, and

\smallbreak b)    $\sigma(f^p g)=f \sigma(g)$.

Observe that by (b), $\sigma$   cannot be  $k$-linear. The $k$-space of all the 
 Frobenius-linear endomorphisms of $K$ is denoted by End$_F(K)$.

A  Frobenius-linear endomorphism $\sigma$ is called a {\it splitting}
if  $\sigma(f^p )=f \,,$ for all $f \in K$. Hence a Frobenius-linear endomorphism $\sigma$ is a splitting if and only if $\sigma (1)
=1.$ Let
$I$ be an ideal of
$K$. Then a splitting
$\sigma$ is said to {\it compatibly split } $I$ if and only if $\sigma I \subset I$.

For $K = \oplus_{\lam \in \Zgil^n_+}\,K_\lambda $  a graded $k$-algebra,
 a Frobenius-linear endomorphism $\sigma$ is called {\it  graded} if for all
$\lambda \in \Zgil^n_+$,
$\sigma (K_{p\lambda}) \subset K_\lambda $ and $\sigma (K_{\lambda}) =0 ,$  if
$p$ does not divide $\lambda$. 

For a smooth scheme $X$  over $k$ with an affine open cover $\{X_j\}$, giving a Frobenius 
splitting of $X$ is equivalent to giving splittings 
$\sigma_j$ of the affine coordinate rings $k[X_j]$ and  $\sigma_{\{j,j'\}}$ of
$k[X_j \cap X_{j'}]$ (for unordered pairs $\{j,j'\}$) such that the 
following compatibility is satisfied:
\begin{equation}
 \sigma_j (f)_{|X_j \cap X_{j'}} = \sigma_{\{j,j'\}} (f_{|X_j \cap X_{j'}})\,, 
\,\,\hbox{for  all} \, f \in k[X_j].
\end{equation}
Moreover, under this correspondence, a closed subscheme $Y \subset X$ is 
compatibly split if and only if 
\begin{equation}
 \sigma_j (\Ical_{X_j \cap Y})  \subset  \Ical_{X_j \cap Y}
\,,\,\, \hbox{for  all} \, j\,,
\end{equation}
where $\Ical_{X_j\cap Y} \subset k[X_j]$ is the ideal of $X_j\cap Y$.  

For a smooth scheme $X$ with an affine open cover $\{X_j\}$ and a line 
bundle together with a  section  
 $(\Dcal , \phi)$ on $X$ as above, giving a Frobenius $\Dcal$-splitting of $X$ is 
equivalent to giving
 additive maps 
$\sigma_j^\Dcal : \Gamma (X_j, \Dcal) \to k[X_j]$ and  
\vglue2pt
\centerline{${\displaystyle\sigma_{\{j,j'\}}^\Dcal : \Gamma (X_j \cap X_{j'}, \Dcal) \to 
k[X_j \cap X_{j'}]}$}
\vglue2pt \noindent satisfying the 
following three properties (3)--(5):
\begin{equation}
\sigma_j^\Dcal (f^p s) = f \sigma_j^\Dcal (s)\,,\,\hbox{for}\, f\in  k[X_j]\,
\hbox{and} \, s\in  \Gamma (X_j, \Dcal)\,
\end{equation}
and similarly for $\sigma_{\{j,j'\}}^\Dcal$, 
\vglue-19pt
\begin{equation} 
\sigma_j^\Dcal(\phi_{|X_j}) = 1\,,
\end{equation}
and for any pair $\{j,j'\}$,
\begin{equation}
   \sigma_j^\Dcal (s)_{|X_j \cap X_{j'}} = 
\sigma_{\{j,j'\}}^\Dcal (s_{|X_j \cap X_{j'}})\,, 
\,\,\hbox{for \, all} \, s \in \Gamma (X_j, \Dcal) ,
\end{equation}
where  $ \Gamma (X_j, \Dcal)$ denotes the $k[X_j]$-module of all the regular 
sections of $\Dcal$ on~$X_j$.

Let $H$ be an algebraic group over $k$ and let $H$ act  on a $k$-algebra $K$  algebraically
via $k$-algebra automorphisms. Then we define an $H$-action on
End$_F(K)$ by 
\vglue2pt
\centerline{${\displaystyle ( h * \sigma ) (a) = h\cdot \sigma(h^{-1}\cdot a), \,\hbox{for}\, h \in H,
\sigma \in \, \hbox{End}_F(K),\, \hbox{and} \, a \in K\,.}$ }
\vglue2pt

Let $B$ be a Borel subgroup of a connected, simply-connected, semisimple 
algebraic group
$G$ over $k$ and let $T\subset B$ be a maximal torus.  Let $B$ act on
$K$ algebraically via $k$-algebra automorphisms.  Then a splitting
$\sigma$ of $K$ is called $B$-{\it canonical} if 
\vskip1ex

c)  $t*\sigma =\sigma$, for all $t\in T$, and

\vskip1ex

d)  for each simple root $\alpha_i$, there exist Frobenius-linear
endomorphisms $\{ \sigma_{m,\alpha_i}\}_{0\leq m<p}$ of $K$ such that
  $$
(x_{\alpha_i}(z)*\sigma )(a) = \sum_{m=0}^{p-1} \sigma_{m,\alpha_i}
(z^ma), \;\;\mbox{for all $z\in k$ and $a\in K$},
  $$
where $x_{\alpha_i}(z)$ is the one-parameter subgroup of $B$
corresponding to the simple root $\alpha_i$.

We can easily sheafify and extend the notion of $B$-canonical splitting for a 
variety $X$ with an action of $B$ (cf.\ [Ka, p.\ 42], [M]).
\enddemo

 Mehta-Ramanathan [MR1] have shown that the flag variety $ G/P$ (for  any 
parabolic subgroup $P\subset G$) is
Frobenius split, and the Schubert subvarieties $X(w)_P\subset G/P$ are
compatibly split. Moreover, this splitting is $B$-canonical. The aim of this section is to show
that the map $\Frps$ constructed in Section~3 induces such a splitting in the case 
$P=B$. The general case will be handled in Section~9. 

Recall from Lemma 2.2 (b) and Proposition 3.3 
that, for any  $\lam  \in X$,  there are maps 
$$
\Fr^* : H^0 \bigl(\Xbar, \bar{\chi}_{\lam}\bigr)^{ \Fr}
\rightarrow H^0\bigl(\Xfk, \chi^\xi_ {p\lam}\bigr),\quad
\Frps : H^0\bigl(\Xfk,  \chi^\xi_ {p\lam}\bigr)^{\Fr'_{\bfk^-}} \rightarrow
H^0\bigl(\Xbar, \bar{\chi}_{\lam}\bigr).
$$
(In fact $\Frps$ is defined for any $\mu \in X$, not only for 
$p\lam$, but unless $\mu$ is divisible by $p$ it is identically zero 
by Remark 3.10 (b).)

For a vector space $V$ over $k$, by $V^{[1]}$ we mean the same additive group,
but the scalar multiplication is twisted as $z \odot v := z^p v$. 
 Let $G(\Fgil_p)$ be the connected,
simply-connected,
semisimple algebraic group defined and split over $\Fgil_p$ corresponding to $\gfk$ and let 
$B( \Fgil_p)$ be its
Borel subgroup (corresponding to $\bfk$) defined over $\Fgil_p$ (cf.\ \S 1 for the notation $\frak g$ and $\frak b$). We denote by 
$G$ and $B$ the corresponding $k$-rational points. 
Consider the base change $c: \zi \to 
\Fgil_p$ which takes
$\xi \mapsto 1$. Then, as in Proposition~5.1, we have canonical $k$-linear 
isomorphisms for any $\lam \in X$,
\begin{eqnarray*}
\theta: H^0(\Ufk_k(\gfk )/\Ufk_k(\bfk ), (\chi^\xi_{\lam})_k)& \simeq& H^0(G/B, 
\Lcal 
(\lam))^{[1]} \,,\,\, \hbox{and}\\
 \bar{\theta}: H^0(\Ubar_k(\gfk )/\Ubar_k(\bfk ), 
(\bar{\chi}_\lam )_k) &\simeq&
H^0(G/B, \Lcal
(\lam))\,,
\end{eqnarray*}
obtained by first taking the base field $\Fgil_p$ in Proposition~5.1  and 
then extending $k$-linearly (with respect to the twisted $k$-linear structure
on $H^0(G/B, 
\Lcal 
(\lam))^{[1]}$ in the first case).
Thus, under the above identifications, the maps $\Fr^*$ and $\Frps$ 
(after the base change) become $k$-linear maps (for any $\lam \in X^+$):
 $$\Fcal_\lam : H^0(G/B, 
\Lcal 
(-\lam))^{ \Fr} \to  H^0(G/B, 
\Lcal 
(-p\lam))^{[1]}\,,\,\,\hbox{and}$$
$$\Fcal'_\lam:  \bigl(H^0(G/B, 
\Lcal 
(-p\lam))^{[1]}\bigr)^{\Fr'_{\bfk^-}} \to H^0(G/B, 
\Lcal 
(-\lam))\,.$$

We have shown in [KL, Th.\ 1] that $\Fcal_\lam$ is the  map $s\mapsto s^p$
sending
a section to its $p^{\rm th}$ power, and $\Fcal'_\lam$ provides a splitting of 
this map by Corollary~3.9.

We consider the $ X^{+}$--graded algebra $K:=\bigoplus_{\mu\in X^{+}}
H^0\bigl(G/B, \Lcal (-\mu)\bigr)$\break (under the multiplication of sections).
Then $G$ acts algebraically on $K$ 
via $k$-algebra automorphisms; in particular, so does $B$.

Let $\Fcal': K\rightarrow K$ be the  map defined by $\Fcal'(f)=0$
for $f\in H^0\bigl(G/B, \Lcal (-\mu)\bigr)$ if $\mu\not\in pX^+$,  and 
 $\Fcal' : H^0\bigl(G/B, \Lcal (-p\mu)\bigr) \rightarrow
H^0\bigl(G/B, \Lcal (-\mu)\bigr)$ is the splitting map $\Fcal'_\mu$
defined above (as maps of abelian groups; without regard to the\break $k$-linear 
or $\bar{U}(\bfk^-)$-module structures).

As a first step we show:

\proclaim{Proposition}\label{bmodulemap} \hskip-8pt
$\Fcal'$ is a Frobenius\/{\rm -}\/linear graded endomorphism of~$K${\rm .} Moreover{\rm ,} it is a 
splitting{\rm .}
Further{\rm ,}   for any
$q\in\Zgil_+${\rm ,} $\lam \in X^+$ and 
$f\in H^0\bigl(G/B, \Lcal (-\lam)\bigr)${\rm :}
\begin{eqnarray}
\Fcal'(\bar{E}_i^{(pq)}f)&=&\bar{E}_i^{(q)}\Fcal'(f)\,, \ \ \hbox{ and
}\speqnu{1}
\\\Fcal'(\bar{F}_i^{(pq)}f)&=&\bar{F}_i^{(q)}\Fcal'(f), \speqnu{2}
\end{eqnarray}
where the action of $\bar{E}_i^{(m)}$ and  $\bar{F}_i^{(m)}$ comes from the 
canonical action of $G$ on $ H^0\bigl(G/B, \Lcal (-\lam)\bigr)${\rm .}

In particular{\rm ,} the splitting $\Fcal'$ is $B$\/{\rm -}\/canonical{\rm .}
\endproclaim

{\it Proof}.
The map $\Fcal'$ is clearly additive by definition. Further, 
by Proposition~3.3  and [Lu2, \S 8.15], 
 $\Fcal' $  satisfies (1) and (2). Next we prove the condition (b) of \S 6.1: 
Take $f \in H^0(G/B, \Lcal (-\lam ))$ and $ g\in 
H^0(G/B, \Lcal (-\mu ))$. Thus
 $f^p g\in H^0(G/B, \Lcal (-(p\lam+\mu )))$, so if $\mu$ is not divisible
by $p$, then  we have $\Fcal'(f^p g)=0=f \Fcal'(g)$, which proves the claim
in this case.
Now, we consider the case when $\mu$ is divisible by $p$: It is easy to  show that the following diagram is
commutative: \pagebreak

\figin{37}{1000}%
\noindent where $m (f\otimes g ) := m_f (g)$ (cf.\ Lemma 4.9) and $\bar{m} $ is defined 
similarly and {\rm Id} is the identity map.
 The commutativity
of the
diagram after base change implies  (b)  of \S 6.1 since $\Fcal_\lam f = f^p$. 
Next 
$\Fcal' (1) = 1$ (as is easy to see), showing that $\Fcal'$ is a splitting. 

Finally we prove that ${\cal F}'$ is $B$-canonical.  For any
simple root $\alpha_i$, the corresponding one-parameter subgroup
$x_{\alpha_i}(z)$ in $B$ can be written as $x_{\alpha_i} (z)=
\sum_{m\geq 0} z^m \Ebar_i^{(m)}$.  Then, for any $a\in K$,  
  \begin{eqnarray}
 &&\hskip-.25in (x_{\alpha_i}(z)*{\cal F}' )(a)
\speqnu{3}\\&&=\  x_{\alpha_i}(z) \cdot
{\cal F}'
(x_{\alpha_i}(-z) \cdot a)\nn \\[3pt]  &&=\  x_{\alpha_i}(z)\cdot{\cal F}' \bigl( \sum_{m\geq 0} (-z)^m
\Ebar_i^{(m)}\cdot a\bigr) \nn\\[3pt]  
 &&= \ x_{\alpha_i}(z)\cdot{\cal F}' \left(\left(\left( \sum_{n\geq
0}(-z)^{np}
\Ebar_i^{(np)} \right)  \left( \sum^{p-1}_{m= 0} (-z)^m
\Ebar_i^{(m)} \right)\right) \cdot a \right),\nn \\[3pt]  
 &&\hskip1in\mbox{since $\Ebar_i^{(np)} \Ebar_i^{(m)}=\Ebar_i^{(np+m)}$
over $k$, for $0\leq m<p$} \nn\\[3pt]  &&=\  \biggl(x_{\alpha_i}(z) \bigl( \sum_{n\geq 0} (-z)^n\,
\Ebar_i^{(n)}\bigr)\biggr)\cdot
\, {\cal F}'\, \Bigl(\sum_{m=0}^{p-1} (-z)^m\, \Ebar_i^{(m)}\cdot
a\Bigr),
\quad\mbox{by (1)}\nn\\[3pt]  
 &&=\ \sum^{p-1}_{m= 0} {\cal F}' ((-z)^m\, \Ebar_i^{(m)}\cdot a). \nn
  \end{eqnarray} 
Define $\sigma_m(a) = \Fcal'((-1)^m  \Ebar_i^{(m)}\cdot a).$ Then 
$$
 \sigma_m(z^pa) = z \sigma_m(a).
\tag{4}
$$
Now since $k$ is an infinite field and $x_{\alpha_i}(z)*{\cal F}' \in $
End$_F(K)$, it is easy to see from (3)--(4) that  $\sigma_m \in $ End$_F(K)$. 
This proves the defining property (d) of canonical splitting as in
\S 6.1.  The proof that, for any $t\in T$,  $t*\sigma = \sigma$ is easy
from the fact that $\Frps \bigl( H^0(\Xfk, \chi^\xi_{-p\mu})_{p\lam} \bigr)
\subset H^0(\Xbar, \bar{\chi}_{-\mu})_\lam$.
\hfill\qed\medbreak
 
Fix $\lam\in X^{++}$ (where $ X^{++}$ is the set of dominant regular weights). 
Replacing $\lam$ by $d\lam$ for $d\in\Zgil_+$ big enough, we may assume that
the embedding of $X = G/B \hookrightarrow \Pgil(\bar{V}_k(\lam))$, taking 
$gB$ to the line [$g\bar{v}_+$], is 
projectively normal (cf.\ [H, Chap. II, Ex.~5.14]),  where 
$ \bar{v}_+$ is a highest weight vector of $\bar{V}_k(\lam)$, $\Pgil(\bar{V}_k(\lam))$
denotes the space of lines in $\bar{V}_k(\lam)$ and  $\bar{V}_k(\lam) :=
 H^0(G/B,\Lcal(-\lam))^*$ is the Weyl module.
(Actually, the embedding is always
projectively normal. This has been shown in [RR] as a consequence of the
Frobenius splitting of $X$, or in [Li] as a consequence of standard
monomial theory; but we do not need this result, in fact we will derive 
this.) The homogeneous coordinate ring of $X\subset \Pgil(\bar{V}_k(\lam))$ is
hence
$k[X]=\bigoplus_{n\ge 0} H^0(X,\Lcal (-n\lam))$.

For $\tau\in W$, let $p_\tau\in H^0(G/B,\Lcal(-\lam))$ be a 
nonzero section of weight
$-\tau(\lam)$ (which is unique up to scalar multiples). Let $X_\tau$ be the  affine open subset defined by
$X_\tau:=\{x\in X\mid p_\tau(x)\not=0\}$. Note that $X_\tau$ depends only
on $\tau$ and not on the choice of $\lam$. (In fact, $X_\tau =
\tau U^-B/B \subset G/B$, where $U^-$ is the unipotent radical of the 
opposite Borel subgroup.) 

The affine coordinate ring $k[X_c]$ of the affine open subset $X_\tau$ is the degree
$0$
part  of the localization $k[X]_{(p_\tau)}$ of $k[X]$ at $p_\tau$:
$$
k[X_\tau]=\bigcup_{m\in\Zgil_+}\{\frac{f}{p_\tau^m}\mid f\in
H^0(X,\Lcal ({-m\lam}))\}.
$$
Here $f/p_\tau^m$ is equivalent to $f'/p_\tau^{m'}$ if there exists a
$q$ and nonzero $s \in\break H^0(X,\Lcal ({-q\lam}))$ such
that $s(fp_\tau^{m'}-f'p_\tau^m)=0$. Since $X$ is irreducible,\break this is 
equivalent to $f'p_\tau^{m-m'}=f$, if 
$m\ge m'$. We define now a splitting\break
$\Fcal'_{\tau}:k[X_\tau]\rightarrow k[X_\tau]$ as follows:
$$
\Fcal'_{\tau} ({\frac{f}{p_\tau^m}}):= \frac{\Fcal'(p_\tau^r f)}{p_\tau^{(r+m)/p}},
\ \hbox{\rm\ where\ }r\in\Zgil_+\hbox{\rm\ is such that\ }p\vert (r+m).
$$
To see that the map ${\cal F}'_\tau$ is well-defined, note first that if we fix an $r\in\Zgil_+$
such
that $p\vert (r+m)$, then
\begin{eqnarray*}
p_\tau^r f\in H^0\bigl(X,\Lcal (-(r+m)\lam )\bigr)&\Rightarrow&
\Fcal'(p_\tau^r f)\in  H^0\biggl(X,\Lcal (-((r+m)/ p)\lam)\biggr)\\&\Rightarrow&
 \frac{\Fcal'(p_\tau^r f)}{p_\tau^{(r+m)/p}}\in
k[X_\tau].
\end{eqnarray*}
Next, since $\Fcal'$ is a Frobenius-linear  endomorphism of $K$, it is easy to see 
that
the definition
is independent of the choice of $r$ and of the chosen representative
$\frac{f}{p_\tau^m}$, and moreover  $\Fcal'_\tau$ is a Frobenius-linear  
endomorphism.
Further, observe that $\Fcal'_{\tau}$ is indeed a splitting:
$$
\Fcal'_{\tau}((\frac{f}{p_\tau^m})^p)=
 \frac{\Fcal'(f^p)}{p_\tau^m}=\frac{f}{p_\tau^m}.
$$
It is easy to see that $\Fcal'_{\tau}$ does not depend upon the choice of 
$p_\tau$. However,
{a~priori}, the definition of the map $\Fcal'_{\tau}$ 
seems to depend on the choice of
$\lam$, but:

\proclaim{Lemma}
The definition of the splitting  $\Fcal'_{\tau}:k[X_\tau]\rightarrow k[X_\tau]$ is
independent of the choice of $\lam \in X^{++}$ {\rm (}\/such that $X\subset 
\Pgil(\bar{V}_k(\lam))$ is projectively normal\/{\rm ).}
\endproclaim

\demo{Proof}
Let $\mu\in X^{++}$ be such that the embedding
$X\subset \Pgil(\bar{V}_k(\mu))$ is projectively normal. 
Let $q_\tau\in H^0(X,\Lcal(-\mu ))$ be a nonzero weight vector of weight
$-\tau(\mu)$. Then, as mentioned  above, 
$X_\tau$ coincides with $\{x\in X\mid q_\tau(x)\not=0\}$, and its affine coordinate ring 
can be identified as the degree
$0$ part of the localization of $\bigoplus_{n\ge 0} H^0(X,\Lcal(-n\mu ))$
at $q_\tau$.

Take a function $h\in k[X_\tau]$ and write 
$h={f\over p_\tau^n}={g}/{q_\tau^m}$ with $f\in\break H^0(X,\Lcal {-(n\lam}))$ and 
 $g\in H^0(X,\Lcal{(-m\mu}))$. By replacing
${f}/{p_\tau^n}$ with ${f p_\tau^r}/{p_\tau^{n+r}}$ and ${g}/{q_\tau^m}$
with ${gq_\tau^t}/{q_\tau^{m+t}}$ if necessary, we may assume that $n=m$
and $n$ is divisible by $p$. Hence
\begin{eqnarray*}
{f\over p_\tau^n}&\hskip-5pt =\hskip-5pt&{g\over q_\tau^n}\Leftrightarrow
f q_\tau^n=g p_\tau^n
\quad \hbox{ on the whole  of } 
 X\,\\
 \Rightarrow \Fcal'(f q_\tau^n)&\hskip-5pt=\hskip-5pt&\Fcal'
(g p_\tau^n)
\Rightarrow q_\tau^{n/p}\Fcal'(f)= p_\tau^{n/p}\Fcal'(g).
\end{eqnarray*}
The last equality  finishes the proof of the lemma.
\enddemo

To glue  $\{\Fcal'_{\tau}\}_{\tau \in W}$ together to a global splitting, it 
remains to
show that the definitions are compatible  on the intersections;  i.e., we need 
to check   condition (1) of Definition~6.1: 

Consider the affine open set $X_{\{\tau,\kappa \}} :=X_\tau\cap X_{\kappa}=
\{x\in X\mid (p_\tau p_\kappa)(x)\not=0\}$. We define, as above, the map
$\Fcal'_{\{\tau, \kappa \}}:k[X_{\{\tau,\kappa\}}] \rightarrow k[X_{\{\tau,\kappa\}}]$; i.e., expressing
$$
k[X_{\{\tau,\kappa\}}]=\bigcup_{m\in\Zgil_+}\{\frac{f}{(p_\tau p_\kappa)^m}\mid f\in
H^0(X,\Lcal (-2m\lam ))\},
$$
we define
$$
\Fcal'_{\{\tau , \kappa\}} ({\frac{f}{(p_\tau p_\kappa)^m}}):= \frac{\Fcal'
((p_\tau p_\kappa)^r f)}
{(p_\tau p_\kappa)^{(r+m)/p}},\ \hbox{\rm\ where\ }r\in\Zgil_+\hbox{\rm\ is such
that\ }
p\vert (r+m).
$$
As above, one sees
that $\Fcal'_{\{\tau , \kappa \}}$ is a well defined Frobenius-linear 
endomorphism of $k[X_{\{\tau,\kappa\}}]$ and  is a splitting.

Take a regular function $h$ on $X_{\{\tau,\kappa\}}$ which is the restriction 
of a
regular function on $X_{\tau}$. So we can represent $h$ as
${f}/{p_\tau^n}$ as well as ${g}/{(p_\tau p_\kappa)^m}$ (for $f \in
H^0(X,\Lcal (-n\lam ))$ and $g \in
H^0(X,\Lcal (-2m\lam ))$), and these
functions coincide on $X_{\{\tau ,\kappa\}}$. Since $X$ is irreducible, 
this implies $f (p_\tau p_\kappa)^m=g p_\tau^n$ on the whole of $X$. 
By multiplying
${f}/{p_\tau^n}$ with $p_\tau^a/p_\tau^a$ and ${g}/{(p_\tau p_\kappa)^m}$
with $(p_\tau p_\kappa)^b/(p_\tau p_\kappa)^b$ if necessary, we  may assume
that
$n,m$ are divisible by $p$. Now,
$$
(p_\tau p_\kappa)^{m/p}\Fcal'(f)=\Fcal'((p_\tau p_\kappa)^m
f)=\Fcal'(p_\tau^n g)
=p_\tau^{n/p}\Fcal'(g),
$$
and hence, as functions on  $X_{\{\tau ,\kappa\}}$, 
$$
\Fcal'_{\tau} (\frac{f}{p_\tau^n})= \frac{ \Fcal'(f)}{p_\tau^{n/p}}=
 \frac{ \Fcal'(g)}{(p_\tau p_\kappa)^{m/p}}=\Fcal'_{\{\tau , \kappa \}}
 (\frac{g}{(p_\tau
p_\kappa)^{m}}).
$$
This proves the compatibility condition (1) of Definition~6.1.

Since the subsets  $X_\tau$, $\tau\in W$, form an affine open cover of $X$,
the above compatibility implies that one can glue the splittings $\Fcal'_{\tau}$ together to 
get a Frobenius splitting $\Theta$ 
on the whole of $X$. Moreover, since $\Fcal'$ is $B$-canonical, so is $\Theta$.

Summarizing, we have:

\proclaim{Theorem}
 The Frobenius\/{\rm -}\/linear graded endomorphism   $\Fcal'$ of $K$\break
{\rm (}\/cf.\ Proposition {\rm 6.2)}
induces a Frobenius splitting $\Theta$ of the flag variety $G/B$ by the method
described above{\rm .} Moreover{\rm ,} this splitting is $B$\/{\rm -}\/canonical{\rm .}\hfill\qed\endproclaim

By an  argument similar to the proof of the above theorem, using Proposition 
4.6  instead of 
Proposition~3.3, we obtain the following stronger result originally due 
to Ramanan-Ramanathan [RR]. (The commutativity of 
the diagram, analogous to that in the proof of Proposition 6.2, follows from 
ideas  similar to the proof of Proposition 4.11.)

\proclaim{Theorem}
 The flag variety $G/B$ is Frobenius $\Dcal$\/{\rm -}\/split{\rm ,} where $\Dcal$
is the line bundle $\Lcal (-2(p-1)\rho)$ together with the section 
$\phi := \theta (\sigma_o)${\rm ,} where ${\sigma}_o$ is as 
defined in \S {\rm 4.11,}  and $\theta$ is the 
isomorphism as in \S {\rm 6.1.} \hfill\qed
 \endproclaim

For an element $w\in W$ denote by $e_w:=w.{\rm id}\in G/B$ the corresponding
$T$-fixed
point. The closure $\overline{B.e_w}$ of the $B$-orbit ${B.e_w}$ is called a
{\it Schubert
variety} in $G/B$ and denoted $X(w)$. The closure $\overline{B^-.e_w}$
of the orbit with respect to the opposite Borel subgroup $B^-$ is called
an {\it opposite Schubert variety} in $G/B$ and denoted $X(w)^-$. Using
representation
theoretic arguments, we  show below that the above  splitting $\Theta $
is compatible with all the
Schubert varieties and opposite Schubert varieties.

As in Section~1, for $\lam \in X^+$, let $V_\xi(\lam)$ be
the Weyl module of highest weight $\lam$ for $\Ufk (\gfk)$ (over $\Zgil_\xi$).
Similarly, let $\bar{V}(\lam)$ be the Weyl  module for $\Ubar (\gfk)$ 
(again over $\Zgil_\xi$). For a base change $\Zgil_\xi \to \Bcal $, we denote 
$\Bcal \otimes V_\xi(\lam)$ by $V_\Bcal (\lam)$ and similarly  
$\bar{V}_\Bcal (\lam)$. 
In particular, we have $V_k(\lam)$ and $\bar{V}_k(\lam)$ for the base change 
$\Zgil_\xi \to k$ ($ \xi \mapsto 1$). Then  $V_k(\lam)$ and $\bar{V}_k(\lam)$
are canonically isomorphic modules for $\bar{U}_k (\gfk)$ via 
[Lu2, \S 8.15].

\vglue5pt

For any $w \in W$, the {\it Demazure module} 
$V_\xi(\lam)_w\subset V_\xi(\lam)$ is, by definition,  the
$\Ufk (\bfk)$--submodule $\Ufk (\bfk) v_{w\lam}$ generated by a primitive
extremal weight vector $v_{w\lam}\in V_\xi(\lam)$ of weight $w(\lam)$,
and the {\it opposite Demazure module} is the \, $\Ufk (\bfk^-)$-submodule 
$V_\xi(\lam)^-_w :=
\Ufk (\bfk^-) v_{w\lam}  \subset V_\xi(\lam)$. Similarly, define the Demazure 
module $\bar{V}(\lam)_w$ and  $\bar{V}(\lam)_w^- \subset \bar{V}(\lam).$

As in Section~3, for a $\Ufk(\gfk)$-module which is a $\Ufk^0$-weight module
$V=\bigoplus_{\mu\in X} V_\mu$,  denote by $V^{{1\over p}}$ the direct sum
of weight spaces $\bigoplus_{\mu\in pX} V_\mu$. It is easy to see that 
 $V^{{1\over p}}$ is equipped with  a $\Ubar(\bfk)$-module structure via 
$\Fr'_\bfk$ (as well as a $\Ubar(\bfk^-)$-module
structure via  $\Fr'_{\bfk^-}$).

By Proposition~3.3, the map $\Frps: 
H^0(\Xfk ,\chi^\xi_{-p\lam})^{\Fr'_\bfk}\to H^0(\Xbar , \bar{\chi}_{-\lam}),
 $ under the identification $\beta$ ((1) of \S 1.4) and a similar 
identification $\bar{\beta }$, decomposes as the composite of the 
 restriction $V_\xi(p\lam)^*\rightarrow (V_\xi(p\lam)^{\frac{1}{p}})^*$
followed by a map $(V_\xi(p\lam)^{\frac{1}{p}})^*\rightarrow \bar{V}(\lam)^*$. We get the
dual maps $\bar{V}(\lam)\hookrightarrow V_\xi(p\lam)^{\frac{1}{p}}
\hookrightarrow
V_\xi(p\lam)$. (The injectivity of the first map follows, e.g., 
from
Corollary 3.9.)
These maps are $\bar{U}(\bfk)$ as well as  $\bar{U}(\bfk^-)$-module maps
(where the latter two modules are equipped with $\Ubar (\bfk), \Ubar (\bfk^-)$-module structures via $\Fr'_{\bfk}$ and
$\Fr'_{\bfk^-}$ respectively); in particular, the $\mu$-weight space (for any $\mu \in X$ ) 
 $\bar{V}(\lam)_\mu$ is mapped to $V_\xi(p\lam)_{p\mu}$. For any $w \in W$,
 since 
the extremal weight space $V_\xi(p\lam)_{pw\lam}$ is of rank $1$, 
we have the inclusions of Demazure modules (obtained by the restrictions of the above inclusions):
  $$
\bar{V}(\lam)_w\hookrightarrow (V_\xi(p\lam)_w)^{\frac{1}{p}}
\hookrightarrow V_\xi(p\lam)_w,\quad
\bar{V}(\lam)_w^-\hookrightarrow (V_\xi(p\lam)_w^-)^{\frac{1}{p}}
\hookrightarrow V_\xi(p\lam)_w^- \,,\tag{1}
  $$
and the associated dual maps
\begin{eqnarray}
&&\bigl(V_\xi(p\lam)/V_\xi(p\lam)_w \bigr)^*\rightarrow 
\bigl( \bar{V}(\lam)/\bar{V}(\lam)_w\bigr)^*, \hbox{ and}\speqnu{2}\\
&&\bigl( V_\xi(p\lam)/V_\xi(p\lam)_w^-\bigr)^*\rightarrow 
\bigl(\bar{V}(\lam)/\bar{V}(\lam)_w^-\bigr)^*\,,\nn \end{eqnarray}
which are restrictions of $\Frps$ under the identifications $\beta$ and 
$\bar{\beta}$. 

\proclaim{Lemma}  Let $\lam\in X^+${\rm .}  Then the composite map $r_w
\circ\bar\theta \circ \bar \beta${\rm :}
  $$
\bar{V}_k(\lam )^*\, \lower3pt\hbox{$\stackrel{\bar{\scriptstyle\beta}}{\buildrel{\rightarrow}\over {\scriptstyle\sim}}$}\,
H^0
\bigl( \Ubar_k(\gfk )/\Ubar_k (\bfk ), (\bar{\chi}_{-\lam })_k\bigr) \,
\lower3pt\hbox{$\stackrel{\bar{\theta}}{\buildrel{\rightarrow}\over {\scriptstyle\sim}}$}\, H^0 (G/B, \Lcal
(-\lam )) \stackrel{r_w}{\rightarrow} H^0 (X_w, \Lcal (-\lam ))
  $$
has kernel precisely equal to $(\bar{V}_k(\lam )/\bar{V}_k(\lam
)_w)^*${\rm ,} where $r_w$ is the restriction map and $\bar{\theta}$ is
the isomorphism of \S {\rm 6.1.}
 
Similarly{\rm ,} the composite map  $r_w\circ\theta\circ \beta${\rm :}
  $$
V_k(\lam )^* \lower3pt\hbox{$\stackrel{{\scriptstyle\beta}}{\buildrel{\rightarrow}\over {\scriptstyle\sim}}$}\, H^0
\bigl( \Ufk_k(\gfk )/\Ufk_k (\bfk ), (\chi^\xi_{-\lam})_k\bigr) 
\,
\lower3pt\hbox{$\stackrel{{\theta}}{\buildrel{\rightarrow}\over {\scriptstyle\sim}}$}\, H^0 (G/B, \Lcal
(-\lam ))^{[1]} \stackrel{r_w}{\rightarrow} H^0 (X_w, \Lcal
(-\lam ))^{[1]}  
  $$
has kernel precisely equal to $(V_k(\lam )/V_k(\lam)_w)^*${\rm .}

Similar statements are true with $X_w$ replaced by $X^-_w$ and
$\bar{V}_k(\lam )_w$, $V_k(\lam )_w$ replaced by
$\bar{V}_k(\lam )^-_w${\rm ,}
$V_k (\lam )_w^-$ respectively{\rm .}
  \endproclaim 

 \demo{Proof}
We prove the first assertion.  The remaining ones are proved
similarly.  It is easy to see that the composite $\bar{\gamma} =
\bar{\theta}\circ\bar{\beta}$ is given by $\bar{\gamma}(f)(gB) = (g,
f(g\bar{v}_{\lam}) \bar{v}^*_{\lam})$, for $f\in\bar{V}_k(\lam )^*$, where
$\bar{v}_{\lam}$ is a nonzero vector of $\bar{V}_k(\lam )_\lam$
and  $\bar{v}^*_{\lam} \in (\bar{\chi}_{-\lam})_k = (\bar{V}_k(\lam
)_{\lam})^*$ is given by $\bar{v}_{\lam}^*(\bar{v}_{\lam} )=1$.  So
$r_w\circ \bar{\gamma} (f)=0 \Leftrightarrow f(bw\bar{v}_{\lam})
= 0$
for all $b\in B \Leftrightarrow f(\bar{V}_k(\lam )_w)\equiv 0$,
since
the $B$-module span and $\Ubar_k(\bfk )$-module span of
$w\bar{v}_{\lam}$ are the same.  This proves the lemma.
  \enddemo

The following result is originally due to Mehta-Ramanathan [MR1].\pagebreak

  \proclaim{Theorem}
Let $Z\subset G/B$ be a subscheme obtained from\break $\{X (w) , X(w)^- \}_{w\in W}$
by repeatedly taking scheme theoretic unions{\rm ,} intersections and
irreducible components{\rm .} Then $Z \subset G/B$ is compatibly split under the 
splitting $\Theta$ of $G/B$ given in Theorem~{\rm 6.4.} In particular{\rm ,} $Z$ is a reduced 
scheme.\endproclaim 

\demo{Proof}
We first show that, for any $w\in W$, $X(w)\subset G/B$ is compatibly
split:  Fix $\tau\in W$ such that $X_{\tau}\cap X(w)\ne
\phi$ and let $I(w)_{\tau}\subset k[X_{\tau}]$ be the ideal of
$X_{\tau}\cap X(w)$ in $k[X_{\tau}]$.  Choose $\lam\in X^{++}$ such
that the embedding $X\hookrightarrow \Pgil (\bar{V}_k(\lam ))$ is
projectively normal and let $p_{\tau}\in H^0(X, \Lcal (-\lam ))$ be
a nonzero section of weight $-\tau (\lam )$.  Then, by Lemma~6.6,
for $f\in H^0(X, \Lcal (-m\lam ))$,\break $\frac{f}{p^m_{\tau}}\in
I(w)_{\tau} \Leftrightarrow  f \in\gamma \bigl( (V_k(m\lam )/V_k(m\lam
)_w)^*\bigr)$, where $\gamma :=\theta\circ \beta$ is as in Lemma~6.6.
Hence, if $p|m$, for   $\frac{f}{p^m_{\tau}}\in
I(w)_{\tau}$,
$$ \Fcal'_\tau ( \frac{f}{p^m_{\tau}}) = \frac{\Fcal'(f)}{p^{m/p}_{\tau}}\in
I(w)_{\tau}\,,$$
by (2) of \S6.5 (with the base change $\Zgil_\xi \to k$).
This proves that $I(w)_{\tau}$ is stable under $\Fcal'_\tau$. The same argument  shows that the ideal $I(w)^-_\tau$ of
$X(w)_\tau^-$ in $k[X_\tau]$ is  again stable under $\Fcal'_\tau$. Hence,  any repeated sum and intersection of 
these ideals is stable under   $\Fcal'_\tau$. This shows that the splitting 
$\Theta$ of $G/B$ compatibly splits any unions and intersections of 
 $X(w)$ and  $X(w)^-$. Since any  irreducible component of a compatibly 
split subscheme is again compatibly 
split (cf.\ [R2, Prop.~1.9]), the theorem follows. 
\enddemo

 \numbereddemo{{R}emark} The map $\Frps$ (resp.\ $\Frps_\gamma) $ can
hence be viewed as a characteristic zero lift of the Frobenius
splitting (resp.\ Frobenius $\Dcal$-splitting, for the line bundle $\Dcal := 
\Lcal (-2(p-1)\rho))$ of $G/B$  on the level of quantum groups. In fact, on the level of quantum groups, as we saw in Sections~3
and~4,   the maps
$\Frps$ and  $\Frps_\gamma $ are defined for any odd integer $\ell > 1$
 (not necessarily prime
numbers). Also see [Li], where a `standard monomial' basis of
 $H^0(X(w),\Lcal(\lam))$ has been constructed. There the maps
$\Fcal'_{\tau}$ are  used to define the ``$\ell^{\rm th}$ root'' of a product
of extremal weight vectors in $V_\xi(\lam)_w^*$.
  \enddemo

  \section{Splitting of the diagonal in $G/B\times G/B$}

As in Sections  1 to 4, we continue to assume that $\ell > 1$ is an odd integer, which,
in addition, is coprime to $3$ if $G_2$ is a component of $\gfk$. 

\numbereddemo{Definition} Let $\bar{M}, \bar{N} \in \bar{\Ccal} (\bfk)$ be two $\Ubar (\bfk)$-modules. 
 By Theorem~2.3, we have  a
$\Ufk(\bfk)$-module map (for any $j \geq 0$):
$$
\pi_1: [\bar{M}\otimes H^j(\Xbar,\bar{N})]^{\Fr}
\longrightarrow \bar{M}^{\Fr}\otimes
H^j(\Xfk,\bar{N}^{\Fr} ),
\ a\otimes f\mapsto a \otimes \Frs(f).
$$
Similarly, 
by Theorem~3.8,  we have a $\bar{U}(\bfk)$-module 
map:
$$
\pi_2: [\bar{M}^{\Fr}\otimes
H^j(\Xfk ,\bar{N}^{\Fr})]^{\Fr'_\bfk}
\longrightarrow \bar{M}\otimes
H^j(\bar{X}, \bar{N}),
\ a\otimes g\mapsto a\otimes \Frps(g).
$$
The composition of the above two maps is the identity map by
Corollary~3.9.
By Theorem~2.3  and inducing $\pi_1$ we get, for all $i\in\Zgil_+$, 
$\Ufk(\gfk)$-module maps:
\figin{43a}{900}
\vglue-8pt\noindent 
where $\hat{\pi}_1$ is induced from $\pi_1$ and $\Fr^*_\Delta$ is by definition
the  
composite map $\hat{\pi}_1 \circ \Fr^*$.

Similarly, by inducing the map $\pi_2$ and the splitting
given by Theorem~3.8,
we get $\bar{U}(\bfk)$-module  maps:
\vglue-24pt
\figin{43b}{900}
\vglue-8pt\noindent where we denote by  $\Frdps$ the composition of the two maps $\Frps$ and 
$\hat{\pi}_2$.
\enddemo

From the functoriality of $\Frps$, we get $\Frps \circ \hat{\pi}_1 = 
 \hat{\bar{\pi}}_1 \circ \Frps$, where the last $\Frps$ is the map 
$H^i\bigl(\Xfk ,[\bar{M}\otimes
H^j(\bar{X},\bar{N})]^{\Fr}\bigl) \to H^i\bigl(\bar{X},\bar{M}\otimes
H^j(\bar{X} ,\bar{N})\bigl)$ and $\bar{\pi}_1 $ is the 
$\Ubar(\bfk)$-module map $ \bar{M} \otimes
H^j(\bar{X},\bar{N}) \to 
[\bar{M}^{\Fr}\otimes
H^j(\Xfk ,\bar{N}^{\Fr})]^{\Fr'_\bfk}$ obtained by applying 
$[ \,\,]^{\Fr'_\bfk}$ to $\pi_1$. Hence
$$  \Frdps\circ\Frds = \hat{\pi}_2 \circ  \hat{\bar{\pi}}_1 \circ \Frps \circ
\Fr^* = {\rm Id}\,, $$
by application of  Corollary~3.9  twice. 
 Summarizing, we have: 

  \proclaim{Proposition}\label{deltasplitt}
For any $\bar{M}, \bar{N} \in \bar{\Ccal} (\bfk)$ and $i,j \in \Zgil_+${\rm ,} 
there exists a functorial 
$\bar{U}(\bfk)$\/{\rm -}\/module  map
  $$
\Frdps:
H^i\bigl(\Xfk ,\bar{M}^{\Fr}\otimes
H^j(\Xfk ,\bar{N}^{\Fr} )\bigr)^{\Fr'_{\bfk}}
\rightarrow
H^i\bigl(\bar{X} ,\bar{M}\otimes
H^j(\bar{X} ,\bar{N})\bigl)
  $$
which is a splitting of the functorial $\Ufk(\km)$\/{\rm -}\/module  map
  $$
\Frds:
H^i\bigl(\bar{X} , \bar{M}\otimes
H^j(\bar{X} , \bar{N})\bigl)^{\Fr}
\rightarrow
H^i\bigl(\Xfk ,\bar{M}^{\Fr}\otimes
H^j(\Xfk ,\bar{N}^{\Fr})
\bigl).
  $$
Moreover{\rm ,} for any $m\geq 0${\rm ,} $1 \leq t \leq n$ and
$f\in H^i\bigl(\Xfk ,\bar{M}^{\Fr}\otimes
H^j(\Xfk ,\bar{N}^{\Fr} )\bigl)${\rm ,}
  $$
\Fbar_t^{(m)}\cdot (\Frdps f) = \Frdps (F_t^{(m\ell )}\cdot f), \tag{1}
  $$
where the action of $F_t^{(m)}$ on  $H^i\bigl(\Xfk ,\bar{M}^{\Fr}\otimes
H^j(\Xfk ,\bar{N}^{\Fr} )\bigl)$ comes from its action on the cohomology 
$ H^i\bigl(\Xfk ,M\bigr)$ for any $M \in \Ccal (\bfk)$ and similarly for the action
of $\Fbar_i^{(m)}${\rm .} \hfill\qed
  \endproclaim

{From now on, to the end of this section, use the same notation and assumptions as in Sections 5 and~6.
 In particular, assume that  $\ell=p$ is an odd prime and
$p >3$ if $\km$ has a component of type ${\tt G}_2$, and let $k$
be an algebraically closed field of characteristic~$p$. 
Also  let $G$ and $B$ be  as in Section~6.}

We view $Y:=G/B\times G/B$ as a $G$-variety
via the diagonal action. Note that the canonical map $G\times_B G/B
\rightarrow Y$, defined by $(g,g'B)\mapsto (gB,gg'B),$ is  a
$G$-equivariant isomorphism. It follows immediately that the $G$-orbits
in $Y$ are precisely the images of $G\times_B C(w)$, where $C(w):=B e_w$ is the
$B$-orbit of the $T$-fixed point $e_w\in G/B$ for  $w\in W$, and the
closure $Y(w)$ of the image of $G \times_B C(w)$   in $Y$ is the image of $G\times_B X(w)$,
where $X(w)\subset G/B$ is  the  Schubert variety $C(w)$. The diagonal
$G/B\subset G/B\times G/B$ corresponds here to $G\times_B X({\rm id})$, where
$X({\rm id})=e_{{\rm id}}$ is the one point Schubert variety. The 
varieties $Y(w)$ are called the {\it $G$-Schubert varieties} in~$Y$.

Of course, $G\times G$ satisfies the conditions of the preceding section,
so the flag variety $(G\times G)/(B\times B)\simeq G/B\times G/B$ is 
Frobenius split. However,  in general,  the $G$-Schubert varieties $Y(w)$
 are not compatibly split with respect
to this splitting. So we need to consider a different splitting. 

Recall [APW, Prop.\ 2.16] that for $\Ufk (\gfk )$-module
$M\in\Ccal (\gfk )$ which is $\zi$-flat and $N\in \Ccal (\bfk )$,
there is a $\Ufk (\gfk )$-module isomorphism
  $$
\delta : H^0(\Xfk , N)\otimes M \simarr H^0(\Xfk ,N\otimes M)
  $$
given by
  $$
\delta (f\otimes m)(x) = \sum_if(x'_i)\otimes (x_i''m), 
  $$
for $f\in H^0 (\Xfk ,N)$, $m\in M$ and $x\in \Ufk (\gfk )$, where
$\Delta x = \sum_i x'_i\otimes x_i''$.
\vglue5pt
There is a similar $\Ubar(\gfk )$-module isomorphism for
$\Nbar\in\Cbar (\bfk )$ and $\Zgil_\xi$-flat $\bar{M}\in\Cbar (\gfk )$:
  $$
\bar{\delta} : H^0(\Xbar ,\bar{N} )\otimes\Mbar \simarr H^0(\Xbar
,\bar{N}\otimes\Mbar ).
$$
In particular, for $\Mbar ,\Nbar\in\bar{\Ccal}(\bfk )$ such that
$H^0(\Xfk , \Nbar^{\Fr})$ and $H^0(\bar{X}, \Nbar )$ are both
 $\zi$-flat, the $\Ubar(\bfk )$-module map $\Frps_{\Delta}$ of \S7.1 (for
$i=j=0$) under the identifications $\delta$ and $\bar{\delta}$ can be
rewritten as 
  $$
\Frps_{\Delta} : [ H^0(\Xfk ,\Mbar^{\Fr}) \otimes H^0(\Xfk ,\Nbar^{\Fr})]
^{\Fr'_\bfk} \to H^0 (\Xbar ,\Mbar )\otimes H^0(\Xbar ,\Nbar ).
  $$
Similarly, we can rewrite the  $\Ufk (\gfk )$-module map $\Fr^\ast_\Delta$ as  
 $$
\Fr^*_{\Delta} :  H^0(\Xbar ,\Mbar)^{\Fr} \otimes H^0(\Xbar ,\Nbar)^{\Fr}
 \to H^0 (\Xfk ,\Mbar ^{\Fr}) \otimes H^0(\Xfk ,\Nbar^{\Fr} ).
  $$

In particular, taking $\Mbar = \bar{\chi}_{-\lam}, \Nbar = \bar{\chi}_{-\mu}$
(for $\lam, \mu \in X^+$) and considering the base change $\Zgil_\xi \to k \,
(\xi \mapsto 1),$ the maps $\Frps_{\Delta}$ and $\Fr^*_{\Delta}$   under the 
identifications  $\theta$ and $\bar{\theta}$ of \S 6.1, correspond respectively to the maps 
\begin{eqnarray*}
 &&\Fcal'_{(\lam , \mu)}: [H^0(G/B,\Lcal({-p\lam}))^{[1]}\otimes H^0(G/B,\Lcal 
({-p\mu}))^{[1]}]^{\Fr'_\bfk} \\&&\phantom{\Fcal'_{(\lam , \mu)}:}
\longrightarrow
H^0(G/B,\Lcal({-\lam}))\otimes H^0(G/B,\Lcal({-\mu}))
\\
\noalign{\noindent and}  
&&\Fcal_{(\lam , \mu)}:H^0(G/B,\Lcal(-{\lam}))^{\Fr}\otimes 
H^0(G/B,\Lcal(-{\mu}))^{\Fr}\\ &&\phantom{\Fcal'_{(\lam , \mu)}:}
\longrightarrow
H^0(G/B,\Lcal({-p\lam}))^{[1]}\otimes H^0(G/B,\Lcal(-{p\mu}))^{[1]}\,.
\end{eqnarray*}
Moreover, $\Fcal'_{(\lam , \mu)} \circ \Fcal_{(\lam , \mu)} = {\rm Id}.$
Observe that the map $\Fcal_{(\lam , \mu)}$ is a $\Ufk (\gfk)$-module map and
$\Fcal'_{(\lam , \mu)}$ is a $\Ubar (\bfk)$-module map (under the diagonal 
actions). 

Consider the $X^{+}\times X^{+}$-graded algebra
$$
K_\Delta:= \bigoplus_{\lam,\mu\in X^{+}}
H^0(G/B,\Lcal ({-\lam}))\otimes H^0(G/B,\Lcal ({-\mu})),
$$
under the multiplication $ (s\otimes t)\cdot (f \otimes g) := (sf)\otimes tg.$ 
We abbreviate $$H^0(G/B,\Lcal ({-\lam}))\otimes H^0(G/B,\Lcal (-\mu))$$
to $K_{\lam,\mu}$. Denote by $\Fcal'_\Delta:K_\Delta \rightarrow K_\Delta$ 
the graded map defined by ${\Fcal'_\Delta} (f\otimes g) =0 $ for 
$f\otimes g \in K_{\lam,\mu}$ if $(\lam, \mu) \notin pX^+ \times pX^+$, and 
let 
${\Fcal'_\Delta}_{| K_{p\lam,p\mu}}$
be the splitting map $\Fcal'_{(\lam ,\mu)}$ (as  maps of abelian groups;
without regard to the $k$-linear or $\Ubar (\bfk)$-module structures).

Similar to Proposition~6.2,  we have the following:  

\proclaim{Proposition} \hskip-8pt
 $\Fcal'_\Delta$
is a Frobenius\/{\rm -}\/linear graded  endomorphism of~$K_\Delta${\rm .} Moreover{\rm ,}
it is a splitting{\rm .} Further{\rm ,}  for any
 $q\in\Zgil_+$ and $f\otimes g\in K_\Delta${\rm ,}  
$$
\Fcal'_\Delta (\bar{E}_i^{(pq)}\cdot (f\otimes g))=\bar{E}_i^{(q)}\cdot
\Fcal'_\Delta (f\otimes g)\,, \hbox{ and}\tag{1}
$$
$$
\Fcal'_\Delta (\bar{F}_i^{(pq)}\cdot (f\otimes g))=\bar{F}_i^{(q)}\cdot
\Fcal'_\Delta (f\otimes g)\,,\tag{2}
$$
where $\bar{E}_i^{(q)}$ and $\bar{F}_i^{(q)}$ act diagonally{\rm .}

In particular{\rm ,} $\Fcal'_\Delta $ is $B$\/{\rm -}\/canonical for the diagonal action of 
$B$ on $K_\Delta${\rm .} 
\endproclaim

\demo{Proof}
The map $\Fcal'_\Delta $ is clearly additive and the  properties (1) and (2)
 follow from Proposition~7.2.
We proceed now to prove   property (b) of Definition~6.1  (following
 the proof of Proposition~\ref{bmodulemap}):

Let $s\otimes t\in K_{\lam,\mu}$ and $f\otimes g\in K_{\eta,\nu}$.
Since $s^p f\otimes t^p g\in K_{p\lam+\eta,p\mu+\nu}$, 
$$
\Fcal'_\Delta (s^p f\otimes t^p g)=0=(s\otimes t) \cdot \Fcal'_\Delta (f\otimes g)
\hbox{{\rm\ \ if\ \ }} (\eta,\nu)\not\in pX^{+}\times pX^{+}.
$$
Assume now that $f\otimes g\in K_{p\eta,p\nu}$.
Consider the commutative diagram:
  \figin{cdill}{890}

\noindent where $H$ denotes $H^0$,  $H(\chi^{\xi}_{\lam})$ denotes $ H^0(\Xfk ,\chi^{\xi}_\lam)$ (similarly 
$H(\bar{\chi}_{\lam})$),  
$\hat{m}$ is the map induced from the $\Ufk (\bfk )$-module map 
$$
\chi^{\xi}_{-p\lam}\otimes H^0(\Xfk ,\chi^{\xi}_{-p\mu}) 
\otimes \chi^{\xi}_{-p\eta} \otimes H^0(\Xfk ,\chi^{\xi}_{-p\nu}) \to
\chi^{\xi}_{-p\lam -p\eta} \otimes H^0(\Xfk ,\chi^{\xi}_{-p\mu -p\nu})$$
 taking
$(a\otimes b)\otimes (c\otimes d) \mapsto (a\otimes c)\otimes m_b(d)$ 
(where $m_b$ is
the map defined in Lemma 4.9) and $m(\sigma\otimes \sigma')=m_{\sigma}(\sigma')$.  Let
$\zeta$
be the composite map $\hat{m}\circ m$; $\bar{\zeta}$ is analogously defined as
$\hat{\bar{m}}\circ\bar{m}$.  (Checking the  commutativity of the above diagram is   routine
 if we keep track of the definitions of various maps involved.)

The commutativity of the above diagram after base change implies 
  property (b)
of Definition~6.1  since
  $$
\Fcal_{(\lam ,\mu)} (s\otimes t) = s^p\otimes t^p, \tag{3}
  $$
for $s\in H^0(G/B, \Lcal (-\lam ))$ and $t\in H^0(G/B, \Lcal (-\mu ))$.  
Observe that
(3) follows from the corresponding property: $\Fcal_{\lam}(s) = s^p$ 
(cf.\ $\S$6.1)
together with the identity:
  $$
\Fcal_{(\lam ,\mu)}(s\otimes t) = \Fcal_{\lam}(s)\otimes \Fcal_{\mu}(t). \tag{4}
  $$

Next, it is easy to see that $\Fcal'_{\Delta}(1)=1$ and hence $\Fcal'_{\Delta}$ is a
splitting.  The assertion that $\Fcal'_{\Delta}$ is $B$-canonical  follows from (1)
by the same argument as that used in the proof of Proposition~6.2. 
\enddemo

Analogous to Lemma~6.6, we get the following:
\proclaim{Lemma} Let $\lam , \mu \in X^+${\rm .} Then the composite map\/{\rm :}
\begin{eqnarray*}
 [V_k(\mu) \otimes  V_k(\lam)]^* \stackrel{\nu}\simarr   V_k(\lam)^*
\otimes  V_k(\mu)^* &\!\!  \stackrel{\tilde{\gamma}}\simarr\!\! &  H^0(G/B \times G/B, \Lcal(-\lam) \boxtimes 
\Lcal(-\mu))^{[1]}\\
&\!\! \to\!\!& H^0(Y(w),  \Lcal(-\lam) \boxtimes 
\Lcal(-\mu))^{[1]}
\end{eqnarray*}
has kernel precisely equal to $ [(V_k(\mu) \otimes  V_k(\lam))/ (\Ufk_k (\gfk).
(v_{w\mu} \otimes v_\lam))
]^*${\rm ,} where $\Ufk_k (\gfk)$ acts diagonally{\rm ,} $v_{w\mu}$ is a nonzero vector 
of weight $w\mu${\rm ,}  the inverse of the first isomorphism $\nu$ is given by
$\nu^{-1} (f \otimes g) (y \otimes x) = f(x) g(y)$ for $f \in  V_k(\lam)^*,
g\in  V_k(\mu)^*, x\in  V_k(\lam),$  $y \in  V_k(\mu),$ and $\tilde{\gamma}$
is induced by the isomorphism $\gamma = \theta \circ \beta$ as in Lemma~{\rm 6.6.}

 A similar statement
is true with $V_k$ replaced by $\bar{V}_k${\rm .}  \hfill\qed
\endproclaim

Dualize the map 
$$\Frps_\Delta :    H^0(\Xfk ,\chi^\xi_{-p\lam}) \otimes 
H^0(\Xfk ,\chi^\xi_{-p\mu})
 \to H^0 (\Xbar ,\bar{\chi}_{-\lam} )\otimes H^0(\Xbar ,\bar{\chi}_{-\mu} )
  $$
to get the map 
$$ \kappa : \bar{V}(\mu) \otimes \bar{V}(\lam) \to V_\xi(p\mu) \otimes V_\xi(p\lam).$$
The map $\kappa$ commutes with $\Ubar (\bfk)$ and $\Ubar (\bfk^-)$-actions, 
where  $ V_\xi(p\mu) \otimes V_\xi(p\lam)$ is equipped with  the diagonal $\Ubar (\bfk)$ 
and $\Ubar (\bfk^-)$-actions via $\Fr'$. It can be seen that 
$$\Frps_\Delta(\beta(v_{p\lam}^*) \otimes g) = \bar{\beta}(\bar{v}^*_\lam) 
\otimes \Frps g\,, \,\hbox{for any} \, g \in H^0(\Xfk ,\chi^\xi_{-p\mu}),$$
where $v_{p\lam}^*\in V_\xi(p\lam)^*$ is defined by $v_{p\lam}^* (v_{p\lam})
=1$ and $v_{p\lam}^* (v) = 0$ for any weight vector of weight $\neq p\lam$
 ($\bar{v}_{\lam}^*$ is defined similarly). Dualizing, we obtain: 
$$ \kappa (\bar{v}_{w\mu} \otimes \bar{v}_{\lam}) \in \Zgil_\xi (v_{pw\mu} 
\otimes v_{p\lam})\,.$$
Hence $$ \kappa (\Ubar(\gfk)\cdot (\bar{v}_{w\mu} \otimes \bar{v}_{\lam}))
\subset \Ufk (\gfk) \cdot (v_{pw\mu} \otimes v_{p\lam}).$$

By the same proof as that of Theorems~6.4  and~6.7, we obtain the 
following from Proposition~7.3  and Lemma~7.4. It was first proved by 
Mehta-Ramanathan [MR2] that  $G/B \times G/B$  (more generally 
  $G/P \times G/P'$) admits a Frobenius splitting 
which compatibly splits all the $G$-Schubert subvarieties.  
 
  \proclaim{Theorem}\label{diagonalsplitttheorem}
Let $Y$ be the $G$\/{\rm -}\/variety $G/B \times G/B$ {\rm (}\/under the diagonal action of $G${\rm ).}
 Then the Frobenius\/{\rm -}\/linear graded endomorphism   
$\Fcal'_\Delta $ of $K_\Delta$ 
{\rm (}\/cf.\ Prop.~{\rm 7.3)}
induces a Frobenius splitting $\Theta_Y$ of $Y$ by a similar   method given 
in Section~{\rm 6.} Moreover{\rm ,} this splitting is $B$\/{\rm -}\/canonical{\rm .}

Further{\rm ,} any subscheme $Z\subset Y$ obtained from the $G$\/{\rm -}\/Schubert varieties 
$\{Y (w) \}_{w\in W}$
by repeatedly taking scheme theoretic unions{\rm ,} intersections and
irreducible components is compatibly split{\rm .}  In particular{\rm ,} $Z$ is a reduced 
scheme{\rm .} \phantom{MW47}
\hfill\qed\endproclaim

\vglue-12pt

  \section{Frobenius splitting of quantized Bott-Samelson desingularization}
\vglue-9pt

In this section, as in Sections~1 to~4, $\ell >1$ is  an odd integer which is assumed 
to be
coprime to 3 if $G_2$ is a component of $\gfk = \gfk (A)$.
For any $1\leq i\leq n$, let $\bfk \subset \pfk_i$ be the minimal parabolic
subalgebra of $\gfk$.

 \advance\theoremcount by 1

\medbreak {\it Definition}  8.1  [APW,  \S 5.1]. For any sequence of simple
reflections $\wfk = (s_{i_1}, \cdots ,s_{i_m})$, define the functor
$D_{\wfk} : \Ccal (\bfk ) \to \Ccal (\bfk )$ inductively by
  $$
D_{\wfk}(M) = H^0 (\Ufk (\pfk_{i_1})/\Ufk (\bfk ), D_{\wfk'}(M)),
  $$
where $\wfk'$ is the subsequence $(s_{i_2}, \cdots ,s_{i_m})$.

Similarly, define $\Dbar_{\wfk} : \Cbar (\bfk ) \to \Cbar (\bfk )$ by
  $$
\Dbar_{\wfk}(M) = H^0(\Ubar (\pfk_{i_1})/\Ubar (\bfk ),
\Dbar_{\wfk'}(M)).
  $$
\

\phantom{goodlook}
\vglue-36pt
Both the functors $D_{\wfk}$ and $\Dbar_{\wfk}$ are left exact.
Denote their right derived functors respectively by
$H^i(\Zcal_{\wfk},-)$ and $H^i(\Zbar_{\wfk},-)$.  Let $k$ be a field which is
a $\zi$-algebra  $\zi \to k$. Let $ \Ccal_k (\bfk )$ (resp.\ $\Cbar_k (\bfk )$) be the analogue of the category $\Ccal (\bfk )$ (resp.\ $\Cbar (\bfk )$), where the base ring $\zi$ is replaced by $k$. For any 
$M \in  \Ccal_k (\bfk )$ we can similarly define $H^*(\Zcal_{\wfk}^k, M)$
as the derived functors of $D_{\wfk}^k(M) := 
H^0(\Ufk_k (\pfk_{i_1})/\Ufk_k (\bfk ), D_{\wfk'}^k(M)).$ (Also,
$H^*(\Zbar_{\wfk}^k, \Mbar)$, for $\Mbar \in \Cbar_k (\bfk ),$ is defined analogously.)

 By [APW, Cor.~2.13] (actually their setting differs a bit, but the 
same proof works),   for any $M\in\Ccal^0_k$, the $\Ufk_k (\bfk
)$-module $N := H^0(\Ufk_k (\bfk )/\Ufk^0_k,M)$ is an injective object of 
 $ \Ccal_k (\bfk )$. In particular, it is 
acyclic for the
functor $D_{\wfk}^k$, i.e., 
  \begin{equation}
H^i(\Zcal_{\wfk}^k, N) = 0, \hbox{  for all $i>0$}.
  \end{equation}

Similarly, for $\bar{M}\in \Cbar^0_k$,
  \begin{equation}
H^i(\Zbar_{\wfk}^k, \bar{N}) = 0, \hbox{  for all $i>0$},
  \end{equation}
where
  $$
\bar{N} := H^0(\Ubar_k (\bfk )/\Ubar^0_k, \bar{M}).  
  $$

\numbereddemo{{R}emark} It is very likely that (1) and (2) above remain true for any 
$\zi$-algebra $k$ (not
only when $k$ is a field).
\enddemo
  
  \proclaim{Theorem} For any sequence $\wfk  = (s_{i_1}, \cdots ,s_{i_m})$ and any $\Mbar \in \Cbar_k (\bfk
)${\rm ,} there exists a  functorial $\Ufk_k (\bfk )$\/{\rm -}\/module map
  $$
\Fr^*_{\wfk} : H^i(\Zbar_{\wfk}^k, \Mbar)^{\Fr} \to H^i(\Zcal_{\wfk}^k,
{\Mbar}^{\Fr}).
  $$\endproclaim

  \advance\eqcount by -2
{\it Proof}.
Consider the standard resolution of $\Mbar$ in the category $\Cbar_k (\bfk
)$ (cf.\ $(\ast )$ of (1.5)):
  \begin{equation}
0 \to \Mbar \to \Qbar_0 \to \Qbar_1 \to \cdots ,
  \end{equation}
and also the standard resolution of $\Mbar^{\Fr}$ in the category $\Ccal_k
(\bfk )$:
  \begin{equation}
0 \to \Mbar^{\Fr} \to Q_0 \to Q_1 \to \cdots .
  \end{equation}
As in the proof of Theorem~2.3, there are $\Ufk_k (\bfk )$-module
homomorphisms\break $\theta_i : \Qbar_i^{\Fr} \to Q_i$, for all $i \geq 0$, making
the diagram (D$'$) of the proof of Theorem~2.3 commutative.

For any $\Qbar\in\Cbar_k (\bfk )$, we first construct a functorial  
$\Ufk_k(\bfk)$-module map
  $$
\Fr^*_{\wfk} : H^0(\Zbar_{\wfk}^k, \Qbar )^{\Fr} \to H^0
(\Zcal_{\wfk}^k, \Qbar^{\Fr}),
  $$
by induction on $\ell (\wfk ) = m$.

By definition,
  $$
H^0 (\bar{Z}_{\wfk}^k, \Qbar ) = H^0 \bigl(\Ubar_k (\pfk_{i_1})/\Ubar_k
(\bfk ), H^0(\Zbar_{\wfk'}^k, \Qbar )\bigr),
  $$
where $\wfk' := (s_{i_2},\cdots ,s_{i_m})$.  By Lemma 2.2(b) (with
$\gfk$ replaced by $\pfk_{i_1}$), we have
a $\Ufk_k (\bfk )$ (in fact a $\Ufk_k (\pfk_{i_1})$)-module
homomorphism
  \begin{eqnarray*}
 && \varphi': H^0\bigl( \Ubar_k (\pfk_{i_1})/\Ubar_k (\bfk ),
H^0(\Zbar_{\wfk'}^k, \Qbar )\bigr)^{\Fr}\\
&&\qquad  \rightarrow H^0\bigl(\Ufk_k (\pfk_{i_1})/\Ufk_k (\bfk ), H^0(\Zbar_{\wfk'}^k, \Qbar
)^{\Fr}\bigr) .
  \end{eqnarray*}
Also, by induction, we have a $\Ufk_k (\bfk )$-module homomorphism
  $$\Fr^*_{\wfk'}:
H^0(\Zbar_{\wfk'}^k, \Qbar )^{\Fr} \rightarrow H^0(\Zcal_{\wfk'}^k,
\Qbar^{\Fr}),
  $$
which induces a $\Ufk_k (\pfk_{i_1})$ (in particular $\Ufk_k (\bfk
)$)-module homomorphism
  $$
\varphi'': H^0\bigl( \Ufk_k (\pfk_{i_1})/\Ufk_k (\bfk ), H^0(\Zbar_{\wfk'}^k,
\Qbar )^{\Fr}\bigr) \rightarrow H^0(\Zcal_{\wfk}^k, \Qbar^{\Fr}).
  $$
Now we set $\Fr^*_{\wfk}$ as  the composition $\varphi:= 
\varphi''\circ \varphi'$ (which is a
$\Ufk_k (\pfk_{i_1})$-module
homomorphism)
  $$
\varphi: H^0(\Zbar_{\wfk}^k, \Qbar )^{\Fr} \rightarrow H^0(\Zcal_{\wfk}^k,
\Qbar^{\Fr}).
  $$
This completes the induction.

Replacing $\Qbar$ by $\Qbar_i$, we get $\Ufk_k (\pfk_{i_1})$-module
homomorphisms
$$
 H^0(\Zbar_{\wfk}^k, \Qbar_i)^{\Fr} \stackrel{\varphi_i }\rightarrow 
H^0(\Zcal_{\wfk}^k,
\Qbar_i^{\Fr})
\stackrel{\theta^*_i}\rightarrow H^0(\Zcal_{\wfk}^k, Q_i),
$$
where $\theta^\ast_i$ is induced from $\theta_i$.
The resolutions (1) and (2) give rise to the cochain complexes
  \begin{eqnarray}
H^0(\Zbar_{\wfk}^k, \Qbar_0)^{\Fr} \to H^0(\Zbar_{\wfk}^k, \Qbar_1)^{\Fr}
\to \cdots , \quad\hbox{ and} \\
H^0(\Zcal_{\wfk}^k, Q_0) \rightarrow H^0(\Zcal_{\wfk}^k, Q_1) \to \cdots
.  
  \end{eqnarray}

The maps $\theta_i^* \circ \varphi_i$ give a cochain map from the cochain complex (3) to the
cochain complex (4).  Taking cohomology,
we get the desired map
\medbreak
\hfill ${\displaystyle
\Fr^*_{\wfk} : H^i(\Zbar_{\wfk}^k, \Mbar)^{\Fr} \to H^i(\Zcal_{\wfk}^k,
\Mbar^{\Fr}).
 }$ 
  \hfill\qed

  \proclaim{Theorem} For any sequence $\wfk$ and any $M\in\Ccal_k (\bfk
)${\rm ,} there exists a  functorial $\Ubar_k (\bfk )$\/{\rm -}\/module map
  $$
\Fr^{\prime\ast}_{\wfk} : H^i (\Zcal_{\wfk}^k, M)^{\Fr'_{\bfk}}
\rightarrow H^i\bigl(\Zbar_{\wfk}^k, M^{\Fr'_{\bfk}}\bigr) .
  $$\endproclaim

 {\it Proof}.
We first define the map $\Frps_{\wfk}$ at the $H^0$-level.  Applying
Proposition~3.3  for $\gfk$ replaced by $\pfk_{i_1}$, we get the
$\Ubar_k(\bfk )$-module map
  \begin{eqnarray*}
\beta_1: H^0(\Zcal^k_{\wfk}, M)^{\Frp_{\bfk}} &= &
H^0\bigl( \Ufk_k(\pfk_{i_1})/\Ufk_k(\bfk ), H^0(\Zcal^k_{\wfk'},
M)\bigr)^{\Frp_{\bfk}} \\
&&\to  H^0\bigl( \Ubar_k(\pfk_{i_1})/\Ubar_k(\bfk ),
H^0(\Zcal^k_{\wfk'}, M)^{\Frp_{\bfk}}\bigr) .
  \end{eqnarray*}

By induction on $\ell (\wfk )$, we have the $\Ubar_k(\bfk )$-module
map
  $$
\Frps_{\wfk'} : H^0(\Zcal^k_{\wfk'}, M)^{\Frp_{\bfk}} \to
H^0(\Zbar^k_{\wfk'}, M^{\Frp_{\bfk}}).
  $$
Inducing $\Frps_{\wfk'}$, we get the $\Ubar_k(\pfk_{i_1})$-module map
  $$
\beta_2: H^0\bigl( \Ubar_k(\pfk_{i_1})/\Ubar_k(\bfk ),
H^0(\Zcal^k_{\wfk'},
M)^{\Frp_{\bfk}}\bigr) 
\to  H^0\bigl( \Ubar_k(\pfk_{i_1})/\Ubar_k(\bfk ),
H^0(\Zbar^k_{\wfk'}, M^{\Frp_{\bfk}})\bigr) .
  $$
Composing $\beta_2\circ\beta_1$, we get the desired map
$\Frps_{\wfk}$ at the $H^0$-level.  Now, by a proof parallel to that
of the proof of Theorem~3.8, we define $\Frps_{\wfk}$ for an
arbitrary $H^i$. \phantom{sometimes}
\hfill\qed

  \proclaim{{C}orollary}For any $\Mbar\in\Cbar_k (\bfk )$ and sequence $\wfk${\rm ,} 
  $$
\Fr^{\prime\ast}_{\wfk}\,\circ\,  \Fr^*_{\wfk} : H^i(\Zbar_{\wfk}^k, \Mbar)
\rightarrow H^i(\Zbar_{\wfk}^k, \Mbar)
  $$
is the identity map{\rm .}
  \endproclaim 

 {\it Proof}. From the definition of the maps involved, it is easy to see 
by induction on 
$\ell (\wfk)$ that  
 the corresponding property is true at the $H^0$-level. Now the validity of the corollary for general $i$ follows by 
the same argument as that of the proof of Corollary~3.9. 
  \hfill\qed

  \numbereddemo{Definition}  Let $\wfk = (s_{i_1},\cdots ,s_{i_m})$ be any
sequence of simple reflections.  We need a certain generalization of
the functor $D_{\wfk}$, still denoted by $D_{\wfk} : \Ccal_k(\bfk
)^{\times m} \to \Ccal_k(\bfk )$, defined as follows:
  $$
D_{\wfk}(M_1,\cdots ,M_m) = H^0\bigl( \Ufk_k(\pfk_{i_1})/\Ufk_k(\bfk
), M_1\otimes D_{\wfk'}(M_2,\cdots ,M_m)\bigr) ,
  $$
where $\wfk' := (s_{i_2}, \cdots, s_{i_m})$.  We similarly define
$\Dbar_{\wfk} : \Cbar_k(\bfk )^{\times m} \to \Cbar_k(\bfk )$.

These functors are again left exact.  Denote their right derived
functors respectively by $H^i\bigl(\Zcal^k_{\wfk}, M_1\boxtimes \cdots 
 \boxtimes
M_m\bigr )$ and $H^i\bigl(\Zbar^k_{\wfk}, \Mbar_1\boxtimes \cdots
\boxtimes \Mbar_m \bigr)$,
for $M_i\in \Ccal_k(\bfk )$ and $\Mbar_i\in\Cbar_k(\bfk )$.  These
are respectively $\Ufk_k(\bfk )$ and $\Ubar_k(\bfk)$-modules.  If
$M_1= \cdots = M_{m-1} = k$ is the trivial representation, then
\advance\eqcount by -4
  \begin{eqnarray}
D_{\wfk}(k,\cdots ,k,M_m) &\cong &D_{\wfk}(M_m), \qquad\hbox{and}  \\
H^i\bigl( \Zcal^k_{\wfk}, k\boxtimes \cdots \boxtimes k \boxtimes
M_m\bigr) &\cong& H^i(\Zcal^k_{\wfk}, M_m). \nonumber
  \end{eqnarray}
Analogous to the definition of $\Fr^*_{\wfk}$ (cf.\ Theorem 8.3), we can define the\break
$\Ufk (\bfk )$-module map  for any $\Mbar_i\in
\Cbar_k(\bfk )$,
  $$
\Fr^*_{\wfk} : H^i\bigl( \Zbar^k_{\wfk}, \Mbar_1\boxtimes \cdots
\boxtimes \Mbar_m\bigr)^{\Fr} \to H^i\bigl(\Zcal^k_{\wfk},
\Mbar_1^{\Fr}\boxtimes \cdots \boxtimes \Mbar_m^{\Fr}\bigr) .
  $$

Similarly, we can define the $\Ubar_k (\bfk )$-module map (cf.\ Theorem 8.4)
  $$
\Frps_{\wfk}: H^i\bigl( \Zcal^k_{\wfk}, \Mbar^{\Fr}_1\boxtimes \cdots
\boxtimes \Mbar_m^{\Fr}\bigr)^{\Frp_{\frak b}} \to H^i\bigl( \Zbar^k_{\wfk},
\Mbar_1\boxtimes \cdots \boxtimes \Mbar_m\bigr) .
  $$
Then
$$ \Fr^{\prime\ast}_{\wfk}\,\,\circ \Fr^*_{\wfk} = {\rm Id} \,\, \hbox{on} \,  
H^i\bigl( \Zbar^k_{\wfk}, \Mbar_1\boxtimes \cdots
\boxtimes \Mbar_m\bigr) .$$
  \enddemo

\phantom{almost time}
\vglue-36pt
 From now on, to the end of this section, assume that $\ell = p$ is 
a prime and $k$ is an algebraically closed field of characteristic $p$ which is
a $\zi$ algebra under $\xi \mapsto 1$.
\medbreak

Let $G,B$ be as in Section~6 and, for any $1\leq i \leq n$, let $B \subset P_i$ be the 
minimal parabolic subgroup containing the simple reflection $s_i$. Recall that
the Bott-Samelson-Demazure-Hansen variety $Z_\wfk$ is defined as
$P_{i_1} \times \cdots \times P_{i_m}/B^{\times m}$ where  $B^{\times m}$
acts on $P_{i_1} \times \cdots \times P_{i_m}$ from the right under 
$$(p_1, \cdots, p_m)\cdot (b_1, \cdots, b_m) := (p_1 b_1, b_1^{-1}p_2b_2, 
\cdots, b_{m-1}^{-1}p_mb_m),$$
for $p_j \in P_{i_j}$ and $b_j \in B$. 
Then  $Z_\wfk$ is a smooth projective variety over $k$. For any $\lam_1,
 \cdots, \lam_m \in X$, the character $e^{\lam_1}\boxtimes \cdots \boxtimes
e^{\lam_m}$ of $B^{\times m}$ gives rise to the line bundle $\Lcal_{\wfk}
({\lam_1}\boxtimes \cdots \boxtimes
{\lam_m})$ on  $Z_\wfk$. 
\smallbreak
Consider the embedding  $Z_\wfk \hookrightarrow (G/B)^{\times m}$ defined by
$$(p_1, \cdots, p_m) \,\hbox{mod} \,B^{\times m} \mapsto (p_1 B, p_1p_2 B,
 \cdots, p_1 \cdots p_m B).$$
Then the line bundle $\Lcal(\lam_1)\boxtimes \cdots \boxtimes \Lcal(\lam_m)$
on $G/B^{\times m}$ restricts to the line bundle  $\Lcal_{\wfk}
({\lam_1}\boxtimes \cdots \boxtimes
{\lam_m})$ on  $Z_\wfk$. In particular, if each of $\lam_1, \cdots, \lam_m$
is dominant and regular, then  $\Lcal_{\wfk}
({-\lam_1}\boxtimes \cdots \boxtimes
{-\lam_m})$ is ample on  $Z_\wfk$.

Define the $k$-algebra (under the multiplication of sections)
  $$
K_{\wfk} := \bigoplus_{(\lam_1,\cdots ,\lam_m)\in
(X^+)^{\times m} }\, H^0\bigl( Z_{\wfk}, \Lcal_{\wfk}(-\lam_1\boxtimes
\cdots \boxtimes -\lam_m)\bigr) . 
  $$
Analogous to the map $\theta$ of $\S$6.1, by induction on $\ell (\wfk
)$, using the Leray spectral sequence for the fibration $Z_{\wfk}\to
P_{i_1}/B$, we get
  $$
\theta_{\wfk}: H^0\bigl( \Zcal^k_{\wfk},
(\chi^{\xi}_{-\lam_1})_k\boxtimes \cdots \boxtimes
(\chi^{\xi}_{-\lam_m})_k\bigr)\simeq H^0\bigl( Z_{\wfk},
\Lcal_{\wfk}(-\lam_1 \boxtimes \cdots \boxtimes -\lam_m )\bigr)^{[1]}  .
  $$
Similarly,
  $$
\bar{\theta}_{\wfk}: H^0\bigl( \Zbar^k_{\wfk},
(\bar{\chi}_{-\lam_1})_k\boxtimes \cdots \boxtimes
(\bar{\chi}_{-\lam_m})_k\bigr) \simeq H^0\bigl( Z_{\wfk},
\Lcal_{\wfk}(-\lam_1 \boxtimes \cdots \boxtimes -\lam_m )\bigr).
  $$

Under the above identifications, the map $\Frps_{\wfk}$ gives rise to
the $k$-linear map 
  \begin{eqnarray*}
\Fcal'_{\wfk}(\lam_1, \cdots,  \lam_m): \Bigl( &H^0\bigl( Z_{\wfk},
\Lcal_{\wfk}(-p\lam_1 \boxtimes \cdots \boxtimes
-p\lam_m)\bigr)^{[1]}\Bigr)^{\Fr'_{\bfk}} \\
 &\to H^0\bigl( Z_{\wfk}, \Lcal_{\wfk}(-\lam_1\boxtimes \cdots
\boxtimes -\lam_m)\bigr) .
  \end{eqnarray*}
Combining these, we get the map (as maps of abelian groups;
without regard to the 
$\Ubar (\bfk )$ or $k$-linear structures)
$\Fcal'_{\wfk} : K_{\wfk}\to K_{\wfk}$, where we take
$\Fcal'_{\wfk |_{H^0(Z_{\wfk}, \Lcal_{\wfk}(-\lam_1\boxtimes \cdots
\boxtimes -\lam_m))}}\equiv 0$ unless $p$ divides each of $\lam_1, \cdots, 
\lam_m$.
\vglue3pt
By an argument similar to the proofs of Proposition~6.2, Theorem~6.4  and Proposition~7.3, we get the following.  It was first
proved by Mehta-Ramanathan [MR1] that $Z_{\wfk}$ is Frobenius split.

  \proclaim{Theorem}
The map $\Fcal'_{\wfk} : K_{\wfk}\to K_{\wfk}$ is a Frobenius\/{\rm -}\/linear
graded endomorphism{\rm .}  Moreover{\rm ,} it is a splitting{\rm .}  Further{\rm ,} for any
$q\in \Zgil_+$ and $f\in K_{\wfk}${\rm ,} 
  $$
\Fcal'_{\wfk} (\Ebar_i^{(pq)}\cdot f) = \Ebar_i^{(q)} \cdot
\Fcal'_{\wfk}(f),
  $$
where the action of $\Ebar_i^{(m)}$ comes from the canonical action
of $B$ on $$H^0\bigl( Z_{\wfk}, \Lcal_{\wfk}(-\lam_1\boxtimes\cdots
\boxtimes -\lam_m)\bigr).$$

In particular{\rm ,} the splitting $\Fcal'_{\wfk}$ is $B$\/{\rm -}\/canonical{\rm .}

The splitting $\Fcal'_{\wfk}$ induces a $B$\/{\rm -}\/canonical Frobenius splitting
of the variety $Z_{\wfk}$ by a method similar to that in Section~{\rm 6.}\endproclaim

  \numbereddemo{{R}emark}
1) For any reduced decomposition of the longest element of the Weyl group  $w_o=s_{i_1} \cdots 
  s_{i_N}$ , consider the sequence $\wfk_o=(s_{i_1},s_{i_2}, \cdots, 
 s_{i_N})$.  Then the above splitting of $Z_{\wfk_o}$ `descends' via
the map $Z_{\frak w_o} \to G/B$, $(p_1,\ldots , p_N)$ mod $B^{\ast N}\mapsto p_1\cdots p_N$ mod~$B$, to give the splitting of
$G/B$ given in Theorem~6.4.
\smallbreak
2) For any subsequence $\vfk$ of $\wfk$, the subvariety
$Z_{\vfk}\subset Z_{\wfk}$ is compatibly split under the splitting of 
$ Z_{\wfk}$ given in the above theorem.
  \enddemo

\vglue-12pt
  \section{Extension of results to the parabolic case}
\vglue-9pt

The aim of this section is to extend various results obtained in the
earlier sections for the Borel case to an arbitrary parabolic case.  We
formulate the extensions but omit the proofs as they are similar to
the proofs given earlier (of the corresponding results in the Borel
case).

Let $\ell$ be as in Sections~1--4 (i.e., it is an odd integer $>1$ assumed to
be coprime to 3 if $G_2$ is a factor of $\gfk$).

For a subset $I\subset \{ 1,\cdots ,n\}$, let $\Ubar (\pfk_I)$
be the parabolic subalgebra of $\Ubar (\gfk )$ generated by
$\Ubar (\bfk )$ and $\{\Fbar_i^{(m)};\; i\in I\hbox{ and }m\geq
0\}$.  Similarly, let $\Ufk (\pfk_I)$ be the parabolic
subalgebra of $\Ufk (\gfk )$ generated by $\Ufk (\bfk )$
and $\{ F_i^{(m)};\; i\in I  \hbox{ and }m\geq 0\}$.  Then $\Delta$
and $S$ keep $\Ufk (\pfk_I)$ stable.

For any subset $I$ as above and $\Ufk (\pfk_I)$-module $M$, we
can analogously define (cf.\ Definition 1.3) 
$$F_{\pfk_I}(M) :=\{
v\in F_{\bfk}(M): F_i^{(m)} v=0, \,\hbox{for all } m\geq m(v)\,\hbox{ and }
\, i\in I\}$$
and thus the category $\Ccal (\pfk_I)$.  Similarly, we can define
the category $\Cbar (\pfk_I)$.  Then Proposition~1.5  is true
with $\Ccal (\bfk)$ (resp.\ $\Ufk (\bfk )$) replaced by
$\Ccal (\pfk_I)$ (resp.\ $\Ufk (\pfk_I)$) and $\Cbar (\bfk)$ (resp.\ $\Ubar 
(\bfk)$) replaced by $\Cbar (\pfk_I)$ (resp.\ $\Ubar(\pfk_I)$).  Hence, we can
define the cohomology $H^i(\Ufk (\gfk )/\Ufk (\pfk_I), M)$,
for $M\in \Ccal (\pfk_I)$.  Similarly, we can define
$H^i(\Ubar (\gfk )/\Ubar (\pfk_I), \Mbar )$, for $\Mbar\in
\Cbar (\pfk_I)$.  We abbreviate $H^i(\Ufk (\gfk
)/\Ufk (\pfk_I), M)$ to $H^i(\Xfk_I, M)$ and similarly
$H^i(\Ubar (\gfk )/\Ubar (\pfk_I), \Mbar )$ to
$H^i(\Xbar_I, \Mbar )$.

For any $M\in \Ccal (\pfk_I)$, there is the canonical map $\pi_I
: H^i(\Xfk_I,M) \to H^i(\Xfk , M)$ and similarly for $\Mbar\in
\Cbar (\pfk_I)$ the map $\bar{\pi}_I: H^i(\Xbar_I, \Mbar)\to
H^i(\Xbar ,\Mbar )$.

Analogous to Theorems (2.3) and (3.8), we have the following:

 \proclaim{Theorem}  For any $\Mbar\in\Cbar (\pfk_I)${\rm ,} there exists a
 functorial $\Ufk (\gfk )$\/{\rm -}\/module map
  $$
\Fr_I^* : H^i (\Xbar_I , \bar{M})^{ \Fr}
\longrightarrow H^i\bigl(\Xfk_I, \bar{M}^{\Fr}\bigr)
  $$
compatible with $\pi_I$ 
in the sense that the following diagram is commutative\/{\rm :}
  $$  
    \begin{array}{ccc}
 H^i\bigl(\bar{X}_I, \bar{M}\bigr)^{ \Fr}&
\stackrel{\Fr^*_I}{\lrar}& H^i\bigl(\Xfk_I, \bar{M}^{\Fr}\bigr) \\[7pt]
 \Big\downarrow {\scriptstyle\bar{\pi}_I}&& \Big\downarrow{\scriptstyle \pi_I}  \\[10pt]
 H^i\bigl(\bar{X} , \bar{M}\bigr)^{ \Fr}
&\lower6pt\hbox{$\stackrel{\textstyle\lrar}{\scriptstyle \Fr^*}$}& H^i\bigl(\Xfk , \bar{M}^{\Fr}\bigr) .\\[-19pt]
     \end{array}  
   $$\endproclaim  

\proclaim{Theorem} For any $\Mbar\in \Cbar (\pfk_I)${\rm ,} there exists a
functorial $\Ubar (\bfk^- )$\/{\rm -}\/module map 
 $$
\Frps_I : H^i\bigl(\Xfk_I ,
\Mbar^{\Fr}\bigr)^{\Fr'_{\bfk^-}} \rightarrow
H^i\bigl(\Xbar_I , \Mbar\bigr) 
  $$ 
such that $\Frps_I \circ \Fr_{I}^*={\rm Id}${\rm .}

Further{\rm ,} 
$$ \bar{E_i}^{(m)} \cdot (\Frps_I f)= \Frps_{I} (E_i^{(m\ell)}\cdot f),$$
for any $f \in  H^i\bigl(\Xfk_I, \bar{M}^{\Fr}\bigr)${\rm .} 
\smallbreak
Moreover{\rm ,} $\Frps_I$ is compatible with $\pi_I${\rm .}\endproclaim

Let $\gamma_I := -2(\ell -1)\rho_I$, where $\rho_I := \sum_{i\notin
I}\omega_i$ and $\omega_i$ is the $i^{\rm th}$ fundamental weight defined by
$\omega_i(h_j)=\delta_{i,j}$.  Observe that the $\Ufk (\bfk )$-module
$\chi^{\xi}_{\gamma_I}$ is, in fact, a module for $\Ufk (\pfk_I)$.
Let $F_o^I := F^{(\ell -1)}_{\beta_{i_m}} \cdots F^{(\ell
-1)}_{\beta_{i_1}}$, where $i_1<\cdots < i_m$, $\{\beta_{i_1}, \cdots
,\beta_{i_m}\} =\Delta_+ \backslash \Delta_+(I)$, and $\Delta_+(I) :=
\Delta_+\cap \sum_{i\in I}\Zgil_+ \alpha_i$ (cf.\ $\S$4.1).  Then observe
that $F_o^I$ is of weight $\gamma_I$.  Decompose $V_{\xi}(-\gamma_I)  
= S(F_o^{I})\, \zi v_+\oplus M$, where $\zi v_+$
is the highest weight space and $M$ is a weight subspace of $V_\xi(-\gamma_I)$.
Let $\widehat{\sigma_o^I} \in V_{\xi}(-\gamma_I)^*$ be defined by
$\widehat{\sigma^I_o} (S(F^I_o)v_+)=1$ and $\widehat{\sigma^I_o}_{|M} \equiv
0$.  

Now, replacing $F_o$ by $F_o^I$ in Section~4, we get the following parabolic
analogue of Theorem~4.7  and Proposition~4.1.1.

  \proclaim{Theorem}  
For any $\Mbar\in\Cbar (\pfk_I )${\rm ,} there exists a functorial $\Ubar
(\bfk^-)$\/{\rm -}\/module map 
\vglue2pt
\centerline{${\displaystyle
\Frps_{\gamma_I} : H^i\bigl(\Xfk_I ,
\chi^{\xi}_{\gamma_I}
\otimes\Mbar^{\Fr}\big)^{\Fr'_{\bfk^-}} \rightarrow H^i(\Xbar_I , \Mbar ).
}$}
\vglue2pt

Further{\rm , }
$$ \bar{E_i}^{(m)} \cdot (\Frps_{\gamma_I} f)= \Frps_{\gamma_I}
(E_i^{(m\ell)}\cdot f),$$
for any $f \in  H^i\bigl(\Xfk_I, \chi^{\xi}_{\gamma_I}\otimes \bar{M}^{\Fr}\bigr)${\rm .} 

Moreover{\rm ,} the composite $\Frps_{\gamma_I} \circ \,
m_{\sigma_o^I}\circ\Fr_I^*$ is the identity map{\rm ,} where $\sigma^I_o\in
H^0\bigl(\Xfk_I, \chi^{\xi}_{\gamma_I}\bigr)$ is given as $\beta
(\widehat{\sigma^I_o})$ {\rm (}\/cf.\ Prop.~{\rm 4.11),} and $m_{\sigma_o^I}:
H^i\bigl(\Xfk_I, \Mbar^{\Fr}\bigr) \to H^i\bigl( \Xfk_I,
\chi^{\xi}_{\gamma_I}\otimes \Mbar^{\Fr}\bigr)$ is defined similarly
to Lemma~{\rm 4.9.}\endproclaim

Similar to Proposition~7.2, we have the following:

 \proclaim{Proposition}\label{deltasplitt} For any subsets $I, I'\subset
\{ 1,\cdots ,n\}${\rm ,} $\Mbar\in \Cbar (\pfk_I)$, $\Nbar\in\Cbar
(\pfk_{I'})$ and $i,j\in\Zgil_+${\rm ,} there exists a functorial
$\bar{U}(\bfk)$\/{\rm -}\/module map
  $$
\Frps_{\Delta}(I,I'):
H^i\bigl(\Xfk_I ,\Mbar^{\Fr}\otimes
H^j(\Xfk_{I'} ,\Nbar^{\Fr} )\bigr)^{\Fr'_{\bfk}}
\rightarrow
H^i\bigl(\bar{X}_I ,\bar{M}\otimes
H^j(\bar{X}_{I'} ,\bar{N})\bigl)
  $$
which is a splitting of the functorial $\Ufk(\gfk)$\/{\rm -}\/module  map
  $$
\Fr^*_{\Delta}(I,I'):
H^i\bigl(\bar{X}_I , \bar{M}\otimes
H^j(\bar{X}_{I'} , \bar{N})\bigl)^{\Fr}
\rightarrow
H^i\bigl(\Xfk_I ,\bar{M}^{\Fr}\otimes
H^j(\Xfk_{I'} ,\bar{N}^{\Fr})
\bigl) .
  $$
Moreover\/{\rm ,}\/ the analogue of {\rm (1)} of Proposition~{\rm 7.2} holds{\rm .} 
  \endproclaim

From now on, until the end of this section, we take $\ell = p$ to be a 
prime and $k$ an algebraically closed field of characteristic $p$, which is a $\Bbb Z_\xi$-algebra under $\Bbb Z_\xi\to k$,
$\xi\mapsto 1$.  Let $G,B$ and the Schubert varieties $X(w)\subset G/B$ be as in
Section~6.  For any subset $I\subset \{ 1,\cdots ,n\}$, let $B\subset
P=P_I$ be the parabolic subgroup containing the simple reflections
$\{ s_i\}_{i\in I}$.   (In particular, when  $I$ is the singleton $\{
i\}$, $P_I$ is the minimal parabolic subgroup $P_i$.)  Let $X_P$ be
the character group of $P$.  Then $X_P\cong \{ \lam\in X: \lam
(h_i)=0$, for all $i\in I\}$.  We set $X^+_P = X_P\cap
X^+$.
For any  $w\in W$, let $$X(w)_P := \overline{BwP/P} \subset G/P$$ be the
Schubert subvariety. Also define the opposite Schubert variety  $$X(w)_P^- := 
\overline{B^-wP/P} \subset G/P.$$  For any $\lam\in X_P$, the associated
homogeneous line bundle on $G/P$ is denoted by $\Lcal_P(\lam )$.
When there is no cause for confusion, we denote its restriction to
$X(w)_P$ again by $\Lcal_P(\lam )$.

  \numbereddemo{Definition}
Analogous to the isomorphisms $\theta$ and $\bar{\theta}$ of $\S$6.1,
we have the isomorphisms (for any $\lam\in X_P$)
\begin{eqnarray*}
\theta_I: H^0(\Ufk_k(\gfk )/\Ufk_k(\pfk_I ), (\chi^\xi_{\lam})_k)
&\simeq &H^0(G/P, 
\Lcal_P 
(\lam))^{[1]} \;\; \hbox{and} 
 \\ \bar{\theta}_I: H^0(\Ubar_k(\gfk )/\Ubar_k(\pfk_I ), 
(\bar{\chi}_\lam )_k) &\simeq &
H^0(G/P, \Lcal_P
(\lam)).
\end{eqnarray*}

Define the $X^+_P$-graded algebra
  $$
K_I := \bigoplus_{\mu\in X^+_P} H^0(G/P, \Lcal_P(-\mu )),
  $$
and let $\Fcal'_I : K_I\to K_I$ be the map defined by
$\Fcal'_{I|H^0(G/P, \Lcal_P(-\mu ))}\equiv 0$ if $\mu\notin pX^+_P$
and $\Fcal'_I: H^0(G/P, \Lcal_P(-p\mu )) \to H^0(G/P, \Lcal_P(-\mu
))$ is the splitting map $\Frps_I$ (as maps of abelian groups) under
the identifications $\theta_I$ and $\bar{\theta}_I$.
  \enddemo

Then analogous to Proposition~6.2, Theorems~6.4,  6.5  and 6.7,
we obtain the following. 

  \proclaim{Theorem}\label{diagonalsplitttheorem}
The map $\Fcal'_I$ is a Frobenius\/{\rm -}\/linear graded endomorphism   
of $K_I${\rm .}  Moreover{\rm ,} it is a $B$\/{\rm -}\/canonical splitting{\rm .}

This induces a $B$\/{\rm -}\/canonical Frobenius splitting of the flag
variety $G/P$ which compatibly splits any subscheme 
$Z_P\subset G/P$ obtained from $\{ X(w)_P,\break X(w)^-_P\}_{w\in
W}$ by repeatedly taking scheme theoretic unions{\rm ,}
intersections and irreducible components{\rm .}

In fact{\rm ,} $G/P$ is Frobenius $\Dcal_I$\/{\rm -}\/split{\rm ,} where
$\Dcal_I$ is the
line bundle $\Lcal_P(\gamma_I)$ together with the section $\phi_I :=
\theta_I(\sigma^I_0)${\rm .}\endproclaim

  \numbereddemo{{R}emark}
Since $G/P$ is Frobenius $\Lcal_P (\gamma_I)$-split, in particular, it
is Frobenius $\Lcal_P (-(p-1)\rho_I)$-split with respect to an
appropriate section of\break $H^0\bigl( G/P, \Lcal_P(-(p-1) \rho_I )\bigr)$.
Further, choosing the section appropriately, we can   easily show  that 
each Schubert variety $X(w)_P$ is
compatibly $\Lcal_P (-(p-1)\rho_I)$-split (under the splitting of $G/P$ obtained above).  This was originally proved
by Ramanathan (cf.\ [R2, Th.~3.5]).
  \enddemo

For any $I,I'\subset \{ 1,\cdots ,n\}$, define the $k$-algebra
  $$
K_{\Delta}(I,I') := \bigoplus_{\lam\in X_P, \,\mu\in X_{P'}} H^0(G/P,
\Lcal_P(-\lam )) \otimes H^0(G/P', \Lcal_{P'}(-\mu )), 
  $$
where $P := P_I$ and $P' := P_{I'}$.  

Analogous to Proposition~7.3, using the maps  \,$\Frps_{\Delta}(I,I')$
of $\S$9.4, we obtain the Frobenius-linear graded endomorphism
  $$
\Fcal'_{\Delta}(I,I') : K_{\Delta}(I,I') \to K_{\Delta}(I,I').
  $$
Moreover, it is a $B$-canonical splitting under the diagonal action of
$B$ on $K_{\Delta}(I,I')$.

Define the $G$-Schubert variety $Y(w)_{P,P'}$ as the image of
$Y(w)$ under the canonical projection map $G/B\times G/B \to G/P
\times G/P'$.  Analogous to Theorem~7.5, we obtain the following:

  \proclaim{Theorem}  The endomorphism $\Fcal'_{\Delta}(I,I')$ induces a
$B$\/{\rm -}\/canonical Frobenius splitting of $G/P\times G/P'${\rm .}

Further{\rm ,} any subscheme $Z\subset G/P\times G/P'$ obtained from the
$G$\/{\rm -}\/Schubert varieties $\{ Y(w)_{P,P'}\}_{w\in W}$ by repeatedly taking
unions{\rm ,} intersections and irreducible components is compatibly split{\rm .}\endproclaim

\centerline{\bf Appendix: Applications}
\vglue24pt

We follow the notation and assumptions of Section~9 (just above Definition 9.5).
In particular, $\ell =p$ is  a prime and $k$ is an algebraically closed field 
of characteristic $p$. 

For completeness, we collect some important (and standard)
consequences of Frobenius splitting of the flag varieties and their
Schubert subvarieties\break (cf.\ [MR1], [RR], [R1], [R2]; and also [A], 
[Jo], 
[S]).

  \nonumproclaim{Theorem A.1}  For any $w\in W$ and $\lam\in X^+_P${\rm ,} 
\begin{itemize}
\ritem{(a)} $H^i(X(w)_P, \Lcal_P(-\lam ))=0${\rm ,} for all $i>0${\rm .}

\ritem{(b)} The restriction map $H^0(G/P, \Lcal_P(-\lam )) \to H^0(X(w)_P,
\Lcal_P(-\lam ))$ is surjective{\rm .}
\end{itemize}

\endproclaim

 {\it Proof}.  Since $G/P$ is Frobenius $\Lcal_P(-
(p-1)\rho_I)$-split  and $X(w)_P$ is compatibly
$\Lcal_P(-(p-1)\rho_I)$-split  by Remark~9.7, the theorem follows immediately from
the standard properties of Frobenius splitting (cf.\ [R2, Prop.\
1.13(ii)]).   \phantom{hithereyou}
  \hfill\qed\medbreak

See [RR, Th.\ 3] for the following result. 

  \nonumproclaim{Theorem A.2}   Any Schubert variety
$X(w)_P\subset G/P$ is normal{\rm .}

Moreover{\rm ,} for any homogeneous ample line bundle $\Lcal =\Lcal_P(-\lam
)$ on $G/P${\rm ,} $X(w)_P$ is projectively normal in the projective
embedding given by $\Lcal${\rm .}\endproclaim

  \demo{Proof}  To prove the normality, we can of course assume
that $P=B$.  We prove the normality of $X(w)$ by induction on the length $\ell
(w)$.  If $\ell (w)=0$, there is nothing to prove.  So take $\ell
(w)>0$ and write $w=w's_i$ for a simple reflection $s_i$ such that
$w'<w$.  Under the canonical map $\pi : G/B\to G/P_i$ , $X(w)$ and
$X(w')$ have the same image $X(w)_{P_i}$.  Moreover, $\pi_{|_{X(w)}}:
X(w)\to X(w)_{P_i}$ is a $\Pgil^1$-fibration and $\pi_{|_{X(w')}}:
X(w')\to X(w)_{P_i}$ is a birational map.  By induction, $X(w')$ is
normal.  For any $\lam\in X^+_{P_i}$, we have the commutative diagram:
$$  \begin{array}{ccc}
H^0\bigl( G/P_i, \Lcal_{P_i}(-\lam )\bigr) & \twoheadrightarrow
&  H^0\bigl( X(w)_{P_i},    \Lcal_{P_i}(-\lam )\bigr) \\
 \quad\quad\downarrow\!{\wr} & &  \downarrow\! \pi^*\hskip.55in \\[6pt]
H^0\bigl( G/B , \Lcal (-\lam )\bigr) & \twoheadrightarrow
&   H^0\bigl( X(w'),  \Lcal(-\lam )\bigr),\quad
  \end{array}
$$ 
where the vertical maps are induced by $\pi$ and the horizontal maps
are induced by inclusions.  The horizontal maps are surjective by
Theorem~A.1  and clearly the left vertical map is an isomorphism,
hence the right vertical map $\pi^*$ is surjective.  
Since $\pi^*$ is surjective for all $\lam\in X^+_{P_i}$ (in
particular, for all large enough positive powers of an ample line
bundle $\Lcal$ on $X(w)_{P_i}$), and $H^1(X(w)_{P_i},
\Lcal_{P_i}(-\lam ))=0$ (by Theorem A.1), we get $(\pi_{|_{X(w')}})_*
\Ocal_{X(w')}=\Ocal_{X(w)_{P_i}}$.  But since $X(w')$ is normal, so
is $X(w)_{P_i}$ and hence $X(w)$ is normal.
\smallbreak
Now we come to the projective normality:  It suffices to show that
the multiplication map
  $$
  H^0\bigl( X(w)_P, \Lcal^{\otimes m}\bigr) \otimes H^0\bigl( X(w)_P,
\Lcal^{\otimes n}\bigr) \longrightarrow H^0\bigl( X(w)_P,
\Lcal^{\otimes (n+m)}\bigr) \tag{1}
  $$
is surjective for all $m,n \geq 1$ (cf.\ [H, Chap. II, Ex.\
5.14(d)]).

By the compatible Frobenius splitting of the diagonal $G/P
\hookrightarrow G/P\times G/P$ (cf.\ Th.\ 9.8), we get that
  $$
H^0\bigl( G/P, \Lcal^{\otimes m}\bigr) \otimes H^0\bigl( G/P,
\Lcal^{\otimes n}\bigr) \twoheadrightarrow H^0\bigl( G/P,
\Lcal^{\otimes (n+m)}\bigr) \tag{2}
  $$
is surjective.  Now (1) follows from (2) by Theorem A.1(b).
  \enddemo

For any sequence of simple reflections $\wfk = (s_{i_1},\cdots
,s_{i_m})$, consider the $B$-equivariant morphism  $\psi_{\wfk} : 
Z_{\wfk} \to G/B$ given by $\psi_{\wfk}\bigl(
(p_1,\cdots ,p_m) \hbox{ mod } B^{\times m}\bigr)$ $= p_1\cdots p_mB$.  For
any $\lam\in X$, let $\Lcal_{\wfk}(\lam )$ be the pull-back line
bundle $\psi^*_{\wfk}(\Lcal (\lam ))$ on $Z_{\wfk}$.  
For $\wfk$ as above, let $\theta (\wfk ): = s_{i_1}\cdots s_{i_m} \in
W$.  A sequence $\wfk$ is called {\it reduced} if the above
decomposition of $\theta (\wfk )$ is reduced.  If $\wfk$ is reduced,
$\hbox{Image}\; \psi_{\wfk}=X(\theta (\wfk ))$ and $\psi_{\wfk} :
Z_{\wfk} \to X(\theta (\wfk ))$ is birational.  In the notation of
$\S$8.6, $\Lcal_{\wfk}(\lam ) \simeq \Lcal_{\wfk}(0\boxtimes \cdots
\boxtimes 0\boxtimes \lam )$.

  \nonumproclaim{Lemma A.3}  For any $\lam\in X$ and any reduced $\wfk${\rm ,}
  \begin{eqnarray}
H^0(Z_{\wfk}, \Lcal_{\wfk}(-\lam ))& \simeq& H^0(X(\theta (\wfk
)), \Lcal (-\lam )).   \\[3pt]
\noalign{\noindent Moreover{\rm ,} for any $\lam\in X^+${\rm ,}}
\noalign{\vskip4pt}
H^i(Z_{\wfk}, \Lcal_{\wfk}(-\lam ))&=&0, \quad \;\hbox{ for all }i>0  
  \end{eqnarray}
In particular{\rm ,} by Kempf\/{\rm '}\/s lemma {\rm (}cf.{\rm ,} e.g.{\rm , [D, $\S$5,} Prop.\
{\rm 2]),} $R^i\psi_{\wfk *}(\Ocal_{Z_{\wfk}})=0$ for all $i>0${\rm .}
  \endproclaim 

  \demo{Proof}  Since $X(\theta (\wfk ))$ is normal and
$\psi_{\wfk}: Z_{\wfk} \to X(\theta (\wfk ))$ is birational, (1)
follows.  

Consider the fibration $\eta_{\wfk}: Z_{\wfk}\to P_{i_1}/B$ with fibre
$Z_{\wfk'}$, where $\wfk' := (s_{i_2},\cdots ,s_{i_m})$.  Assume, by
induction on $m$, that $H^i(Z_{\wfk'}, \Lcal_{\wfk'}(-\lam ))=0$ for
all $i>0$.  Hence, by the degenerate Leray spectral sequence, we get:
  $$
H^i(Z_{\wfk}, \Lcal_{\wfk}(-\lam )) \simeq H^i\bigl( P_{i_1}/B, \Lcal
(H^0(Z_{\wfk'}, \Lcal_{\wfk'}(-\lam )))\bigr), \tag{3}
  $$
where, for a $B$-module $M$, $\Lcal (M)$ denotes the associated
homogeneous vector bundle.

Now, by Theorem A.1(b), we have the surjective map
  $$
H^0(G/B, \Lcal (-\lam )) \twoheadrightarrow H^0(X(w'), \Lcal (-\lam
)) \simeq H^0(Z_{\wfk'}, \Lcal_{\wfk'}(-\lam )) \tag{4}
  $$
with kernel, say $K$, where $w' := \theta (\wfk')$.
From the long exact cohomology sequence associated to (4), we get
  $$
H^1\bigl( P_{i_1}/B, \Lcal (H^0(G/B, \Lcal (-\lam )))\bigr)
\twoheadrightarrow H^1\bigl( P_{i_1}/B, \Lcal (H^0(Z_{\wfk'},
\Lcal_{\wfk'}(-\lam )))\bigr) , \tag{5}
  $$
since $H^2(P_{i_1}/B, \Lcal (K))=0$ from the dimensional consideration.  But since  $H^0(G/B, \Lcal (-\lam ))$
is a $G$-module, $H^1\bigl( P_{i_1}/B$,  $\Lcal (H^0(G/B$, $\Lcal
(-\lam )) )\bigr) =0.$   So, by (3) and (5), we get $$H^i(Z_{\wfk},
\Lcal_{\wfk}(-\lam ))=0$$ for all $i>0$.  This proves (2).
  \enddemo
    
  \nonumproclaim{Theorem A.4}  
Any Schubert variety $X(w)_P\subset G/P$ is Cohen\/{\rm -}\/Macaulay{\rm .}

Moreover{\rm ,} for any homogeneous ample line bundle $\Lcal =\Lcal_P(-\lam
)$ on $G/P${\rm ,} $X(w)_P$ is projectively Cohen\/{\rm -}\/Macaulay in the
projective embedding given by $\Lcal${\rm .}
 \endproclaim

 {\it Proof}.
We first prove that $X(w)_P$ is Cohen-Macaulay.  We can clearly
assume that $P=B$.  By the standard characterization of Cohen-Macaulay
schemes (cf.\ [H, Chap. III, Th.\ 7.6 and its proof]), it suffices to
show that
  $$
H^i(X(w), \Lcal (\lam ))=0, \tag{1}
  $$
for all $i<\ell (w)$ and all dominant regular $\lam$.

Take a reduced sequence $\wfk$ with $\theta (\wfk )=w$.  Then, by Lemma~A.3,
  $$
H^i(X(w), \Lcal (\lam )) \simeq H^i(Z_{\wfk}, \Lcal_{\wfk}(\lam
)). \tag{2}
  $$

Assume, by induction, that $H^i(Z_{\wfk'}, \Lcal_{\wfk'}(\lam ))=0$ for
all $i<m-1$, where $m:=\ell (w)$.  Then, by the Leray spectral sequence
for the fibration $Z_{\wfk }\to P_{i_1}/B$, 
  $$
H^i(Z_{\wfk}, \Lcal_{\wfk}(\lam ))=0 \,\,\,\hbox{unless} \,\, i=m-1\,\hbox{ or}
\,\, m.\tag{3}
  $$
  Now, by   Serre duality,
  $$
H^{m-1}(Z_{\wfk}, \Lcal_{\wfk}(\lam )) \simeq H^1(Z_{\wfk}, K_{Z_{\wfk}}
\otimes \Lcal_{\wfk}(-\lam ))^* \,,\tag{4}
  $$
where $K_{Z_{\wfk}}$ is the canonical bundle of $Z_{\wfk}$.  By [R1,
Prop.~2], 
  $$
K_{Z_{\wfk}} \simeq \Ocal_{Z_{\wfk}}[-\partial Z_{\wfk}] \otimes
\Lcal_{\wfk}(\rho ),
  $$
where $\rho$ is the half sum of positive roots,
$\partial Z_{\wfk} := \psi^{-1}_{\wfk}(\partial X(w))$, and $\partial X(w) :=
X(w)\backslash (BwB/B)$.

\advance\eqcount by 2
From the sheaf exact sequence:
  $$
0 \to \Ocal_{Z_{\wfk}}[-\partial Z_{\wfk}]\to \Ocal_{Z_{\wfk}} \to
  \Ocal_{\partial Z_{\wfk}} \to 0
  $$
tensored with $\Lcal_{\wfk}(\rho -\lam ),$ we get the exact sequence:
  \begin{eqnarray}
&&H^0(Z_{\wfk}, \Lcal_{\wfk}(\rho -\lam )) \stackrel{r} \to 
 H^0\bigl(\partial
  Z_{\wfk}, \Lcal_{\wfk}(\rho -\lam )_{|\partial Z_{\wfk}}\bigr) 
 \\
&&\qquad \to H^1\bigl( Z_{\wfk}, \Ocal_{Z_{\wfk}}[-\partial Z_{\wfk}]\otimes
  \Lcal_{\wfk}(\rho -\lam )\bigr) \to  H^1(Z_{\wfk},
  \Lcal_{\wfk}(\rho -\lam ))=0 .\nonumber
  \end{eqnarray}
(By Lemma A.3, the last term is 0.)  We now prove that the restriction 
  map $r$ is surjective:

By the following lemma (and Theorem 8.7, Remark 8.8(2) and Theorem 6.7),  the map $\psi'_{\wfk} : \partial Z_{\wfk}\to \partial
   X(w)$, gotten by restricting $\psi_{\wfk}$, satisfies
  $$
( \psi'_{\wfk})_\ast \Ocal_{\partial Z_{\wfk}} = \Ocal_{\partial X(w)}. \tag{6}
  $$
(Observe that $\psi'_{\wfk}$ has connected fibres since $\psi_{\wfk}$
   has connected fibres by Zariski's main theorem, as $X(w)$ is
   normal by Theorem A.2.)

From (6) we get (for any $\mu \in X$)
  $$
H^0\bigl( \partial Z_{\wfk}, \Lcal_{\wfk}(\mu )_{|\partial Z_{\wfk}}\bigr)
  \simeq H^0\bigl( \partial X(w), \Lcal (\mu )_{|\partial X(w)}\bigr) .
  $$
Now the surjectivity of $r$ follows from (1) of Lemma~A.3, since 
 $G/B$ is Frobenius $\Lcal (-(p-1)\rho)$-split and $\partial X(w)$
  is compatibly  $\Lcal (-(p-1)\rho)$-split (cf.\ Remark 9.7).
  Thus, from
  the exact sequence (5), we get
  $$
H^1\bigl( Z_{\wfk}, \Ocal_{Z_{\wfk}}[-\partial Z_{\wfk}] \otimes
  \Lcal_{\wfk} (\rho -\lam )\bigr) =0.
  $$
So, from (3) and (4), we get
  $$
H^i(Z_{\wfk}, \Lcal_{\wfk}(\lam ))=0 \;\hbox{ unless }i=m.
  $$
By  (2), this proves (1) and hence  $X(w)$ is Cohen-Macaulay.

To prove that $X(w)_P$ is projectively Cohen-Macaulay, it suffices to
   show  (in view of Theorem A.2) that
  $$
H^i(X(w)_P, \Lcal^n)=0 \;\hbox{ for all $0<i<\dim X(w)_P$ and all
   $n\in\Zgil$}.  \tag{7}
  $$
Since $X(w)_P$ is Cohen-Macaulay, we get (7) for all $0\leq i <\dim
   X(w)_P$ and all $n<0$ (by [H, Chap. III, Th.\ 7.6(b)]).  
The vanishing (7) for $n\geq 0$
   follows from Theorem A.1.  (a)  This proves the theorem.
  \hfill\qed\medbreak

We recall the following simple lemma due to Mehta-Srinivas [MS, Lemma~2].

\nonumproclaim{Lemma A.5} Let $\pi : X \to Y $ be a proper and surjective morphism between 
Frobenius split schemes such that all its fibres are connected{\rm .} Assume further
 that for each irreducible component $Y'$ of $Y$ there exists a component 
$X'$ of $X$  such that  $\pi_{\vert X'} : X' \to Y' $
is birational{\rm .} Then $\pi_*\Ocal_X = \Ocal_Y${\rm .}
\endproclaim

\demo{{R}emark {\rm A.6}} The above proof (of Theorem A.4) is a minor simplification 
of the proof given in [MS].
\enddemo

\nonumproclaim{{C}orollary A.7 {\rm (of Theorem A.4)}} For any  reduced $\wfk${\rm ,} the resolution 
$\psi_\wfk : Z_\wfk \to X(w)$ is rational{\rm ,} where $w :=\theta 
(\wfk)${\rm .} 
\endproclaim
 
\demo{Proof} In view of Lemma~A.3, it suffices to show that 
 $$
 R^q\psi_{\wfk *}(K_{Z_{\wfk}})=0\;\hbox{ for all $0< q$}.  \tag{1}
  $$
Fix an ample line bundle $\Lcal (-\lam)$ on $G/B$. To prove (1), it suffices 
to show that 
$$
 H^0\bigl(X(w), R^q\psi_{\wfk *}(K_{Z_{\wfk}} \otimes 
\Lcal_\wfk(-n\lam))\bigr)=0\;\hbox{ for all $n \gg 0$}.  \tag{2}
 $$
We choose $n_o$ large enough so that for all $n \geq n_o$ and $q\geq 0$
 $$
 H^p\bigl(X(w), R^q\psi_{\wfk *}(K_{Z_{\wfk}} \otimes 
\Lcal_\wfk(-n\lam))\bigr)=0\;\hbox{ for all $p >  0$},  
 $$
(cf.\ [H, Chap. III, Prop.\ 5.3]). Then, by the degenerate Leray spectral
sequence for $\psi_\wfk$, we get (for all $n \geq n_o$)
$$
 H^q\bigl(Z_\wfk , K_{Z_{\wfk}} \otimes 
\Lcal_\wfk(-n\lam))\simeq
 H^0\bigl(X(w), R^q\psi_{\wfk *}(K_{Z_{\wfk}} \otimes 
\Lcal_\wfk(-n\lam))\bigr).  \tag{3}
 $$
By the Serre duality,
\begin{eqnarray*}
 H^q\bigl(Z_\wfk , K_{Z_{\wfk}} \otimes 
\Lcal_\wfk(-n\lam)\bigr)&\simeq &
 H^{\ell(w)-q}\bigl(Z_\wfk ,
\Lcal_\wfk(n\lam)\bigr)^*\\
&\simeq &
 H^{\ell(w)-q}\bigl(X(w) ,
\Lcal (n\lam)\bigr)^*\,, \,\,\hbox{by Lemma~A.3}\\
&=& 0\,, \,\,\hbox{by (1) of Theorem~A.4}.
 \end{eqnarray*}
Combining this with (3), we get (2), proving the corollary.
\enddemo

\demo{{R}emark {\rm A.8(a)}} It is possible that the restriction on $\ell$ in this paper 
(i.e.  $\ell > 1 $ is an odd integer and coprime to $3$ if $G_2$ is a factor of $\gfk$) can be removed by using the results of Kaneda [Kan] and Andersen-Paradowski [AP] (also see [Li] and [KL, Remark 2]).

(b) It is natural to ask if the results of this paper can be extended to the symmetrizable Kac-Moody algebras.
\enddemo

\AuthorRefNames [APW]


\begin{references}

\bibitem{A} 
  \name{H.\ H.\ Andersen},
Schubert varieties and Demazure's character formula,  {\it Invent.\ Math.}\
{\bf 79} (1985),  611--618.

\bibitem{AP}
 \name{H.\ H.\ Andersen} and \name{J.\ Paradowski},
Fusion categories arising from semisimple Lie algebras,  
{\it Comm.\ Math.\ Phys.}\/
{\bf 169} (1995), 563--588.
 
\bibitem{APW} 
  \name{H.\ H.\ Andersen, P.\  Polo}, and \name{K. Wen},
Representations of quantum algebras,  {\it Invent.\ Math.}\/ 
{\bf 104} (1991),  1--59.

\bibitem{AW} 
  \name{H.\ H.\ Andersen} and \name{K. Wen},
Representations of quantum algebras. The mixed case,  {\it J. Reine Angew.\
 Math.}\/
{\bf 427} (1992),  35--50.


\bibitem{CPS}
   \name{E.\ Cline, B.\ Parshall}, and \name{L.\ Scott},  Cohomology,
hyperalgebras, and  representations,  
{\it J.\ Algebra} {\bf 63}  (1980), 98--123.

\bibitem{D}
\name{M.\ Demazure}, D\'{e}singularisation des vari\'{e}t\'{e}s de Schubert 
g\'{e}n\'{e}ralis\'{e}es,
{\it Ann.\ Sci.\ \'{E}cole Norm.\ Sup\/}.\  {\bf 7} (1974), 53--88.

\bibitem{H}
   \name{R.\ Hartshorne}, {\it Algebraic Geometry}, {\it Grad.\ Texts
   in Math\/}.\ {\bf 52}, 
Springer-Verlag, New York  (1977).

\bibitem{J}
   \name{J.\ C.\ Jantzen}, {\it  Representations of Algebraic
   Groups}, {\it Pure and Applied Math\/}.\ {\bf 131}, 
Academic Press, Inc., Boston, MA (1987).  

\bibitem{Jo}
 \name{A.\ Joseph},
On the Demazure character formula,  
{\it Ann.\ Sci.\ \'{E}cole Norm.\ Sup\/}.\  {\bf 18} (1985), 389--419.

\bibitem{Ka}
\name{W.\ van der Kallen}, {\it Lectures on Frobenius Splittings and
$B$\/{\rm -}\/modules}, published for
{\it T.I.F.R.{\rm ,}  Bombay}, Springer-Verlag, New York (1993).

\bibitem{Kan}
\name{M.\ Kaneda}, 
Cohomology of infinitesimal quantum algebras,
{\it J.\ Algebra} {\bf 226} (2000), 250--282.
 
\bibitem{Kas}
\name{M.\ Kashiwara}, 
The crystal base and Littelmann's refined Demazure character formula, 
{\it Duke Math.\ J. } {\bf 71} (1993),  839--858.

\bibitem{K}
\name{G.\ Kempf}, 
Linear systems on homogeneous spaces, 
{\it Annals of Math\/}.\ {\bf 103} (1976), 557--591.

\bibitem{Ku}
\name{S.\ Kumar}, 
Representations of quantum groups at roots of unity, in
{\it Quantum Topology} (D.\ N.\ Yetter, ed.),  World Scientific, Singapore,
  187--224 (1994).


\bibitem{KL}
 \name{S.\ Kumar} and \name{P.\ Littelmann}, Frobenius splitting
in characteristic
zero and the quantum Frobenius map,  {\it J. of Pure and Applied Algebra}
{\bf 152} (2000),  201--216.

\bibitem{Li}
   \name{P.\ Littelmann}, Contracting modules and standard monomial
   theory for symmetrizable Kac-Moody algebras,
{\it JAMS} {\bf 11}  (1998),  551--567.

\bibitem{Lu1}
\name{G.\ Lusztig}, Modular representations and quantum groups,
in {\it   Classical Groups and Related Topics\/},   
{\it Contemp.\ Math.}\/ {\bf 82} (1989), 59--77.

\bibitem{Lu2}
  \bibline, Quantum groups at roots of $1$,
{\it Geom.\ Dedicata} {\bf 35} (1990),  89--113.

\bibitem{Lu3}
  \bibline, {\it Introduction to Quantum Groups},   {\it Progress 
in Math\/}.\ {\bf 110},  Birkh{\rm \"{\it a}}user,
Boston (1993).


\bibitem{M}
\name{O.\ Mathieu}, Filtrations of $G$-modules,
{\it Ann.\ Sci.\ {\rm \'{\it E}}cole Norm.\ Sup\/}.\ {\bf 23} (1990),  625--644.


  \bibitem{MR1}
\name{V.\ Mehta} and \name{A.\ Ramanathan}, Frobenius splitting and cohomology 
vanishing for Schubert varieties,
{\it Annals of Math}.\/ {\bf 122} (1985), 27--40.

\bibitem{MR2}
\bibline, Schubert varieties in $G/B \times
G/B$,
{\it Compositio Math\/}.\ {\bf 67} (1988), 355--358.


  \bibitem{MS}
\name{V.\ Mehta} and \name{V.\ Srinivas}, A note on 
 Schubert varieties in $G/B$,
{\it Math.\ Ann}.\/ {\bf 284} (1989),  1--5.

\bibitem{Mi}
 \name{J.\ Milnor}, {\it Introduction to Algebraic $K$-theory},
 {\it Ann.\ of Math.\ Studies\/} {\bf 72}, 
  Princeton Univ.\  Press, Princeton, NJ 
 (1971).

\bibitem{RR}
\name{S.\ Ramanan} and \name{A.\ Ramanathan}, 
Projective normality of flag varieties 
and Schubert varieties, 
{\it Invent.\ Math\/}.\ {\bf 79} (1985), 217--224.

\bibitem{R1}
\name{A.\ Ramanathan},   Schubert varieties are arithmetically 
Cohen-Macaulay, 
{\it Invent.\ Math}.\/ {\bf 80} (1985), 283--294.

\bibitem{R2}
\name{A.\ Ramanathan},  Equations defining  Schubert varieties and Frobenius
splitting of diagonals, 
{\it IHES Publ.\ Math\/}.\  {\bf 65} (1987),  61--90.

\bibitem{S}
\name{C.\ S. Seshadri}, Line bundles on  Schubert varieties, in 
{\it Vector Bundles on Algebraic Varieties\/} (T.I.F.R Colloquium, Bombay, 1984), 
Oxford University Press
(1987), 499--528.

  \end{references}
   \end{document}